\documentclass[leqno,11pt]{amsart}
\usepackage{amsmath,amscd}
\usepackage{epsfig,graphics}
\usepackage{amssymb,latexsym}
\usepackage{mathrsfs}
\usepackage{eucal}
\usepackage{hyperref}

\newtheorem{thm}{\bf Theorem}[section]
\newtheorem{df}[thm]{\bf Definition}
\newtheorem{prop}[thm]{\bf Proposition}
\newtheorem{cor}[thm]{\bf Corollary}
\newtheorem{lem}[thm]{\bf Lemma}
\newtheorem{rem}[thm]{\bf Remark}
\newtheorem{ex}[thm]{\bf Example}

\newcommand{\bs}{\boldsymbol}
\newcommand{\A}{\mathcal{A}}

\newcommand{\B}{\mathbf{B}}
\newcommand{\cB}{\mathcal{B}}
\newcommand{\W}{\mathcal{W}}
\newcommand{\cP}{\mathscr{P}}

\newcommand{\pf}{\noindent{\bfseries Proof. }}
\newcommand{\ov}{\overline}
\newcommand{\w}{{\bf w}}
\newcommand{\bi}{\bs{\rm i}}
\newcommand{\bj}{\bs{\rm j}}

\newcommand{\bk}{\bs{\rm k}}
\newcommand{\bl}{\bs{\rm l}}

\newcommand{\M}{{\mathcal{M}}}

\newcommand{\hf}{\frac{1}{2}}

\newcommand{\gl}{\mathfrak{gl}}
\newcommand{\Z}{\mathbb{Z}}
\newcommand{\C}{\mathbb{C}}
\newcommand{\h}{\mathfrak{h}}
\newcommand{\te}{\widetilde{e}}
\newcommand{\tf}{\widetilde{f}}

\newcommand{\td}{\widetilde}

\newcommand{\mf}{\mathfrak}
\newcommand{\La}{\Lambda}
\newcommand{\I}{\mathcal{I}}

\numberwithin{equation}{section}

\begin{document}
\title[ ]
{Littlewood identity and Crystal bases}
\author{JAE-HOON KWON}
\address{Department of Mathematics \\ University of Seoul   \\  Seoul 130-743, Korea }
\email{jhkwon@uos.ac.kr }

\thanks{This work was partially supported by the National Research Foundation (NRF) grant funded by the Korean government (MEST) (No. 2011-0006735) and by LG Yonam foundation.}

\begin{abstract}
We give a new combinatorial model for the crystals of integrable highest weight modules over the classical Lie algebras of type $B$ and $C$ in terms of classical Young tableux. We then obtain a  new description of  its Littlewood-Richardson rule and a maximal Levi branching rule in terms of  classical Littlewood-Richardson tableaux, which extends in a bijective way the well-known stable formulas at large ranks.
We also show that  this tableau model admits a natural superization and it produces the  characters of  irreducible highest weight modules over orthosymplectic Lie superalgebras, which  correspond to the integrable highest weight modules over the classical Lie algebras of type $B$ and $C$ under the Cheng-Lam-Wang's super duality. 
\end{abstract}

\maketitle


\section{Introduction}
\subsection{Motivation}
Let $\frak{g}$ be a symmetrizable Kac-Moody algebra, and let $V_\lambda$ be an irreducible highest weight $\frak{g}$-module with highest weight $\lambda$. For dominant integral weights $\lambda,\mu,\nu$, the multiplicity of $V_\lambda$ in the tensor product $V_\mu\otimes V_\nu$ is called a {\it generalized LR (Littlewood-Richardson) coefficient}.

When $\frak{g}=\gl_n$, the multiplicity  $c^\lambda_{\mu\nu}$ of $V_\lambda$ in $V_\mu\otimes V_\nu$ is called a (classical) LR coefficient and it is given by the number of certain semistandard Young tableaux of skew shape called  LR tableaux (see \cite{Fu,Mac95}). When $\frak{g}$ is the classical Lie algebra of  type  $B,C,D$, the generalized LR coefficients and various branching coefficients can also be expressed in terms of classical LR coefficients. There are nice formulas, where the multipicities are given as a  sum of products of  $c^\lambda_{\mu\nu}$'s when the highest weights of irreducible representations are in a certain stable range. The formulas of this type originated in the Littlewood's restriction formula for the representations of orthogonal and symplectic groups to general linear groups \cite{Lw-1,Lw-2} (we refer the reader to \cite{HTW,Su''} which have nice survey on the related works).  One may understand these formulas as stable limits of the generalized LR  and branching coefficients at large ranks of $\frak{g}$, where we fix partitions parametrizing highest weights and take the rank of $\frak{g}$ to infinity. But the multiplicities outside a stable range are given as an alternating sum of products of  $c^\lambda_{\mu\nu}$'s, which may contain many cancellations \cite{EW,Kng-b,KoT90}, and it is not known  in general  how to extend these stable limits in a subtraction-free way for arbitrary dominant integral highest weights (see \cite{Su'} for a combinatorial extension of the Littlewood's restriction formula in the symplectic case for arbitrary partitions).

One of  our main motivations in this work is to find a new combinatorial model, which explains the stable limits of generalized LR and branching coefficients in type  $B,C,D$, and  also extends them to the case of arbitrary highest weights in a bijective or subtraction-free way. We will approach this problem using the Kashiwara's crystal base theory \cite{Kas90}, which provides a powerful combinatorial tool for computing the LR rule for the Drinfeld and Jimbo's quantized enveloping algebra of ${\mf g}$ (cf.\cite{BZ,Li-2,Li,Na}).

\subsection{Littlewood identity and crystals} Let $\cP$ be the set of partitions and $s_\lambda({\bf x})$  the Schur function in ${\bf x}=\{x_1,x_2,\ldots \}$. We start with considering the following {\it Littlewood identity} 
\begin{equation}\label{Littlewood identity}
\dfrac{1}{\prod_i(1-x_i^\epsilon)\prod_{i<j}(1-x_ix_j)}=\sum_{\tau\in \cP}s_{\epsilon\tau}({\bf x}),
\end{equation}
where $\epsilon=1,2$ \cite{Lw-2}.

Let ${\mf b}_\infty$ and ${\mf c}_\infty$ be the affine Lie algebras of type $B_\infty$ and  $C_\infty$, respectively, and let ${\mf l}_\infty$ be a Levi subalgebra of type $A_{+\infty}$  obtained by removing a short  (resp. long) simple root $\alpha_0$. Note that the identity (\ref{Littlewood identity}) multiplied by $s_\lambda({\bf x})$ is the character of a generalized Verma module over $U_q({\mf x}_\infty)$ induced from a highest weight ${\mf l}_\infty$-module with highest weight corresponding to $\lambda$, and it is also the character of the set of pairs of semistandard Young tableaux $(S,T)$ of shape $\epsilon\tau$ and $\lambda$, respectively. Here  we assume $\epsilon=1$ for ${\mf x}={\mf b}$ and $\epsilon=2$ for ${\mf x}={\mf c}$.

Let $\B({\mf x}_\infty,\Lambda)$ denote the crystal associated to the integrable highest weight $U_q({\mf x}_\infty)$-module with highest weight $\Lambda$. We parametrize the dominant integral weights for ${\mf x}_\infty$  by
$(\lambda,n)\in \cP\times\Z_{>0}$ such that $2\lambda_1\leq \epsilon\,n$ (see Section \ref{affine Lie algebra}), and denote by $\Lambda^{\mf x}(\lambda,n)$ the corresponding dominant integral weight.
Motivated by the above observation on the Littlewood identity multiplied by a Schur function, we consider a crystal structure on the set of bitableaux associated with a generalized Verma module with ${\mf l}_\infty$-highest weight $\lambda$, and then characterize $\B({\mf x}_\infty,\La^{\mf x}(\lambda,n))$ as its subcrystal.
More precisely, we show as the main result in this paper (Theorem \ref{main result}) that $\B({\mf x}_\infty,\La^{\mf x}(\lambda,n))$ can be realized as the set of pairs of semistandard Young tableaux $(S,T)$ such that
\begin{itemize}
\item[(1)] ${\rm sh}(S)=\epsilon\tau$ for some $\tau\in \cP$ and   ${\rm sh}(T)=\lambda$,

\item[(2)] $\Delta_{\mf x}((T\rightarrow S)_R)\leq n$,
\end{itemize}
where $(T\rightarrow S)_R$ is the recording tableau of the insertion $T\rightarrow S$ \cite{Thomas} and  $\Delta_{\mf x}((T\rightarrow S)_R)$ is a new combinatorial statistic on the set of LR tableaux (see Definition \ref{Def of Delta}). The highest weight crystals of type $B_k$ and $C_k$ can be obtained by restricting the entries of $S$ and $T$ in $\B({\mf x}_\infty,\La^{\mf x}(\lambda,n))$ to $\{\,1,\ldots,k\,\}$.

The statistic $\Delta_{\mf x}$ is an essential part in the realization of $\B({\mf x}_\infty,\La^{\mf x}(\lambda,n))$, and it can be viewed as  a combinatorial realization of the crystal datum $\varepsilon_0^\ast$ which comes from the Kashiwara's $\ast$-crystal structure \cite{Kas90} on the negative part of $U_q({\mf x}_\infty)$ (see Remark \ref{epsilon star}). We remark that when $\lambda$ is the empty partition $(0)$, the tableau description of $\B({\mf x}_\infty,\Lambda^{\mf x}((0),n))$ combined with its Jacobi-Trudi or Weyl-Kac character formula recovers  partial Littlewood identities, which appeared in the study of symmetric plane partitions \cite{Mac95,St1,St2} (see also \cite{K09}).

We introduce another statistic $\nabla_{\mf x}$ on the set of LR tableaux (Definition \ref{LR coefficient}), which is a combinatorial realization of the crystal datum $\varepsilon_0$. Then as immediate consequences, we obtain a new description of the generalized LR coefficients and branching coefficients with respect to ${\mf l}_\infty$ (Theorems \ref{LR rule} and \ref{branching}). The multiplicities are given by the number of pairs of classical LR tableaux  with  certain restrictions on the statistics $\Delta_{\mf x}$ and $\nabla_{\mf x}$, which are always satisfied by highest weights in a stable range. The stable limits of the generalized LR coefficients and branching coefficients (here we fix $\lambda$ and take $n$ of $\Lambda^{\mf x}(\lambda,n)$ to infinity) look different from but simpler than the known results. Also, in case of $\mf{c}_\infty$, we can recover and hence extends in a bijective way  the well-known stable formulas using the reciprocity law via the $({\mf c}_\infty, {\rm Sp}(2n))$-Howe duality on a Fock space \cite{Wa} (Remarks \ref{classical stable formula} and \ref{classical stable formula'}).

We expect  a similar description for the highest weight crystals of type $D_\infty$, which is related with the product in (\ref{Littlewood identity})  when $\epsilon=0$. But our method of proof  based on folding crystals of type $A_\infty$ is no longer available for type $D_\infty$ and hence a different method is required. We hope to find a new approach  in the near future which does not depend on type.

\subsection{Super duality}
Lie superalgebras and their representations came from the study of super symmetries in theoretical  physics. The problem of finding their finite dimensional irreducible characters was one of the most important problems in this area after the  classification of finite dimensional simple Lie superalgebras by Kac \cite{K1}. In case of a general linear Lie superalgebra $\mf{gl}(m|n)$, it was solved by Serganova \cite{Se} and Brundan \cite{Br}.

Recently, the finite dimensional irreducible character problem for general linear and orthosymplectic Lie superalgebras was solved by  Cheng-Lam \cite{CL} and Cheng-Lam-Wang \cite{CLW}, respectively by a completely different method called {\it super duality}, which was conjectured in \cite{CWZ,CW} in case of general linear Lie superalgebras. The super duality is an equivalence between a parabolic category of  classical Lie algebras and that of general linear and orthosymplectic Lie superalgebras, respectively, and it gives a new insight to the representations of Lie superalgebras. In particular, it reveals a natural connection with the Kazhdan-Lusztig theory of Lie algebras and therefore provides  an irreducible character of  a Lie superalgebra in these categories (including finite dimensional irreducible representations) in terms of classical Kazhdan-Lusztig polynomials.

Another main motivation of the work in this paper came from the study of a combinatorial model for irreducible characters of Lie superalgebras.
A tableaux model is known for a special class of finite dimensional irreducible representations of
$\mf{gl}(m|n)$ or a queer Lie superalgebra ${\mf q}(n)$ appearing in a tensor power of their natural representations  \cite{BR,Kng-1,Sa,Wo,Sergeev}. But  very little is known about the other classes of simple Lie superalgebra as far as we know. There have been several works on tableau description of characters of orthosymplectic Lie superalgebras, which however are not necessarily irreducible (see \cite{BLR} and  also \cite{Kng-2} for a survey on the related works), and the previously known tableau models for the classical Lie algebras of type $B,C,D$ \cite{Be,KN,Kng,KoT87,Pr,Su}  do not seem to have a natural super analogue which works for orthosymplectic Lie superalgebras.

We note that a super Schur function or super symmetric function (in infinite variables) is  the  character of an irreducible module over a general linear Lie superalgebra of infinite rank, which  corresponds to an integrable highest weight module over the general linear Lie algebra of type $A_{+\infty}$ under the super duality functor \cite{CL}. Hence one may understand super duality as a categorical interpretation of superization in combinatorics.
Moreover, from a viewpoint of super duality, it is natural to define super Schur functions of type $B,C,D$ by the characters of irreducible modules over  orthosymplectic Lie superalgebras corresponding to  integrable highest weight modules over the Kac-Moody algebras of type $B_\infty$, $C_\infty$, $D_\infty$ under super duality, respectively.

We show that a superization of the bitableaux $(S,T)$ in $\B(\mf{x}_{\infty},\Lambda^{\mf x}(\lambda,n))$ is compatible with super duality and therefore it gives a tableau model for the super Schur functions of type $B$ and $C$  (Theorem \ref{application to super duality}). This can be viewed as an orthosymplectic analogue of  a hook Young tableau model for super symmetric functions \cite{BR,Kng-1}. Also by restricting the entries of  bitableaux  to a finite $\Z_2$-graded set, we get a tableau model for infinite dimensional irreducible characters of orthosymplectic Lie superalgebras of finite rank, which were studied in \cite{CKW,CZ} including unitarizable highest weight representations of Hermitian symmetric pairs.

We remark that we do not have a tableau description of finite dimensional irreducible representations of orthosymplectic Lie superalgebras of finite rank. By super duality, this problem is equivalent to finding a tableau description of  infinite dimensional irreducible representations of classical Lie algebras of finite rank in parabolic categories, which is in general much more difficult than the case of integrable highest weight modules.

\subsection{Organization}  The paper is organized as follows. In Section \ref{Preliminary}, we recall necessary background on crystals, Young tableaux and LR tableaux. In Section \ref{A Young tableau model for crystals of type $B$ and $C$}, we state the main result, a new combinatorial model for highest weight crystals of type $B_\infty$ and $C_\infty$. In Section \ref{Decomposition rule}, we discuss its generalized LR rule and branching rule as a crystal of type $A_{+\infty}$. In Section \ref{Character formula for orthosymplectic Lie superalgebras}, we give an  application to Lie superalgebras. Finally in Section \ref{Proof for BC}, we give a proof of the main theorem.

\subsection{Acknowledgement}
Most of this work was done while the author was visiting University of California, Berkeley during 2010-2011. He would like to thank N. Reshetikhin for his invitation and kind hospitality. He also thanks S.-J. Cheng and W. Wang who kindly explained their works on super duality, and A. Schilling for introducing SAGE  which helped him with computation and drawing of crystal graphs.

\section{Preliminary}\label{Preliminary}
\subsection{Review on crystals}
Let us briefly recall the notion of crystals  (see \cite{Kas94} for a general review and references therein).

Let $I$ be an index set. Let $\mathfrak{g}$ be a symmetrizable Kac-Moody algebra associated with a generalized Cartan matrix $A=(a_{ij})_{i,j\in I}$. Denote the weight lattice of $\mathfrak{g}$ by $P$, the set of simple roots by $\Pi=\{\,\alpha_i\,|\,i\in I\,\}\subset P$, and the set of simple coroots by $\Pi^\vee=\{\,h_i\,|\,i\in I\,\}\subset P^\vee$ with $\langle\alpha_j,h_i\rangle=a_{ij}$, where $\langle \cdot,\cdot\rangle$ is a natural pairing on $P^\vee\times P$.

A {\it $\mathfrak{g}$-crystal} (or {\it crystal} for short) is a set
$B$ together with the maps ${\rm wt} : B \rightarrow P$,
$\varepsilon_i, \varphi_i: B \rightarrow \mathbb{Z}\cup\{-\infty\}$ and
$\te_i, \tf_i: B \rightarrow B\cup\{{\bf 0}\}$ ($i\in I$) such that
for $b\in B$ and $i\in I$
\begin{itemize}
\item[(1)] $\varphi_i(b) =\langle {\rm wt}(b),h_i \rangle +
\varepsilon_i(b),$

\item[(2)]  $\varepsilon_i(\te_i b) = \varepsilon_i(b) - 1$, $\varphi_i(\te_i b) =
\varphi_i(b) + 1$, ${\rm wt}(\te_ib)={\rm wt}(b)+\alpha_i$ if $\te_i b \neq {\bf 0}$,

\item[(3)] $\varepsilon_i(\tf_i b) = \varepsilon_i(b) + 1$, $\varphi_i(\tf_i b) =
\varphi_i(b) - 1$, ${\rm wt}({\tf_i}b)={\rm wt}(b)-\alpha_i$ if $\tf_i b \neq {\bf 0}$,

\item[(4)] $\tf_i b = b'$ if and only if $b = \te_i
b'$ for $b, b' \in B$,

\item[(5)] $\te_ib=\tf_ib={\bf 0}$ if $\varphi_i(b)=-\infty$,
\end{itemize}
where ${\bf 0}$ is a formal symbol. Here we assume that $-\infty+n=-\infty$ for all $n\in\Z$.
Note that $B$
is equipped with an $I$-colored oriented graph structure, where
$b\stackrel{i}{\rightarrow}b'$ if and only if $b'=\tf_{i}b$ for
$b,b'\in B$ and $i\in I$.

We call $B$ {\it connected} if it is connected as a graph, and  call  $B$ {\it normal} if $\varepsilon_i(b)={\rm
max}\{\,k\,|\,\te_i^kb\neq {\bf 0}\,\}$ and $\varphi_i(b)={\rm
max}\{\,k\,|\,\tf_i^kb\neq {\bf 0}\,\}$ for $b\in B$ and $i\in I$.
The {\it dual crystal $B^\vee$ of $B$} is defined to be the set
$\{\,b^\vee\,|\,b\in B\,\}$ with ${\rm wt}(b^\vee)=-{\rm wt}(b)$, $\varepsilon_i(b^\vee)=\varphi_i(b)$,
$\varphi_i(b^\vee)=\varepsilon_i(b)$, $\te_i(b^\vee)=\left(\tf_i b \right)^\vee$ and
$\tf_i(b^\vee)=\left(\te_i b \right)^\vee$
for $b\in B$ and $i\in I$. We assume that ${\bf 0}^\vee={\bf
0}$.

Let $B_1$ and $B_2$ be crystals. A {\it morphism}
$\psi : B_1 \rightarrow B_2$ is a map from $B_1\cup\{{\bf 0}\}$ to
$B_2\cup\{{\bf 0}\}$ such that
\begin{itemize}
\item[(1)] $\psi(\bf{0})=\bf{0}$,

\item[(2)] ${\rm wt}(\psi(b))={\rm wt}(b)$,
$\varepsilon_i(\psi(b))=\varepsilon_i(b)$, and
$\varphi_i(\psi(b))=\varphi_i(b)$ if $\psi(b)\neq \bf{0}$,

\item[(3)] $\psi(\te_i b)=\te_i\psi(b)$ if $\psi(b)\neq \bf{0}$ and
$\psi(\te_i b)\neq \bf{0}$,

\item[(4)] $\psi(\tf_i
b)=\tf_i\psi(b)$ if $\psi(b)\neq
\bf{0}$ and $\psi(\tf_i b)\neq \bf{0}$,
\end{itemize}
for $b\in B_1$ and $i\in I$.
We call $\psi$ an {\it embedding} and $B_1$ a {\it subcrystal of}
$B_2$ when $\psi$ is injective, and call $\psi$ {\it strict} if
$\psi : B_1\cup\{{\bf 0}\} \rightarrow B_2\cup\{{\bf 0}\}$ commutes
with $\te_i$ and $\tf_i$ for all $i\in I$, where we assume that $\te_i{\bf
0}=\tf_i{\bf 0}={\bf 0}$.


A {\it tensor product $B_1\otimes B_2$} of crystals $B_1$ and $B_2$ is defined to be $B_1\times B_2$ as a set with elements  denoted by $b_1\otimes b_2$, where  {\allowdisplaybreaks
\begin{equation*}
\begin{split}
{\rm wt}(b_1\otimes b_2)&={\rm wt}(b_1)+{\rm wt}(b_2), \\
\varepsilon_i(b_1\otimes b_2)&= {\rm
max}\{\varepsilon_i(b_1),\varepsilon_i(b_2)-\langle {\rm
wt}(b_1),h_i\rangle\}, \\
\varphi_i(b_1\otimes b_2)&= {\rm max}\{\varphi_i(b_1)+\langle {\rm
wt}(b_2),h_i\rangle,\varphi_i(b_2)\},\\
{\te}_i(b_1\otimes b_2)&=
\begin{cases}
{\te}_i b_1 \otimes b_2, & \text{if $\varphi_i(b_1)\geq \varepsilon_i(b_2)$}, \\
b_1\otimes {\te}_i b_2, & \text{if
$\varphi_i(b_1)<\varepsilon_i(b_2)$},
\end{cases}\\
{\tf}_i(b_1\otimes b_2)&=
\begin{cases}
{\tf}_i b_1 \otimes b_2, & \text{if  $\varphi_i(b_1)>\varepsilon_i(b_2)$}, \\
b_1\otimes {\tf}_i b_2, & \text{if $\varphi_i(b_1)\leq
\varepsilon_i(b_2)$},
\end{cases}
\end{split}
\end{equation*}
\noindent for $i\in I$. Here we assume that ${\bf 0}\otimes
b_2=b_1\otimes {\bf 0}={\bf 0}$.} Then $B_1\otimes B_2$ is a
crystal. Note that $(B_1\otimes B_2)^\vee\simeq
B_2^\vee\otimes B_1^\vee$.

Let $P^+$ be the set of dominant integral weights for $\mathfrak{g}$. For $\Lambda\in P^+$, we denote by $\B(\mathfrak{g},\Lambda)$ the $\mathfrak{g}$-crystal associated with the highest weight module over the quantized enveloping algebra $U_q(\mathfrak{g})$ with highest weight $\Lambda$ and highest weight element $u_\Lambda$. Note that $\B(\mathfrak{g},\Lambda)$ is a connected normal $\mathfrak{g}$-crystal. For $\lambda\in P$, let
$T_{\lambda}=\{\,t_{\lambda}\,\}$ be a $\mathfrak{g}$-crystal with
${\rm wt}(t_{\lambda})=\lambda$  and
$\varepsilon_{i}(t_{\lambda})=\varphi_{i}(t_{\lambda})=-\infty$ for
$i\in I$.

Let $\Z[P]$ be the group ring of $P$ with a $\Z$-basis $\{\,e^\lambda\,|\,\lambda\in P\,\}$. For a crystal $B$, we define the character of $B$ by  ${\rm ch}B=\sum_{b\in B}e^{{\rm wt}(b)}$. Note that ${\rm ch}B$ is well-defined if and only if $B_\lambda=\{\,b\in B\,|\,{\rm wt}(b)=\lambda\,\}$ is finite for $\lambda\in P$.

\subsection{Affine Lie algebras of infinite rank}\label{affine Lie algebra}
Let us recall affine Lie algebras of infinite rank (see \cite{K} for detailed exposition).
Let $\Z^\times=\Z\setminus\{0\}$ and let $\gl_\infty$ be the Lie algebra spanned by complex matrices
$(a_{i,j})_{i,j\in \Z^\times}$ with finitely many non-zero entries, which is
spanned by the elementary matrix $E_{i,j}$ ($i,j\in\Z^\times$). Let
$\h=\bigoplus_{i\in \Z^\times}\C E_{i,i}$  and $\langle\cdot,\cdot\rangle$ denote the natural
pairing on $\frak{h}^*\times\frak{h}$. Denote  the set of simple coroots by $\Pi^{\vee}=\{\, h_{-i}=E_{-i-1, -i-1}-E_{-i, -i},\
h_i=E_{i,i}-E_{i+1,i+1}\, (i\in \Z_{>0}),\ h_0=E_{-1,-1}-E_{1,1} \, \}$, and
the set of simple roots  by $\Pi=\{\,
\alpha_{-i}=\epsilon_{-i-1}-\epsilon_{-i}\,\
\alpha_i=\epsilon_i-\epsilon_{i+1} \ (i\in \Z_{>0}),
\alpha_0=\epsilon_{-1}-\epsilon_1\, \}$, where $\epsilon_i\in\h^*$ is
given by $\langle \epsilon_i,E_{j,j}\rangle=\delta_{ij}$.
The Dynkin diagram associated with the
Cartan matrix
$\left(\langle\alpha_j,h_i\rangle\right)_{i,j\in\mathbb{Z}}$ is of type
\begin{center}
\hskip 2cm\setlength{\unitlength}{0.45cm}
\begin{picture}(15,2.5)(0,-0.5)
\put(-4,.2){\makebox(0,0)[c]{$A_{\infty}\ :$}}

\put(-.85,.5){\line(1,0){1}}
\put(0,0){\makebox(1,1){$\bigcirc$}}\put(.85,.5){\line(1,0){1}}
\put(5,0){\makebox(1,1){$\bigcirc$}}\put(3.85,.5){\line(1,0){1.2}}\put(5.85,.5){\line(1,0){1.2}}
\put(7,0){\makebox(1,1){$\bigcirc$}}\put(7.85,.5){\line(1,0){1.2}}
\put(9,0){\makebox(1,1){$\bigcirc$}}\put(9.85,.5){\line(1,0){1.2}}
\put(14,0){\makebox(1,1){$\bigcirc$}}\put(12.85,.5){\line(1,0){1.2}}\put(14.85,.5){\line(1,0){1.2}}
\put(2.5,.2){$\cdots$}\put(11.5,.2){$\cdots$}
\put(-2,.2){$\cdots$}\put(16.5,.2){$\cdots$}

\put(0.5,-.5){\makebox(0,0)[c]{\tiny ${\alpha}_{-n}$}}
\put(5.5,-.5){\makebox(0,0)[c]{\tiny ${\alpha}_{-1}$}}
\put(7.5,-.5){\makebox(0,0)[c]{\tiny ${\alpha}_0$}}
\put(9.5,-.5){\makebox(0,0)[c]{\tiny ${\alpha}_1$}}
\put(14.5,-.5){\makebox(0,0)[c]{\tiny ${\alpha}_n$}}

\end{picture}\hskip 1.2cm .
\end{center}\vskip 2mm
The weight lattice of $\gl_{\infty}$ is $P=\Z\Lambda_0\oplus
\bigoplus_{i\in\Z}\Z\epsilon_i=\bigoplus_{i\in\Z}\Z\Lambda_i$, where $\Lambda_i$ is the $i$-th fundamental weight, that is $\langle \Lambda_i,h_j\rangle =\delta_{ij}$.
Note that
\begin{align}
\Lambda_i=
\begin{cases}
\Lambda_0- \epsilon_{-1}-\cdots-\epsilon_{i}, & \text{if  $i< 0$}, \\
\Lambda_0+ \epsilon_1+\cdots+\epsilon_i, & \text{if  $i> 0$}.
\end{cases}
\end{align}
For a subset $S$ of $\Z^\times$, we denote by $\gl_{S}$ the subalgebra of $\gl_\infty$ spanned by $E_{ij}$ for $i,j\in S$. For simplicity, we put $\gl_{>0}=\gl_{\Z_{>0}}$ and  $\gl_{<0}=\gl_{\Z_{<0}}$, which are of type $A_{+\infty}$.

Let us realize the root systems of type $B_\infty$ and $C_\infty$  by  folding the diagram of $A_\infty$ \cite{Kas96}.
From now on, we assume that $\epsilon\in \{1,2\}$ and the symbol ${\mf x}$ denotes ${\mf b}$ (resp. ${\mf c}$) if $\epsilon=1$ (resp. $\epsilon=2$).
Let
\begin{equation}
\begin{split}
P_{\mf x}&=\{\,\lambda\in
P\,|\,\tfrac{1}{\epsilon}\langle\lambda,h_0\rangle\in\Z,\
\langle\lambda,h_i\rangle=\langle\lambda,h_{-i}\rangle\ (i\in\mathbb{Z}_{> 0})\,\},\\
\Pi^\vee_{\mf x}&= \{\, \widehat{h}_0= \tfrac{1}{\epsilon} h_0\ ,\  \widehat{h}_i=h_i \ (i\in \Z_{>0}) \,\}\subset \Pi^\vee, \\
{\Pi}_{\mf x}&=\{\,\widehat{\alpha}_0=\epsilon\alpha_0, \ \
\widehat{\alpha}_i=\alpha_{i}+\alpha_{-i} \ \ (i\in \Z_{>0})\,\}\subset \Pi.\\
\end{split}
\end{equation}
Then the Dynkin diagram associated with the
Cartan matrix
$A_{\mf x}=\left(\langle\widehat{\alpha}_j,\widehat{h}_i\rangle\right)_{i,j\in \Z_{\geq 0}}$ is of type
\begin{center}
\hskip -1cm  \setlength{\unitlength}{0.16in}
\begin{picture}(24,4)
\put(2,2){\makebox(0,0)[c]{$B_{\infty}\ :$}}

\put(5.6,2){\makebox(0,0)[c]{$\bigcirc$}}
\put(8,2){\makebox(0,0)[c]{$\bigcirc$}}
\put(10.4,2){\makebox(0,0)[c]{$\bigcirc$}}
\put(14.85,2){\makebox(0,0)[c]{$\bigcirc$}}
\put(17.25,2){\makebox(0,0)[c]{$\bigcirc$}}
\put(19.4,2){\makebox(0,0)[c]{$\bigcirc$}}
\put(8.35,2){\line(1,0){1.5}} \put(10.82,2){\line(1,0){0.8}}
\put(13.2,2){\line(1,0){1.2}} \put(15.28,2){\line(1,0){1.45}}
\put(17.7,2){\line(1,0){1.25}} \put(19.81,2){\line(1,0){0.9}}
\put(6.8,1.97){\makebox(0,0)[c]{$\Longleftarrow$}}
\put(12.5,1.95){\makebox(0,0)[c]{$\cdots$}}
\put(21.5,1.95){\makebox(0,0)[c]{$\cdots$}}
\put(5.6,1){\makebox(0,0)[c]{\tiny $\widehat{\alpha}_0$}}
\put(8,1){\makebox(0,0)[c]{\tiny $\widehat{\alpha}_1$}}
\put(10.4,1){\makebox(0,0)[c]{\tiny $\widehat{\alpha}_2$}}
\put(15,1){\makebox(0,0)[c]{\tiny $\widehat{\alpha}_{n-1}$}}
\put(17.4,1){\makebox(0,0)[c]{\tiny $\widehat{\alpha}_n$}}
\put(19.8,1){\makebox(0,0)[c]{\tiny $\widehat{\alpha}_{n+1}$}}

\put(25.5,2){\makebox(0,0)[c]{($\epsilon=1$)}}

\end{picture}
\end{center}

\begin{center}
\hskip -1cm  \setlength{\unitlength}{0.16in}
\begin{picture}(24,4)
\put(2,2){\makebox(0,0)[c]{$C_{\infty}\ :$}}

\put(5.6,2){\makebox(0,0)[c]{$\bigcirc$}}
\put(8,2){\makebox(0,0)[c]{$\bigcirc$}}
\put(10.4,2){\makebox(0,0)[c]{$\bigcirc$}}
\put(14.85,2){\makebox(0,0)[c]{$\bigcirc$}}
\put(17.25,2){\makebox(0,0)[c]{$\bigcirc$}}
\put(19.4,2){\makebox(0,0)[c]{$\bigcirc$}}
\put(8.35,2){\line(1,0){1.5}} \put(10.82,2){\line(1,0){0.8}}
\put(13.2,2){\line(1,0){1.2}} \put(15.28,2){\line(1,0){1.45}}
\put(17.7,2){\line(1,0){1.25}} \put(19.81,2){\line(1,0){0.9}}
\put(6.8,1.97){\makebox(0,0)[c]{$\Longrightarrow$}}
\put(12.5,1.95){\makebox(0,0)[c]{$\cdots$}}
\put(21.5,1.95){\makebox(0,0)[c]{$\cdots$}}
\put(5.6,1){\makebox(0,0)[c]{\tiny $\widehat{\alpha}_0$}}
\put(8,1){\makebox(0,0)[c]{\tiny $\widehat{\alpha}_1$}}
\put(10.4,1){\makebox(0,0)[c]{\tiny $\widehat{\alpha}_2$}}
\put(15,1){\makebox(0,0)[c]{\tiny $\widehat{\alpha}_{n-1}$}}
\put(17.4,1){\makebox(0,0)[c]{\tiny $\widehat{\alpha}_n$}}
\put(19.8,1){\makebox(0,0)[c]{\tiny $\widehat{\alpha}_{n+1}$}}

\put(25.5,2){\makebox(0,0)[c]{($\epsilon=2$)}}
\end{picture}
\end{center}\vskip 5mm
Let ${\mf x}_\infty$ be the Kac-Moody algebra associated with  $(A_{\mf x},P_{\mf x},\Pi^\vee_{\mf x}, \Pi_{\mf x})$ (see \cite{K} or Section \ref{Osp} for an explicit construction of ${\mf x}_\infty$).
The fundamental weights $\Lambda^{\mf x}_i$ ($i\in \Z_{\geq  0}$) for ${\mf x}_\infty$ (that is, $\langle \Lambda^{\mf x}_i,\widehat{h}_j\rangle=\delta_{ij}$) are given by
\begin{equation}
\begin{split}
{\Lambda}_i^{\mf x}=
\begin{cases}
\epsilon\Lambda_0, & \text{if $i=0$},\\
\Lambda_i+\Lambda_{-i}=2\La_0+\widehat{\epsilon}_1+\cdots+\widehat{\epsilon}_i, & \text{if $i\geq 1$},
\end{cases}
\end{split}
\end{equation}
where $\widehat{\epsilon}_i=\epsilon_i-\epsilon_{-i}$.
Note that 
\begin{equation}
P_{\mf x}=\Z{\Lambda}_0^{\mf x}\oplus
\bigoplus_{i\in\Z_{> 0}}\Z\widehat{\epsilon}_i.
\end{equation}
\vskip 2mm

For $k\geq 2$, we denote by ${\mf x}_k$  the subalgebra of ${\mf x}_\infty$, which is a finite dimensional simple Lie algebra of type $B_k$ ($\epsilon=1$) or $C_k$ ($\epsilon=2$) associated to $\{\,\widehat{h}_i\,|\, 0\leq i\leq k-1\,\}$ and $\{\,\widehat{\alpha}_i\,|\, 0\leq i\leq k-1\,\}$. Also, let $\mathfrak{l}_\infty$ be the subalgebra of ${\mf x}_\infty$ associated to $\{\,\widehat{h}_i\,|\, i\in\Z_{>0}\,\}$ and $\{\,\widehat{\alpha}_i\,|\, i\in\Z_{>0}\,\}$, which is isomorphic to  $\gl_{>0}$.
\vskip 2mm

Let us introduce some notations. Let $\cP$ denote the set of partitions.
Following \cite{Mac95}, we denote the length of $\lambda=(\lambda_i)_{i\geq 1}\in \cP$ by
$\ell(\lambda)$ and the conjugate of
$\lambda$ by $\lambda'=(\lambda'_i)_{i\geq 1}$.
We identify
$\lambda$ with a Young diagram. But here we use a different convention for Young diagrams, where in each $i$-th row from the bottom there are $\lambda_i$'s  boxes placed in a right-justified way. We assume that the rows (resp. columns) in a Young diagram are enumerated from the bottom (resp. right) unless otherwise specified. We denote by $\lambda/\mu$ the skew Young diagram for $\mu=(\mu_i)_{i\geq 1}\in\cP$ with $\mu_i\leq \lambda_i$ ($i\geq 1$). We denote by $|\lambda/\mu|$ the number of boxes in $\lambda/\mu$.   For example,

$$(5,3,2) \ \ = \ \
{\def\lr#1{\multicolumn{1}{|@{\hspace{.6ex}}c@{\hspace{.6ex}}|}{\raisebox{-.3ex}{$#1$}}}\raisebox{-.6ex}
{$\begin{array}{ccccc}
\cline{4-5}
& & &\lr{\ \ } & \lr{\ \ } \\
\cline{3-5}
& & \lr{\ \ } & \lr{\ \ }& \lr{\ \ }\\
\cline{1-5}
\lr{\ \ } & \lr{\ \ } & \lr{\ \ } & \lr{\ \ }& \lr{\ \ }\\
\cline{1-5}
\end{array}$}}
\ \ \  \ \ \
(5,3,2)/(3,1) \ \ = \ \
{\def\lr#1{\multicolumn{1}{|@{\hspace{.6ex}}c@{\hspace{.6ex}}|}{\raisebox{-.3ex}{$#1$}}}\raisebox{-.6ex}
{$\begin{array}{ccccc}
\cline{4-5}
& & &\lr{\ \ } & \lr{\ \ } \\
\cline{3-5}
& & \lr{\ \ } & \lr{\ \ }& \cdot \\
\cline{1-4}
\lr{\ \ } & \lr{\ \ } & \cdot & \cdot & \cdot \\
\cline{1-2}
\end{array}$}}\ \ \ .
$$

Now, let
\begin{equation}
\begin{split}
\cP({\mf x})&=
\{\,(\lambda,n)\,|\,\lambda\in\cP, n\in\Z_{\geq 0}, \ 2\lambda_1\leq  \epsilon\,n\,\} .
\end{split}
\end{equation}
For $(\lambda,n)\in \cP ({\mf x})$, define
\begin{equation}\label{Lambda lambda n}
\La^{\mf x}({\lambda,n})=
n\Lambda^{\mf x}_{0}+\sum_{i\geq 1}\lambda_i\widehat{\epsilon}_i.
\end{equation}
Then $P_{\mf x}^+=\{\,\La^{\mf x}(\lambda,n)\,|\, (\lambda,n)\in \cP({\mf x})\,\} $ is the set of dominant integral weights for ${\mf x}_\infty$. Also, we put
\begin{equation}
\cP_{\mf x} =\{\,\epsilon\lambda=(\epsilon\lambda_i)_{i\geq 1}\,|\, \lambda=(\lambda_i)_{i\geq 1}\in\cP \,\}.
\end{equation}

\subsection{Young tableaux}
Let $(\A,<)$ be a linearly
ordered  set with a $\mathbb{Z}_2$-grading $\A=\A_0\sqcup\A_1$.  For
$a\in \A_p$ ($p\in\Z_2$), we put $|a|=p$.
Unless otherwise stated, we assume that $\Z$ or its subset is equipped with a usual linear ordering, where all the elements are of degree $0$.

For a skew Young diagram $\lambda/\mu$, a tableau $T$ obtained by
filling $\lambda/\mu$ with entries in $\A$ is called
$\A$-semistandard  if (1) the entries in each row (resp. column) are
weakly increasing from left to right (resp. from top to bottom), (2)
the entries in $\A_0$ (resp. $\A_1$) are strictly increasing in each
column (resp. row). When $\A=\Z_{>0}$, $T$ is usually  called a (semistandard) Young tableau. We say that $\lambda/\mu$ is the shape of
$T$, and write ${\rm sh}(T)=\lambda/\mu$. For $a\in \A$, let $m_a(T)$ be  the number of occurrences of $a$ in $T$.
We denote by ${ SST}_{\A}(\lambda/\mu)$ the set of all
$\A$-semistandard tableaux of shape $\lambda/\mu$.

Let ${\bf x}_{\A}=\{\,x_a\,|\,a\in \A \,\}$ be the set of formal commuting variables. We let $S_{\lambda/\mu}({\bf x}_{\A})=\sum_{T\in  SST_\A(\lambda/\mu)}x_{\A}^T$ be the super Schur function associated to $\lambda$, where $x_{\A}^T=\prod_{a\in\A}x_a^{m_a(T)}$.  Note that $S_{\lambda}({\bf x}_{\Z_{>0}})$ is the usual Schur function for $\lambda\in\cP$, which we will denote by $s_\lambda({\bf x})$  for simplicity.

Let $\W_{\A}$ be the set of finite words with letters in $\A$. For
$T\in  SST_{\A}(\lambda/\mu)$, we denote by $w_{\rm col}(T)\in \W_\A$
the word obtained by reading the entries of $T$ column by column
from right to left, and from top to bottom in each column. Also, we
denote by $w_{\rm row}(T)\in\W_\A$ the word obtained by reading the
entries of $T$ row by row from top to bottom, and from right to left
in each row. We will denote by $w(T)$ the row or column word of $T$ when a statement holds for both words of $T$.

For $a\in\A$ and $S\in  SST_{\A}(\mu)$, let us do the Schendsted's row insertion of $a$ starting from the bottom row of $S$, and then proceed to the upper rows as far as possible. We denote the resulting tableau by $a\rightarrow S$ (cf. \cite{BR,Fu}). For $w=w_1\ldots w_r\in\W_\A$, we let
$(w\rightarrow S)=( w_r\rightarrow(\cdots(w_1\rightarrow S)\cdots))$.  For $T\in  SST_{\A}(\nu)$, we define $(T\rightarrow S)=(w(T)\rightarrow S)$.

Note that $\Z^{\times}$ is a normal $\gl_{\infty}$-crystal given by

\begin{equation}
\cdots \ \stackrel{-2}{\longrightarrow} -2 \
\stackrel{-1}{\longrightarrow} -1 \stackrel{0} {\longrightarrow} 1
\stackrel{1}{\longrightarrow} 2\stackrel{2}{\longrightarrow} \cdots
\end{equation}\vskip 2 mm
\noindent
with ${\rm wt}(i)=\epsilon_i$. If $\A$ is a subset of $\Z^\times$, we may regard $\A$ as a normal $\gl_\A$-crystal whose associated  graph is induced from that of $\Z^{\times}$ above by restricting the vertices to $\A$.
Then by identifying $w_1\ldots w_r\in\W_\A$ with $w_1\otimes\cdots\otimes w_r\in\A^{\otimes r}$, $\W_\A$ is a $\gl_\A$-crystal under tensor product of crystals. Also, the image of $SST_\A(\lambda/\mu)$ in $\W_\A$
under the map $T \mapsto w(T)=w_1\ldots w_r$  together with
$\{{\bf 0}\}$ is invariant under $\te_i,\tf_i$ for all $i$ such that $\alpha_i$ is a simple root of $\gl_\A$. Hence
$SST_\A(\lambda/\mu)$ is a subcrystal of $\W_\A$ \cite{KN}.

\subsection{LR tableaux}
For $\lambda,\mu,\nu\in\cP$ with $|\lambda|=|\mu|+|\nu|$, let ${\rm \bf LR}^\lambda_{\mu
\nu}$ be the set of tableaux $U$ in
$ SST_{\Z_{>0}}(\lambda/\mu)$ such that
\begin{itemize}
\item[(1)] $m_i(U)=\nu_i$  for $i\geq 1$,

\item[(2)] the number of occurrences of each $i\geq 1$ in
$w_1\ldots w_k$ is no less than that of $i+1$ in $w_1\ldots w_k$ for $1\leq k\leq |\nu|$,
where $w(U)=w_1\ldots w_{|\nu|}$.
\end{itemize}
We call   ${\rm \bf LR}^\lambda_{\mu \nu}$ the set of {\it
LR tableaux of shape $\lambda/\mu$ with content
$\nu$} and put $c^\lambda_{\mu \nu}=\left|{\rm \bf LR}^\lambda_{\mu
\nu}\right|$. Our definition is not the same as that of a usual LR  tableau  \cite{Mac95}, but there is a natural one-to-one correspondence between them. (For example, ${\rm \bf LR}^\lambda_{\mu \nu}$ is in one-to-one correspondence with $\ov{\rm \bf LR}^\lambda_{\mu \nu}$ in \cite[Section 3.2]{K10} under $180^\circ$ rotation and replacing $i$ with $-i$. Then use \cite[Proposition 3.2]{K10}.)

For $S\in   SST_{\A}(\mu)$ and $T\in {SST}_{\A}(\nu)$, suppose that
$\lambda={\rm sh}(T\rightarrow S)$ and $w(T)=w_1\ldots
w_{|\nu|}$. Define $(T\rightarrow S)_R$ to be the tableau of shape
$\lambda/\mu$ such that for $1\leq k\leq |\nu|$ ${\rm sh}(w_1\ldots w_k\rightarrow S)/{\rm
sh}(w_1\ldots w_{k-1}\rightarrow S)$ is filled with $i$ when ${w_k}$
is an entry in the $i$-th row of $T$. Then the map
$(S,T)\mapsto \left((T\rightarrow S),(T\rightarrow S)_R\right)$
gives a bijection \cite{Thomas}
\begin{equation}\label{LR}
 SST_{\A}(\mu)\times  SST_{\A}(\nu)\longrightarrow
\bigsqcup_{\lambda\in\cP} SST_{\A}(\lambda)\times {\rm \bf LR}^{\lambda}_{\mu\,\nu}.
\end{equation}
This implies that  $S_\mu({\bf x}_\A)S_\nu({\bf x}_\A)=\sum_{\lambda}c^\lambda_{\mu \nu}S_\lambda({\bf x}_\A)$. Let $Sym_\A$ be the $\mathbb{Q}$-span of $S_\lambda({\bf x}_\A)$'s and put $Sym=Sym_{\Z_{>0}}$ for simplicity. Then $Sym_{\A}$ is a $\mathbb{Q}$-algebra with a basis $\{\,S_\lambda({\bf x}_\A)\,|\,\lambda\in \cP \text{ with } SST_\A(\lambda)\neq \emptyset\,\}$ \cite{BR}, and  there is a surjective algebra homomorphism
\begin{equation}\label{Omega_A}
\omega_{\A} : Sym \longrightarrow  Sym_\A
\end{equation}
sending $s_\lambda({\bf x})$ to $S_\lambda({\bf x}_\A)$. Also $\omega_\A$ is an isomorphism when $\A$ is infinite.

Let ${\rm LR}^\lambda_{\mu\nu}$ be the set of $V\in
 SST_{\Z_{>0}}(\nu)$ such that $(V\rightarrow H_\mu)= H_\lambda$, where $H_{\lambda}$ denotes the highest weight element in a crystal $ SST_{\Z_{>0}}(\lambda )$ that is, the $k$-th entry from the top of each column  is $k$.
There is a bijection
\begin{equation}\label{imath}
{{\rm \bf LR}}^\lambda_{\mu\nu} \stackrel{(\cdot)^\sharp}{\longrightarrow} {\rm \large LR}^\lambda_{\mu \nu},
\end{equation}
where
$U^\sharp=V$ with $U=(V\rightarrow H_{\mu})_R$.
Indeed,  the number of
$k$'s in the $i$-th row of $V$ is equal to the number of $i$'s in
the $k$-th row of $U$ for $i, k\geq 1$ \cite{Na}.

\begin{ex}{\rm\mbox{}
\begin{center}
$U=$\ \ \
${\def\lr#1{\multicolumn{1}{|@{\hspace{.6ex}}c@{\hspace{.6ex}}|}{\raisebox{-.3ex}{$#1$}}}\raisebox{-.6ex}
{$\begin{array}{ccccc}
\cline{5-5}
& & & & \lr{1} \\
\cline{3-5}
& & \lr{1} & \lr{1}& \lr{2}\\
\cline{2-5}
  & \lr{1} & \lr{2} & \lr{3}& \cdot \\
\cline{1-4}
\lr{2} & \lr{3} & \cdot & \cdot &  \cdot\\
\cline{1-2}
\end{array}$}}$
\ \ \ $\in {\rm \bf LR}^{(5,4,3,1)}_{(3,1)\ (4,3,2)} $\ \ \ \ \,
$U^\sharp=$
  \ \ \
${\def\lr#1{\multicolumn{1}{|@{\hspace{.6ex}}c@{\hspace{.6ex}}|}{\raisebox{-.3ex}{$#1$}}}\raisebox{-.6ex}
{$\begin{array}{cccc}
\cline{3-4}
 & &\lr{1} & \lr{2} \\
\cline{2-4}
 & \lr{1} & \lr{2}& \lr{3}\\
\cline{1-4}
 \lr{2} & \lr{3} & \lr{3}& \lr{4}\\
\cline{1-4}
\end{array}$}}$\ \ \ $\in {\rm LR}^{(5,4,3,1)}_{(3,1)\ (4,3,2)}$
\end{center}\vskip 2mm}
\end{ex}

\section{A Young tableau model for crystals of type $B$ and $C$}\label{A Young tableau model for crystals of type $B$ and $C$}
\subsection{Crystals associated with Littlewood identities}\label{crystal Txlambda}
\begin{df}{\rm
For $\lambda\in \cP$, let
\begin{equation*}
{\mathcal{T}}^{\mf x}(\lambda)=
\bigsqcup_{\tau\in\cP_{\mf x}}SST_{\Z_{>0}}(\tau)\times SST_{\Z_{>0}}(\lambda).\end{equation*}}
\end{df}\vskip 2mm

Let us define an ${\mf x}_\infty$-crystal structure on ${\mathcal{T}}^{\mf x}(\lambda)$.
First, we may regard ${\mathcal{T}}^{\mf x}(\lambda)=\bigsqcup_{\tau} SST_{\Z_{>0}}(\tau)\otimes SST_{\Z_{>0}}(\lambda)$ as an $\mathfrak{l}_\infty$-crystal (Recall that $\mathfrak{l}_\infty \simeq \gl_{>0}$). We denote by $\td{E}_i$ and $\td{F}_i$ the operators corresponding to $\widehat{\alpha}_i$ for $i\in \Z_{>0}$.
As an $\mathfrak{l}_\infty$-crystal, ${\mathcal{T}}^{\mf x}(\lambda)$ decomposes as follows:
\begin{equation}\label{l infty decomposition}
{\mathcal{T}}^{\mf x}(\lambda)\simeq \bigsqcup_{\sigma\in\cP}\bigsqcup_{\tau\in\cP_{\mf x}} SST_{\Z_{>0}}(\sigma)\times {\bf LR}^{\sigma}_{\tau \lambda},
\end{equation}
where the isomorphism is given by sending $(S,T)$ to $((T\rightarrow S),(T\rightarrow S)_R)$ and $\td{E}_i$, $\td{F}_i$ ($i\in \Z_{>0}$) act on the first component of the right-hand side (see for example \cite[Proposition 4.17]{KK}).

Next, let us define the operators $\td{E}_0$ and $\td{F}_0$ on ${\mathcal{T}}^{\mf x}(\lambda)$ corresponding to $\widehat{\alpha}_0$.

\noindent \textsc{Case 1}. Suppose that ${\mf x}={\mf b}$ and  $(S,T)\in \mathcal{T}^{\mf b}(\lambda)$ is given with ${\rm sh}(S)=\tau$ for some $\tau\in\cP$. For each $k\geq 1$, let $s_k$  be the entry
in the top of the $k$-th column of $S$. We put
\begin{equation}\label{signs}
\sigma_k=
\begin{cases}
+ \ , & \text{if the $k$-th column is empty}, \\
+ \ ,& \text{if   $s_k>1$}, \\
- \ ,& \text{if   $s_k=1$}.
\end{cases}
\end{equation}
In the sequence $\sigma=(\ldots,\sigma_2,\sigma_1)$, we replace a
pair $(\sigma_{s'},\sigma_{s})=(-,+)$, where $s'>s$ and
$\sigma_t=\cdot$ for $s'>t>s$, with $(\,\cdot\,,\,\cdot\,)$, and repeat
this process as far as possible until we get a sequence with no $-$
placed to the left of $+$.  We denote  this sequence by $\td{\sigma}$.

 We define $\td{E}_0 (S,T)=(S',T)$, where $S'$ is the tableau
obtained from $S$  by removing  $\boxed{1}$ in the column  corresponding to the left-most $-$ in $\td{\sigma}$. We call this $\boxed{1}$ {\it removable}.
If there is no such $-$ sign or removable $\boxed{1}$, then we define $\td{E}_0 (S,T)={\bf
0}$.

We define $\td{F}_0  (S,T)=(S'',T)$ where $S''$ is the tableau obtained from $S$ by adding   $\boxed{1}$ on top of the
column  corresponding to the right-most $+$ in $\td{\sigma}$. If there is no such $+$ sign, then we
define $\td{F}_0(S,T)={\bf 0}$.  \vskip 2mm

\noindent \textsc{Case 2}. Suppose that ${\mf x}={\mf c}$ and  $(S,T)\in \mathcal{T}^{\mf c}(\lambda)$ is given with ${\rm sh}(S)=2\tau$ for some $\tau\in\cP$. For each $k\geq 1$,  consider $(s_{2k},s_{2k-1})$  the pair of entries  in the top of the $2k$-th  and $(2k-1)$-st columns of $S$, respectively. Note that  $s_{2k-1}$ and $s_{2k}$ are placed in the same row of $2\tau$ and $s_{2k}\leq s_{2k-1}$. We put
\begin{equation}\label{signs}
\sigma_k=
\begin{cases}
+ \ , & \text{if  $s_{2k},s_{2k-1}\geq 2$ or the $(2k-1)$-st column is empty}, \\
- \ ,& \text{if $s_{2k}=s_{2k-1}=1$}, \\
\ \cdot \ \, ,& \text{otherwise}.
\end{cases}
\end{equation}

As in \textsc{Case 1}, we obtain a reduced  sequence $\td{\sigma}$ from $\sigma=(\ldots,\sigma_2,\sigma_1)$.
 We define $\td{E}_0 (S,T)=(S',T)$, where $S'$ is the tableau
obtained from $S$  by removing  a domino
${\def\lr#1{\multicolumn{1}{|@{\hspace{.6ex}}c@{\hspace{.6ex}}|}{\raisebox{-.3ex}{$#1$}}}\raisebox{-.6ex}
{$\begin{array}[b]{cc}
\cline{1-1}\cline{2-2}
\lr{1}&\lr{1}\\
\cline{1-1}\cline{2-2}
\end{array}$}}$
 in the pair of columns  corresponding to the left-most $-$ in $\td{\sigma}$. We call this domino
 ${\def\lr#1{\multicolumn{1}{|@{\hspace{.6ex}}c@{\hspace{.6ex}}|}{\raisebox{-.3ex}{$#1$}}}\raisebox{-.6ex}
{$\begin{array}[b]{cc}
\cline{1-1}\cline{2-2}
\lr{1}&\lr{1}\\
\cline{1-1}\cline{2-2}
\end{array}$}}$
  {\it removable}.
If there is no such $-$ sign or removable domino
${\def\lr#1{\multicolumn{1}{|@{\hspace{.6ex}}c@{\hspace{.6ex}}|}{\raisebox{-.3ex}{$#1$}}}\raisebox{-.6ex}
{$\begin{array}[b]{cc}
\cline{1-1}\cline{2-2}
\lr{1}&\lr{1}\\
\cline{1-1}\cline{2-2}
\end{array}$}}$
, then we define $\td{E}_0 (S,T)={\bf
0}$.

We define $\td{F}_0  (S,T)=(S'',T)$ where $S''$ is the tableau obtained from $S$ by adding  a domino
${\def\lr#1{\multicolumn{1}{|@{\hspace{.6ex}}c@{\hspace{.6ex}}|}{\raisebox{-.3ex}{$#1$}}}\raisebox{-.6ex}
{$\begin{array}[b]{cc}
\cline{1-1}\cline{2-2}
\lr{1}&\lr{1}\\
\cline{1-1}\cline{2-2}
\end{array}$}}$
 on top of the pair of
columns  corresponding to the right-most $+$ in $\td{\sigma}$. If there is no such $+$ sign, then we
define $\td{F}_0(S,T)={\bf 0}$.

\begin{lem}
${\mathcal{T}}^{\mf x}(\lambda)\cup\{{\bf 0}\}$ is invariant under $\td{E}_0$ and $\td{F}_0$.
\end{lem}
\pf It is straightforward to check from the definitions of $\td{E}_0$ and $\td{F}_0$.
\qed

\begin{ex}\label{Nex}{\rm Let
\begin{equation*}
\begin{split}
& (S,T)=\left(\ \
\text{${\def\lr#1{\multicolumn{1}{|@{\hspace{.6ex}}c@{\hspace{.6ex}}|}{\raisebox{-.3ex}{$#1$}}}\raisebox{-.6ex}
{$\begin{array}{ccccc}
\cline{5-5}
& & & & \lr{1 } \\
\cline{3-5}
& & \lr{1 } & \lr{3}& \lr{3}\\
\cline{1-5}
\lr{1 } & \lr{2 } & \lr{3 } & \lr{4 }& \lr{5 }\\
\cline{1-5}
\end{array}$}}$}\ \ \ , \ \ \
\text{${\def\lr#1{\multicolumn{1}{|@{\hspace{.6ex}}c@{\hspace{.6ex}}|}{\raisebox{-.3ex}{$#1$}}}\raisebox{-.6ex}
{$\begin{array}{cc}
& \\
\cline{2-2}
&\lr{2}\\
\cline{1-1}\cline{2-2}
\lr{2}& \lr{3}\\
\cline{1-1}\cline{2-2}
\end{array}$}}$}\ \ \  \right) \in {\mathcal{T}}^{\mf b}(2,1).
\end{split}
\end{equation*}
Then
\begin{equation*}
\begin{split}
&\sigma=(\
\cdots \ + \ +\ -\  + \ -    \  +\   - \ ), \\
&\td{\sigma}=(\
\cdots \ + \ +\ \ \cdot \ \ \cdot \ \ \, \cdot  \  \ \, \cdot \ \ \,  - \ ).
\end{split}
\end{equation*}
Hence we have
\begin{equation*}
\begin{split}
\td{E}_0(S,T)&=
\left(\ \
\text{${\def\lr#1{\multicolumn{1}{|@{\hspace{.6ex}}c@{\hspace{.6ex}}|}{\raisebox{-.3ex}{$#1$}}}\raisebox{-.6ex}
{$\begin{array}{ccccc}
\cline{3-5}
& & \lr{1 } & \lr{3}& \lr{3}\\
\cline{1-5}
\lr{1 } & \lr{2 } & \lr{3 } & \lr{4 }& \lr{5 }\\
\cline{1-5}
\end{array}$}}$}\ \ \ , \ \ \
\text{${\def\lr#1{\multicolumn{1}{|@{\hspace{.6ex}}c@{\hspace{.6ex}}|}{\raisebox{-.3ex}{$#1$}}}\raisebox{-.6ex}
{$\begin{array}{cc}
\cline{2-2}
&\lr{2}\\
\cline{1-1}\cline{2-2}
\lr{2}& \lr{3}\\
\cline{1-1}\cline{2-2}
\end{array}$}}$}\ \ \  \right), \\
\td{F}_0(S,T)&=
\left(\ \
\text{${\def\lr#1{\multicolumn{1}{|@{\hspace{.6ex}}c@{\hspace{.6ex}}|}{\raisebox{-.3ex}{$#1$}}}\raisebox{-.6ex}
{$\begin{array}{cccccc}
\cline{6-6}
& & &  & & \lr{1 } \\
\cline{4-6}
& & & \lr{1 } & \lr{3}& \lr{3}\\
\cline{1-6}
\lr{1 }& \lr{1 } & \lr{2 } & \lr{3 } & \lr{4 }& \lr{5 }\\
\cline{1-6}
\end{array}$}}$}\ \ \ , \ \ \
\text{${\def\lr#1{\multicolumn{1}{|@{\hspace{.6ex}}c@{\hspace{.6ex}}|}{\raisebox{-.3ex}{$#1$}}}\raisebox{-.6ex}
{$\begin{array}{cc}
& \\
\cline{2-2}
&\lr{2}\\
\cline{1-1}\cline{2-2}
\lr{2}& \lr{3}\\
\cline{1-1}\cline{2-2}
\end{array}$}}$}\ \ \  \right)
\end{split}
\end{equation*}}
\end{ex}

\begin{ex}{\rm
Let
\begin{equation*}
\begin{split}
& (S,T)=\left(\ \
\text{${\def\lr#1{\multicolumn{1}{|@{\hspace{.6ex}}c@{\hspace{.6ex}}|}{\raisebox{-.3ex}{$#1$}}}\raisebox{-.6ex}
{$\begin{array}{cccccccc}
\cline{5-8}
& & & &\lr{1 } & \lr{1 } & \lr{3}& \lr{3}\\
\cline{1-8}
\lr{1} & \lr{1} & \lr{1}& \lr{2} & \lr{2} & \lr{3 } & \lr{4 }& \lr{5 }\\
\cline{1-8}
\end{array}$}}$}\ \ \ , \ \ \
\text{${\def\lr#1{\multicolumn{1}{|@{\hspace{.6ex}}c@{\hspace{.6ex}}|}{\raisebox{-.3ex}{$#1$}}}\raisebox{-.6ex}
{$\begin{array}{cc}
\cline{2-2}
&\lr{2}\\
\cline{1-1}\cline{2-2}
\lr{2}& \lr{3}\\
\cline{1-1}\cline{2-2}
\end{array}$}}$}\ \ \  \right) \in {\mathcal{T}}^{\mf c} (2,1).
\end{split}
\end{equation*}
Then
\begin{equation*}
\begin{split}
&\sigma=(\
\cdots \ + \ - \   \ \cdot \ \, -    \ \,  +\ ), \\
&\td{\sigma}=(\
\cdots \ + \ - \   \ \cdot \ \ \,  \cdot    \ \  \, \cdot  \ \ \, ),
\end{split}
\end{equation*}
and
\begin{equation*}
\begin{split}
\td{E}_0(S,T)&=
\left(\ \
\text{${\def\lr#1{\multicolumn{1}{|@{\hspace{.6ex}}c@{\hspace{.6ex}}|}{\raisebox{-.3ex}{$#1$}}}\raisebox{-.6ex}
{$\begin{array}{cccccc}
\cline{3-6}
 & &\lr{1 } & \lr{1 } & \lr{3}& \lr{3}\\
\cline{1-6}
 \lr{1}& \lr{2} & \lr{2} & \lr{3 } & \lr{4 }& \lr{5 }\\
\cline{1-6}
\end{array}$}}$}\ \ \ , \ \ \
\text{${\def\lr#1{\multicolumn{1}{|@{\hspace{.6ex}}c@{\hspace{.6ex}}|}{\raisebox{-.3ex}{$#1$}}}\raisebox{-.6ex}
{$\begin{array}{cc}
\cline{2-2}
&\lr{2}\\
\cline{1-1}\cline{2-2}
\lr{2}& \lr{3}\\
\cline{1-1}\cline{2-2}
\end{array}$}}$}\ \ \  \right),\\
\td{F}_0(S,T)&=
\left(\ \
\text{${\def\lr#1{\multicolumn{1}{|@{\hspace{.6ex}}c@{\hspace{.6ex}}|}{\raisebox{-.3ex}{$#1$}}}\raisebox{-.6ex}
{$\begin{array}{cccccccccc}
\cline{7-10}
& & & & & &\lr{1 } & \lr{1 } & \lr{3}& \lr{3}\\
\cline{1-10}
\lr{1}& \lr{1} & \lr{1} & \lr{1} & \lr{1}& \lr{2} & \lr{2} & \lr{3 } & \lr{4 }& \lr{5 }\\
\cline{1-10}
\end{array}$}}$}\ \ \ , \ \ \
\text{${\def\lr#1{\multicolumn{1}{|@{\hspace{.6ex}}c@{\hspace{.6ex}}|}{\raisebox{-.3ex}{$#1$}}}\raisebox{-.6ex}
{$\begin{array}{cc}
\cline{2-2}
&\lr{2}\\
\cline{1-1}\cline{2-2}
\lr{2}& \lr{3}\\
\cline{1-1}\cline{2-2}
\end{array}$}}$}\ \ \  \right) .
\end{split}
\end{equation*}

 }
\end{ex}

\vskip 2mm

Now, for  $(S,T)\in {\mathcal{T}}^{\mf x}(\lambda)$, we  define
\begin{equation}
\begin{split}
{\rm wt}^{\mf x}(S,T)&=
\sum_{i\geq 1}\left(m_i(S)+m_i(T)\right)\widehat{\epsilon}_i ,\\
{\varepsilon}^{\mf x}_i(S,T)&=\max\{\,k\,|\,\td{E}_i^k(S,T)\neq {\bf 0}\,\},\\
{\varphi}^{\mf x}_i(S,T)&=\langle {\rm wt}^{\mf x}(S,T), \widehat{h}_i \rangle+ {\varepsilon}^{\mf x}_i(S,T),
\end{split}
\end{equation}
for $i\in \Z_{\geq 0}$. Put
\begin{equation}
{\mathcal{H}}_{\lambda}=(\emptyset, H_\lambda),
\end{equation}
where $\emptyset$ is the empty tableau.

\begin{prop}
For $\lambda\in\cP$, ${\mathcal{T}}^{\mf x}(\lambda)$ is an ${\mf x}_\infty$-crystal with respect to ${\rm wt}^{\mf x}$, ${\varepsilon}_i^{\mf x}$, ${\varphi}_i^{\mf x}$, $\td{E}_i$, $\td{F}_i$ $(i\in\Z_{\geq 0})$.
Also we have
$${\mathcal{T}}^{\mf x}(\lambda)=\{\,\td{F}_{i_1}\cdots\td{F}_{i_r}{\mathcal{H}}_{\lambda}\,|\,r\geq
0,\  i_1,\ldots,i_r \in\Z_{\geq 0}\,\}\setminus\{{\bf 0}\}.$$ In particular, ${\mathcal{T}}^{\mf x}(\lambda)$ is a connected ${\mf x}_\infty$-crystal with a unique highest weight element $\mathcal{H}_\lambda$.
\end{prop}
\pf By definition, ${\mathcal{T}}^{\mf x}(\lambda)$ is an ${\mf x}_\infty$-crystal. Let us show that given $(S,T)\in{\mathcal{T}}^{\mf x}(\lambda)$, there exist  $i_1,\ldots,i_r \in\Z_{\geq 0}$ such that $\td{E}_{i_1}\cdots\td{E}_{i_r}(S,T)={\mathcal{H}}_{\lambda}$. We will use induction on $| {\rm sh}(S) |$.

If $|{\rm sh}(S)|=0$, then $S$ is the empty tableau. Since $SST_{\Z_{>0}}(\lambda)$ is a highest weight ${\mf l}_{\infty}$-crystal, there exist  $i_1,\ldots,i_r \in\Z_{>0}$ such that $\td{E}_{i_1}\cdots\td{E}_{i_r}(\emptyset,T)={H}_{\lambda}$, and hence $\td{E}_{i_1}\cdots\td{E}_{i_r}(S,T)={\mathcal{H}}_{\lambda}$.

Suppose that $|{\rm sh}(S)|>0$. We may assume that $(S,T)$ is an ${\mf l}_\infty$-highest weight element. Then by tensor product rule of crystals, we have $S=H_\tau$. Note that there is at least one removable
${\def\lr#1{\multicolumn{1}{|@{\hspace{.6ex}}c@{\hspace{.6ex}}|}{\raisebox{-.3ex}{$#1$}}}\raisebox{-.6ex}
{$\begin{array}[b]{c}
\cline{1-1}
\lr{1}\\
\cline{1-1}
\end{array}$}}$
 or
${\def\lr#1{\multicolumn{1}{|@{\hspace{.6ex}}c@{\hspace{.6ex}}|}{\raisebox{-.3ex}{$#1$}}}\raisebox{-.6ex}
{$\begin{array}[b]{cc}
\cline{1-1}\cline{2-2}
\lr{1}&\lr{1}\\
\cline{1-1}\cline{2-2}
\end{array}$}}$
 in $S$, and hence $\td{E}_0S\neq {\bf 0}$ with $|{\rm sh}(\td{E}_0S)|=|{\rm sh}(S)|-\epsilon$. By induction hypothesis,  there exist   $i_1,\ldots,i_r \in\Z_{\geq 0}$ such that $\td{E}_{i_1}\cdots\td{E}_{i_r}\td{E}_0(S,T)={\mathcal{H}}_{\lambda}$.
\qed

\begin{rem}{\rm
Note that
\begin{equation*}
{\rm ch}\mathcal{T}^{\mf x}(\lambda)=
s_{\lambda}({\bf x})\sum_{\tau\in \cP_{\mf x}}s_\tau({\bf x})=\sum_{\sigma\in \cP}\sum_{\tau\in \cP_{\mf x}}c^\sigma_{\tau \lambda}s_\sigma({\bf x}),
\end{equation*}
where $x_i=e^{\widehat{\epsilon_i}}\in\Z[P_{\mf x}]$ for $i\in\Z_{>0}$.
By the Littlewood identity (\ref{Littlewood identity}), we have
\begin{equation*}
{\rm ch}\mathcal{T}^{\mf x}(\lambda)=\dfrac{s_\lambda({\bf x})}{\prod_i(1-x_i^\epsilon)\prod_{i<j}(1-x_ix_j)}.
\end{equation*}
}
\end{rem}
\vskip 2mm

\subsection{Crystals of integrable highest weight modules}
For ${\bf w}=w_1\ldots w_r\in\W_{\mathbb{Z}_{>0}}$, let us say that $\w$ is weakly decreasing if $w_1\geq \ldots\geq w_r$.

\begin{df}\label{Def of Delta}{\rm
  For $U\in {\rm \bf LR}^\lambda_{\mu\nu}$ with $\mu\in \cP_{\mf x}$, we define
\begin{equation*}
\begin{split}
&\Delta_{\mf x}(U)  =\max
\left\{\, \frac{\mu_k+2\ell(\w)}{\epsilon}\, \Bigg\vert \,
\begin{array}{l}
\text{{\rm (1)} $\w=w_1\ldots w_r$ is a subword of $w(U)$,}       \\
\text{{\rm (2)} $\w$ is weakly decreasing,  }       \\
\text{{\rm (3)} $w_1$ is   placed in the $k$-th row of $\lambda$.}       \\
\end{array}\, \right\}. \\
\end{split}
\end{equation*} }
\end{df}

\begin{rem}\label{remark on Delta}{\rm
For $U\in {\rm \bf LR}^\lambda_{\mu\nu}$ let $\w=w_1\ldots w_r$ be a weakly decreasing subword of $w(U)$. Let $w_k$ be the entry of $U$ at the position $(i_k,j_k)$ , that is, the $i_k$-th row and the $j_k$-th column. If $w(U)$ is a row word, then $i_1\geq \ldots \geq i_r$. Also we have $j_1< \ldots < j_r$ since $w$ is weakly decreasing. If $w(U)$ is a column word, then  $j_1< \ldots < j_r$. Suppose that $i_k<i_{k+1}$ for some $k$. Then we may replace $w_{k+1}$ with the entry of $U$ at $(i_k,j_{k+1})$ to obtain a new weakly decreasing subword of the same length. Repeating this process, we may also assume that $i_1\geq \ldots \geq i_r$. Hence we may always assume that  $i_1\geq \ldots \geq i_r$ and $j_1< \ldots < j_r$. Note that a subword of $w(U)$ which gives a maximal value $\Delta_{\mf x}(U)$ is not necessarily unique.
}
\end{rem}

\begin{ex}{\rm Consider
\begin{center}
 $U \ \ =$ \ \
${\def\lr#1{\multicolumn{1}{|@{\hspace{.6ex}}c@{\hspace{.6ex}}|}{\raisebox{-.3ex}{$#1$}}}\raisebox{-.6ex}
{$\begin{array}{ccccc}
\cline{5-5}
& & & & \lr{1 } \\
\cline{4-5}
& &  & \lr{1}& \lr{2}\\
\cline{2-5}
  & \lr{{\bf 2}} & \lr{{\bf 3}} & \cdot & \cdot \\
\cline{1-3}
\lr{{\bf 1}} & \lr{4} & \cdot & \cdot &  \cdot\\
\cline{1-2}
\end{array}$}}$
\ \
$\in\ \  { \rm \bf LR}^{(5,4,2,1)}_{(3,2)\ (3,2,1,1)}$.
\end{center}\vskip 2mm
The  bold faced letters form a weakly decreasing subword of $w(U)$, which gives a maximal value $2+2\cdot 3=8=\Delta_{\mf b}(U)$.

And for
\begin{center}
 $U' \ \ =$\ \ \
${\def\lr#1{\multicolumn{1}{|@{\hspace{.6ex}}c@{\hspace{.6ex}}|}{\raisebox{-.3ex}{$#1$}}}\raisebox{-.6ex}
{$\begin{array}{cccccc}
\cline{5-6}
& & &  & \lr{1}& \lr{1}\\
\cline{4-6}
&  &  & \lr{2} & \cdot & \cdot \\
\cline{1-2}\cline{4-4}
\lr{{\bf 2}} & \lr{{\bf 3}} & \cdot & \cdot &  \cdot & \cdot\\
\cline{1-2}
\end{array}$}}$
\ \
$\in\ \  { \rm \bf LR}^{(6,3,2)}_{(4,2)\ (2,2,1)}$,
\end{center}\vskip 2mm
we have $\Delta_{\mf c}(U')=(4+2\cdot 2)/2=4$.
}
\end{ex}

\begin{df}{\rm
For  $(\lambda,n)\in \cP({\mf x})$, we define
\begin{equation*}
\begin{split}
&{\mathcal{T}}^{\mf x}(\lambda,n)=  \left\{\, (S,T) \ \Big\vert
\begin{array}{l}
\text{$(S,T)\in {\mathcal{T}}^{\mf x}(\lambda)$}, \
\text{$\Delta_{\mf x}\left((T\rightarrow S)_R\right)\leq  n$.}
\end{array} \right\}.
\end{split}
\end{equation*}}
\end{df}


\begin{figure}
\includegraphics{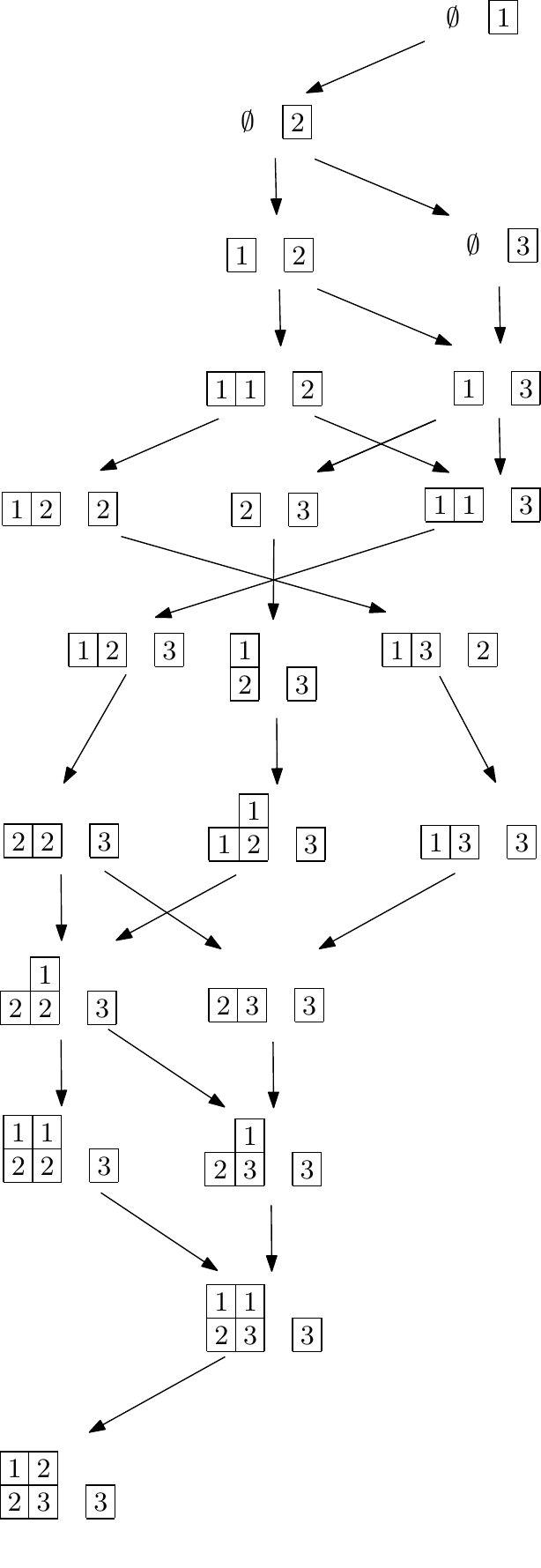}
\caption{The subgraph of $\mathcal{T}^{\mf b}((1),2)$ obtained by applying $\td{F}_i$'s ($i=0,1,2$) to the highest weight vector $\mathcal{H}_{(1)}$ of weight $\Lambda^{\mf b}_1$.}
\end{figure}

\begin{figure}
\includegraphics{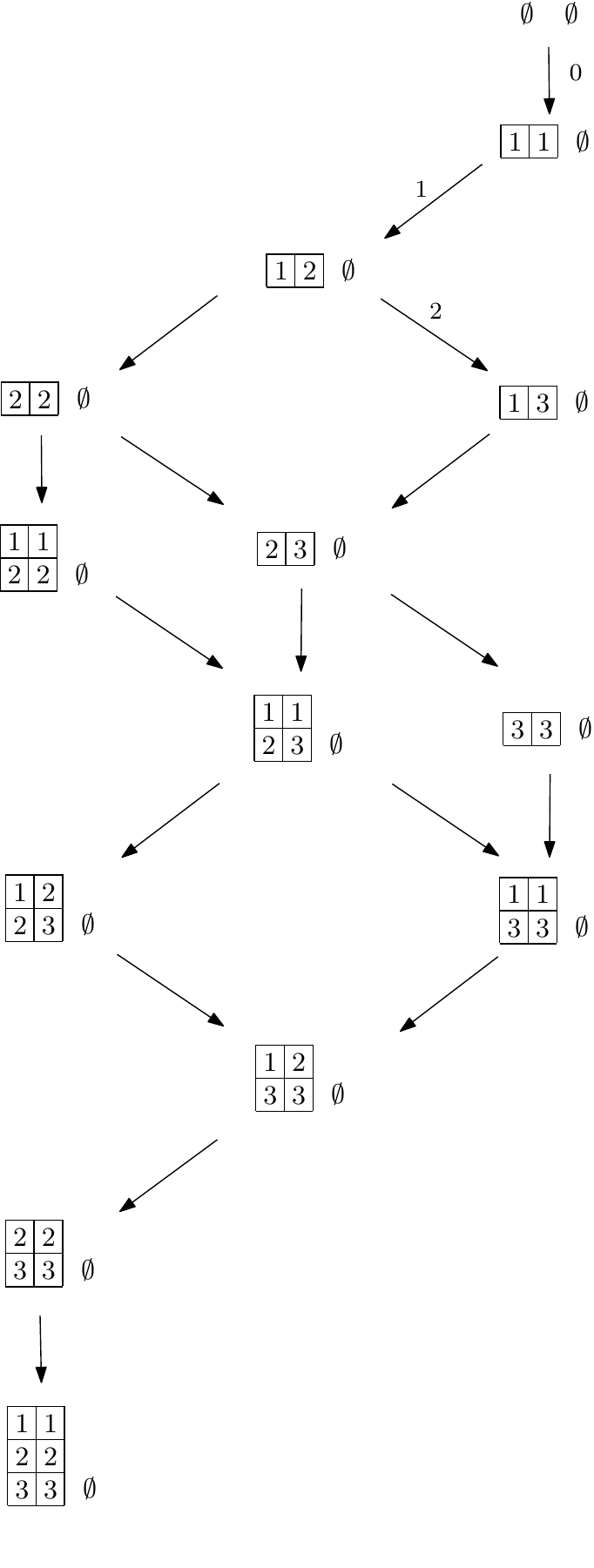}
\caption{The subgraph of $\mathcal{T}^{\mf c}((0),1)$ obtained by applying $\td{F}_i$'s ($i=0,1,2$) to the highest weight vector $\mathcal{H}_{(0)}$ of weight $\Lambda^{\mf c}_0$.}
\end{figure}


Note that $\Delta_{\mf x}(\mathcal{H}_\lambda)= 2\lambda_1/\epsilon$, and hence $\mathcal{H}_\lambda \in {\mathcal{T}}^{\mf x}(\lambda,n)$.
For $(S,T)\in {\mathcal{T}}^{\mf x}(\lambda,n)$ we define
\begin{equation}
{\rm wt}^{\mf x}_n(S,T)={\rm wt}^{\mf x}(S,T)+n\La^{\mf x}_0.
\end{equation}
Then we may view ${\mathcal{T}}^{\mf x}(\lambda,n)$ as an $\mf{x}_\infty$-crystal with respect to ${\rm wt}^{\mf x}_n$, $\varepsilon_i^{\mf x}$, $\varphi^{\mf x}_i$ and $\td{E}_i, \td{F}_i$ ($i\in\Z_{\geq 0}$), which is a subgraph of ${\mathcal{T}}^{\mf x}(\lambda)\otimes T_{n\La^{\mf x}_0}$.  \vskip 2mm

Now we are in a position to state the main result of this paper, which gives a new combinatorial realization of the crystal $\B({\mf x}_\infty,\La^{\mf x}(\lambda,n))$. The proof will be given in Section  \ref{Proof for BC}.

\begin{thm}\label{main result}
For $(\lambda,n)\in  \cP({\mf x})$, we have
$$ {\mathcal{T}}^{\mf x}(\lambda,n)\simeq \B({\mf x}_\infty,\La^{\mf x}(\lambda,n)).$$
\end{thm}

\begin{rem}\label{epsilon star}{\rm
Let $\B(\infty)$ be an ${\mf x}_\infty$-crystal associated with the negative part of $U_q({\mf x}_\infty)$ with a highest weight element $u_\infty$. There exists a weight preserving involution $\ast$ on $\B(\infty)$, which  gives another crystal structure on $\B(\infty)$ with respect to $\td{E}_i^\ast=\ast\circ \td{E}_i \circ \ast$ and $\td{F}_i^\ast=\ast\circ \td{F}_i \circ \ast$.
Note that there exists a crystal embedding $ \B({\mf x}_\infty,\La^{\mf x}(\lambda,n)) \hookrightarrow \B(\infty)\otimes T_{\La^{\mf x}(\lambda,n)}$ sending $u_{\La^{\mf x}(\lambda,n)}$ to $u_\infty\otimes t_{\La^{\mf x}(\lambda,n)}$, and   its image is
\begin{equation}\label{star epsilon}
\{\,b\otimes  t_{\La^{\mf x}(\lambda,n)}\,|\, \varepsilon^*_i(b)\leq  \langle \La^{\mf x}(\lambda,n),\widehat{h}_i \rangle \ (i\in \Z_{\geq 0})\,\},
\end{equation}
where $\varepsilon^*_i(b)=\max\{\,k\,|\,(\td{E}_i^\ast)^kb\neq{\bf 0}\,\}$ (see \cite{Kas94}).

Since $ {\mathcal{T}}^{\mf x}(\lambda,n)\subset   {\mathcal{T}}^{\mf x}(\lambda,n+1)$ for $n\geq 1$, we may regard $\bigcup_n  {\mathcal{T}}^{\mf x}(\lambda,n) = {\mathcal{T}}^{\mf x}(\lambda)\otimes T_{\La^{\mf x}(\lambda,n)}\subset \B(\infty)\otimes T_{\La^{\mf x}(\lambda,n)}$.
Theorem \ref{main result} combined with (\ref{star epsilon}) implies that $$\Delta_{\mf x}\left((T\rightarrow S)_R\right)-(2\lambda_1/\epsilon) =\varepsilon_0^\ast(S,T)$$ for $(S,T)\in {\mathcal{T}}^{\mf x}(\lambda)$. Hence $\Delta_{\mf x}$ can be viewed as a combinatorial realization of $\varepsilon_0^\ast$.}
\end{rem}

\section{Decomposition rule}\label{Decomposition rule}
\subsection{Generalized LR rule}
Let us introduce another statistic on LR tableaux.
\begin{df}\label{LR coefficient}{\rm \mbox{}
\begin{itemize}
\item[(1)] For $U\in {\rm \bf LR}^\lambda_{\mu\nu}$ with $\nu\in\cP_{\mf x}$, we define
\begin{equation*}
\begin{split}
\nabla_{\mf x}(U)
&=
\begin{cases}
\text{the number of removable
${\def\lr#1{\multicolumn{1}{|@{\hspace{.6ex}}c@{\hspace{.6ex}}|}{\raisebox{-.0ex}{$#1$}}}\raisebox{-.6ex}
{$\begin{array}[b]{c}
\cline{1-1}
\lr{1}\\
\cline{1-1}
\end{array}$}}$
's in $U^\sharp$}, & \text{if ${\mf x}={\mf b}$},\\
\text{the number of removable
${\def\lr#1{\multicolumn{1}{|@{\hspace{.6ex}}c@{\hspace{.6ex}}|}{\raisebox{-.0ex}{$#1$}}}\raisebox{-.6ex}
{$\begin{array}[b]{cc}
\cline{1-1}\cline{2-2}
\lr{1}&\lr{1}\\
\cline{1-1}\cline{2-2}
\end{array}$}}$
's in $U^\sharp$}, & \text{if ${\mf x}={\mf c}$}.
\end{cases}
\end{split}
\end{equation*}

\item[(2)] For $(\lambda,m+n), (\mu,m), (\nu,n)\in \cP({\mf x})$, we define
\begin{equation*}
\begin{split}
&{\rm \bf LR}^{(\lambda,m+n)}_{(\mu,m) (\nu,n)}({\mf x})  =\\ &\left\{\, (U,V) \, \Bigg\vert
\begin{array}{l}
\text{(1) $(U,V) \in {{\rm \bf LR}}^{\eta}_{\mu \tau}\times {{\rm \bf LR}}^{\lambda}_{\eta \nu}$ for some $\tau\in \cP_{\mf x},\eta\in\cP$,}\\
\text{(2) $\nabla_{\mf x}(U)\leq  m-(2\mu_1/\epsilon)$\, and \,  $\Delta_{\mf x}\left(( V^\sharp \rightarrow  U^\sharp)_R\right)\leq  n$.}
\end{array} \right\},
\end{split}
\end{equation*}
and let
\begin{equation*}
c^{(\lambda,m+n)}_{(\mu,m) (\nu,n)}({\mf x})=\Big| {\rm \bf LR}^{(\lambda,m+n)}_{(\mu,m) (\nu,n)}({\mf x})  \Big|.
\end{equation*}
\end{itemize}
}
\end{df}
\vskip 2mm
Then we have the following generalized LR rule of ${\mf x}_\infty$-crystals.

\begin{thm}\label{LR rule}
For  $(\mu,m), (\nu,n)\in \cP({\mf x})$, we have {\rm
\begin{equation*}
\begin{split}
\B({\mf x}_\infty,\La^{\mf x}(\mu,m))&\otimes \B({\mf x}_\infty,\La^{\mf x}(\nu,n))\simeq  \\ & \bigsqcup_{(\lambda,m+n) \in\cP({\mf x})}\B({\mf x}_\infty,\La^{\mf x}(\lambda,m+n))\times {{\rm \bf LR}}^{(\lambda,m+n)}_{(\mu,m) (\nu,n)}({\mf x}),
\end{split}
\end{equation*}}
where $\td{E}_i$ and $\td{F}_i$ for $i\in\Z_{\geq 0}$ act on the first component of the right-hand side.
\end{thm}
\pf Let ${\bf m}_1\otimes {\bf m}_2\in \B({\mf x}_\infty,\La^{\mf x}(\mu,m))\otimes \B({\mf x}_\infty,\La^{\mf x}(\nu,n))$ be given. Suppose that ${\bf m}_1\otimes {\bf m}_2$ is a highest weight element with highest weight $\Lambda^{\mf x}(\lambda,m+n)$ for some $(\lambda,m+n)\in\cP(\mf x)$.
Then by tensor product rule of crystals, it is equivalent to saying that ${\bf m}_1=u_{\Lambda^{\mf x}(\mu,m)}$ and for $i\in\Z_{\geq 0}$
\begin{equation}\label{hwvector}
{\varepsilon}^{\mf x}_i({\bf m}_2)\leq \langle \Lambda^{\mf x}(\mu,m), \widehat{h}_i \rangle.
\end{equation}
By Theorem \ref{main result}, we may identify ${\bf m}_1=(\emptyset, H_{\mu})$ and ${\bf m}_2=(S,T)\in SST_{\Z_{>0}}(\tau)\times SST_{\Z_{>0}}(\nu)$  for some $\tau\in\cP_{\mf x}$. Since ${\varepsilon}^{\mf x}_i({\bf m}_2)\leq \langle \Lambda^{\mf x}(\mu,m), \widehat{h}_i \rangle$ for $i\in \Z_{>0}$, equivalently, $\td{E}_i({\bf m}_1\otimes {\bf m}_2)={\bf 0}$ for $i\in\Z_{>0}$ (that is, a highest weight element in an $\frak{l}_\infty$-crystal), we have $(S\rightarrow H_\mu)=H_\eta$ and $(T\rightarrow H_\eta)=H_\lambda$  for some $\eta\in \cP$. Hence $S\in {{\rm LR}}^{\eta}_{\mu\,\tau}$ and $T\in {\rm LR}^{\lambda}_{\eta \nu}$.

Let $U \in {{\rm \bf LR}}^{\eta}_{\mu\,\tau}$ and $V\in {\bf LR}^{\lambda}_{\eta \nu}$ be such that $U^\sharp=S$ and  $V^\sharp=T$. Also by (\ref{hwvector}), we have $${\varepsilon}^{\mf x}_0({\bf m}_2)={\varepsilon}^{\mf x}_0(S,T)=\nabla_{\mf x}(U)\leq  \langle \Lambda^{\mf x}(\mu,m), \widehat{h}_0 \rangle=m-(2\mu_1/\epsilon).$$ Finally, we have
$\Delta_{\mf x}(( V^\sharp \rightarrow  U^\sharp)_R)\leq n$ since $(S,T)\in {\mathcal{T}}^{\mf x}(\nu,n)$. Hence $(U,V)\in {{\rm \bf LR}}^{(\lambda,m+n)}_{(\mu,m) (\nu,n)}({\mf x})$.

Since ${\bf m}_1\otimes {\bf m}_2$ is uniquely determined by $(U,V)$ and the correspondence ${\bf m}_1\otimes {\bf m}_2\mapsto (U,V)$ is reversible, the connected component isomorphic to $\B({\mf x}_\infty,\La^{\mf x}(\lambda,m+n))$ is parametrized by ${{\rm \bf LR}}^{(\lambda,m+n)}_{(\mu,m) (\nu,n)}({\mf x})$.  This completes the proof. \qed\vskip 2mm

\begin{cor}\label{LR rule-character}
For  $(\mu,m), (\nu,n)\in \cP({\mf x})$, we have {\rm
\begin{equation*}
\begin{split}
{\rm ch}\B({\mf x}_\infty,\La^{\mf x}(\mu,m))\ &{\rm ch}\B({\mf x}_\infty,\La^{\mf x}(\nu,n))=  \\ & \sum_{(\lambda,m+n) \in\cP({\mf x})}c^{(\lambda,m+n)}_{(\mu,m) (\nu,n)}({\mf x})\ {\rm ch} \B({\mf x}_\infty,\La^{\mf x}(\lambda,m+n)).
\end{split}
\end{equation*}}
\end{cor}

\begin{ex}{\rm
When $\mu=\nu=(0)$, it is straightforward to check that
$${\rm \bf LR}^{(\lambda,m+n)}_{((0),m) ((0),n)}({\mf x})  = \{\,(H_\lambda,\emptyset)\,\}$$ for $(\lambda,m+n)\in\cP({\mf x})$,
which implies that
\begin{equation*}
\begin{split}
\B({\mf x}_\infty,\La^{\mf x}((0),m))&\otimes \B({\mf x}_\infty,\La^{\mf x}((0),n))\simeq  \bigsqcup_{(\lambda,m+n) \in\cP({\mf x})}\B({\mf x}_\infty,\La^{\mf x}(\lambda,m+n)).
\end{split}
\end{equation*}
This recovers a multiplicity free decomposition of the product of rectangular shaped characters  for $\mf{sp}(2k)$ and $\mf{so}(2k+1)$  by Okada \cite[Theorem 2.5 (1) and (2)]{O} (see Section \ref{finite rank crystals}). Also, a Krattenthaler's generalization to nearly rectangular shaped characters \cite[Theorem 3]{Kr} can be deduced from Theorem \ref{LR rule} without difficulty.}
\end{ex}

The  generalized LR coefficients $c^{(\lambda,m+n)}_{(\mu,m) (\nu,n)}({\mf x})$ has the following stable limit.
\begin{cor}\label{stable formula-1}
Let $(\lambda,m+n), (\mu,m), (\nu,n)\in \cP({\mf x})$ be given. If
$(\lambda,m), (\lambda,n)\in\cP({\mf x})$,
then
$$c^{(\lambda,m+n)}_{(\mu,m) (\nu,n)}({\mf x})=\sum_{\eta\in\cP}\sum_{\tau\in \cP_{\mf x}}c^\eta_{\mu \tau}c^{\lambda}_{\eta \nu}.$$
\end{cor}
\pf Suppose that $(U,V) \in {{\rm \bf LR}}^{\eta}_{\mu \tau}\times {{\rm \bf LR}}^{\lambda}_{\eta \nu}$ is given for some $\tau\in \cP_{\mf x}$ and $\eta\in\cP$.
Note that $\tau_1+\mu_1\leq \eta_1\leq \lambda_1$. By hypothesis, we have
$\tau_1+2\mu_1\leq 2\lambda_1 \leq \epsilon m,$ which implies that
$$\nabla_{\mf x}(U)\leq \tau_1/\epsilon\leq m-(2\mu_1/\epsilon).$$
Similarly, we have
$\tau_1+2\nu_1\leq 2\lambda_1 \leq \epsilon n,$
which implies that
$$\Delta_{\mf x}\left(( V^\sharp \rightarrow  U^\sharp)_R\right)\leq \frac{\tau_1+2\nu_1}{\epsilon} \leq  n.$$ Therefore, $(U,V)\in {\rm \bf LR}^{(\lambda,m+n)}_{(\mu,m) (\nu,n)}({\mf x})$.  We have $${\rm \bf LR}^{(\lambda,m+n)}_{(\mu,m) (\nu,n)}({\mf x})=\bigsqcup_{\eta\in\cP}\bigsqcup_{\tau\in\cP_{\mf x}}{{\rm \bf LR}}^{\eta}_{\mu\, \tau}\times {{\rm \bf LR}}^{\lambda}_{\eta \nu}.$$  \qed

\begin{rem}\label{classical stable formula}{\rm
The stable limit formula in Corollary \ref{stable formula-1} also recovers and extends the known result in the following sense.
Let $\ell=\min\{m,n\}$. Suppose that $2\lambda_1,2\mu_1,2\nu_1\leq \epsilon\ell$. It is not difficult to check that
$$c^{(\lambda,m+n)}_{(\mu,m) (\nu,n)}({\mf x})=\sum_{\eta\in\cP}\sum_{\tau\in \cP_{\mf x}}c^\eta_{\mu \tau}c^{\lambda}_{\eta \nu}=\sum_{\eta\in\cP}\sum_{\tau\in \cP_{\mf x}}c^\gamma_{\mu \nu}c^{\lambda}_{\gamma \tau}.$$
If we apply the reciprocity law via  the $({\mf c}_\infty,{\rm Sp}(2k))$-Howe duality on the level $k$ fermionic Fock space \cite[Theorem 3.4]{Wa}, then it follows that $c^{(\lambda,m+n)}_{(\mu,m) (\nu,n)}({\mf c})$ coincides with the stable branching coefficient with respect to ${\rm Sp}(2m)\times {\rm Sp}(2n)\subset {\rm Sp}(2m+2n)$ given by Koike and Terada \cite[Corollary 2.6]{KoT90}. Also by the $({\mf b}_\infty,{\rm Pin}(2k))$-Howe duality \cite[Theorem 3.3]{Wa}, $c^{(\lambda,m+n)}_{(\mu,m) (\nu,n)}({\mf b})$ gives a stable branching coefficient with respect to ${\rm Pin}(2m)\times {\rm Pin}(2n)\subset {\rm Pin}(2m+2n)$, which seems to be new (cf.\cite[Section 7]{Wa}).}
\end{rem}

\subsection{Branching rule as an ${\mf l}_\infty$-crystal}
Let us describe the decomposition of $\B({\mf x}_\infty,\La^{\mf x}(\lambda,n))$ as an ${\mf l}_\infty$-crystal.
\begin{df}\label{branching coefficient}{\rm
For $(\lambda,n) \in \cP({\mf x})$ and $\sigma\in\cP$, we define
\begin{equation*}
{\rm \bf LR}^{\, \sigma}_{(\lambda,n)}({\mf x})=\{\,U\,|\,\text{$U\in {\rm \bf LR}^{\sigma}_{\tau \lambda}$ for some $\tau\in \cP_{\mf x}$ and $\Delta_{\mf x}(U)\leq n$}\,\},
\end{equation*}
and let
$$c^\sigma_{(\lambda,n)}({\mf x})=\Bigl\vert{\rm \bf LR}^{\, \sigma}_{(\lambda,n)}({\mf x})\Bigr\vert.$$
}
\end{df}\vskip 2mm

For $\sigma\in \cP$, let $\B({\mf l}_\infty,\omega_\sigma)$ be the ${\mf l}_\infty$-crystal associated to the highest weight $U_q({\mf l}_\infty)$-module with highest weight $\omega_\sigma=\sum_{i\geq 1}\sigma_i\widehat{\epsilon}_i$. Note that $\B({\mf l}_\infty,\omega_\sigma)\simeq SST_{\Z_{>0}}(\sigma)$.
Then we have the following decomposition.

\begin{thm}\label{branching} For $(\lambda,n)\in \cP({\mf x})$, we have
{\rm
\begin{equation*}
\B({\mf x}_\infty,\La^{\mf x}(\lambda,n))\simeq \bigsqcup_{\sigma\in\cP}
\B({\mf l}_\infty,\omega_\sigma)\times{\rm \bf LR}^{\,\sigma}_{(\lambda,n)}({\mf x})
\end{equation*}}
\hskip -2mm as an $\mathfrak{l}_\infty$-crystal, where $\td{E}_i$ and $\td{F}_i$  $(i\in \Z_{>0})$ act on the first component of the right-hand side.
\end{thm}
\pf Let $(S,T)\in \B({\mf x}_\infty,\La^{\mf x}(\lambda,n))$ be such that $\td{E}_i(S,T)={\bf 0}$ for $i\in \Z_{>0}$. By tensor product of crystals, we have $S=H_{\tau}$ for some $\tau\in \cP_{\mf x}$ and $(T\rightarrow H_\tau)=H_{\sigma}$ for some $\sigma\in \cP$. Hence $U=(T\rightarrow H_\tau)_R\in {\bf LR}^\sigma_{\tau \lambda}$ with $\Delta_{\mf x}(U)\leq n$.
Since the correspondence $(S,T)\mapsto U$ is reversible, we obtain the above decomposition. \qed

\begin{cor}\label{branching-character}
For $(\lambda,n) \in \cP({\mf x})$, we have {\rm
$${\rm ch} \B({\mf x}_\infty,\La^{\mf x}(\lambda,n))=e^{n {\Lambda}^{\mf x}_0}\sum_{\sigma\in\cP} c^\sigma_{(\lambda,n)}({\mf x}) s_{\sigma}({\bf x}).$$
}
\end{cor}

The branching coefficient $c^\sigma_{(\lambda,n)}({\mf x})$ has the following stable limit.
\begin{cor}\label{stable formula-2}
Let $(\lambda,n) \in \cP({\mf x})$ and $\sigma\in\cP$ be given. If $(\sigma,n)\in\cP({\mf x})$, then
$$c^\sigma_{(\lambda,n)}({\mf x})=\sum_{\tau\in\cP_{\mf x}}c^\sigma_{\tau \lambda}.$$
\end{cor}
\pf Suppose that $U\in {\rm \bf LR}^{\sigma}_{\tau \lambda}$ is given for some $\tau\in \cP_{\mf x}$. Since $\tau_1+\lambda_1\leq \sigma_1$, we have
$$
\Delta_{\mf x}(U)\leq \frac{\tau_1+2\lambda_1}{\epsilon}\leq 2\sigma_1/\epsilon\leq n.
$$
which implies that $U\in  {\rm \bf LR}^{\,\sigma}_{(\lambda,n)}({\mf x})$. Hence, we have ${\rm \bf LR}^{\,\sigma}_{(\lambda,n)}({\mf x})=\bigsqcup_{\tau\in\cP_{\mf x}}{\rm \bf LR}^{\sigma}_{\tau \lambda}$.
\qed
\vskip 2mm

\begin{rem}\label{classical stable formula'}{\rm
When ${\mf x}={\mf c}$, the formula in Theorem \ref{branching} also recovers and extends the classical result in the following sense.
Consider the  the $({\mf c}_\infty,{\rm Sp}(2n))$-Howe duality on a  fermionic Fock space \cite[Theorem 3.4]{Wa}. We may consider $({\mf l}_\infty,{\rm GL}(2n))$-duality on the same Fock space. Then by the reciprocity law associated to the see-saw pairs of $({\mf c}_\infty,{\rm Sp}(2n))$ and $({\mf l}_\infty,{\rm GL}(2n))$, we conclude that $c^\sigma_{(\lambda,n)}({\mf c})$ in Corollary \ref{stable formula-2} gives a branching coefficient with respect to ${\rm Sp}(2n)\subset {\rm GL}(2n)$, which coincides with the Littlewood's restriction formula in the symplectic case \cite{Lw-1,Lw-2}.
}
\end{rem}

\subsection{Restriction to $\mf{so}(2k+1)$ and $\mf{sp}(2k)$-crystals}\label{finite rank crystals}
Fix  $k\geq 2$. For $(\lambda,n)\in\cP({\mf x})$ with $\ell(\lambda)\leq k-1$, let
\begin{equation}
\mathcal{T}^{\mf x}_k(\lambda,n)=\left\{\,(S,T)\,\Bigg\vert\,
\begin{array}{l}
\text{(1) $(S,T)\in\mathcal{T}^{\mf x}(\lambda,n)$}, \\
\text{(2)  the entries in $S$ and $T$ are no more than $k$.}
\end{array}
\,\right\}.
\end{equation}
We regard $\La^{\mf x}(\lambda,n)$ as a dominant integral weight for ${\mf x}_k$ and obtain the following realization of the crystal $\B({\mf x}_k,\La^{\mf x}(\lambda,n))$.
\begin{prop}\label{truncate to x_k}  For $(\lambda,n)\in\cP({\mf x})$ with $\ell(\lambda)\leq k-1$, we have
$${\mathcal{T}}_k^{\mf x}(\lambda,n)\simeq \B({\mf x}_k,\La^{\mf x}(\lambda,n)),$$
as an ${\mf x}_k$-crystal.
\end{prop}
\pf By restriction, we may regard ${\mathcal{T}}^{\mf x}_k(\lambda)$ as an ${\mf x}_k$-crystal. Since
$${\mathcal{T}}^{\mf x}_k(\lambda)=\{\,\td{F}_{i_1}\cdots\td{F}_{i_r}{\mathcal{H}}_{\lambda}\,|\,r\geq
0,\  i_1,\ldots,i_r \in \{0,1,\ldots,k-1\}\,\}\setminus\{{\bf 0}\},$$
we have ${\mathcal{T}}^{\mf x}_k(\lambda)\simeq \B({\mf x}_k,\La^{\mf x}(\lambda,n))$.
\qed\vskip 3mm

Note that $\B({\mf x}_k,\La^{\mf x}(\lambda,n))$ is the dual of a highest weight crystal of type $B_k$ and $C_k$ in usual convention (see for example \cite{KN}). Let us explain it briefly for the reader's convenience.

Let  $\alpha^{\dagger}_i=-\widehat{\alpha}_{k-i}$ for $1\leq i\leq k $.
If we put $\epsilon^{\dagger}_i=\widehat{\epsilon}_{k-i+1}$ for $1\leq i\leq k$, then
$\alpha_i^{\dagger}=\epsilon_i^{\dagger}-\epsilon^{\dagger}_{i+1}$ for $1\leq i\leq k-1$ and $\alpha^{\dagger}_k=\epsilon \epsilon^{\dagger}_k$. Let us denote by ${\mf x}_k^{\dagger}$ the Lie algebra isomorphic to ${\mf x}_k$ with the positive simple roots $\{\,\alpha^{\dagger}_i\,|\,1\leq i\leq k\,\}$. The associated  Dynkin diagrams are
\begin{center}\hskip 2cm
\setlength{\unitlength}{0.16in}
\begin{picture}(24,4)
\put(1,2){\makebox(0,0)[c]{${\mf b}^\dagger_{k}=\mf{so}(2k+1)\ :$}}

\put(5.6,2){\makebox(0,0)[c]{$\bigcirc$}}
\put(12.6,2){\makebox(0,0)[c]{$\bigcirc$}}
\put(10.4,2){\makebox(0,0)[c]{$\bigcirc$}}
\put(14.85,2){\makebox(0,0)[c]{$\bigcirc$}}
\put(6,2){\line(1,0){1.3}} \put(8.7,2){\line(1,0){1.3}} \put(10.82,2){\line(1,0){1.3}}
%
\put(13.7,2){\makebox(0,0)[c]{$\Longrightarrow$}}

\put(8,1.95){\makebox(0,0)[c]{$\cdots$}}
\put(5.6,1){\makebox(0,0)[c]{\tiny ${\alpha}^\dagger_1$}}
\put(12.7,1){\makebox(0,0)[c]{\tiny ${\alpha}^\dagger_{k-1}$}}
\put(10.4,1){\makebox(0,0)[c]{\tiny ${\alpha}^\dagger_{k-2}$}}
\put(15,1){\makebox(0,0)[c]{\tiny ${\alpha}^\dagger_k$}}

\put(20.5,2){\makebox(0,0)[c]{($\epsilon=1$),}}

\end{picture}
\end{center}

\begin{center}\hskip 2cm
\setlength{\unitlength}{0.16in}
\begin{picture}(24,4)
\put(1,2){\makebox(0,0)[c]{${\mf c}^\dagger_{k}=\mf{sp}(2k)\ \ \ \ \ \ :$}}

\put(5.6,2){\makebox(0,0)[c]{$\bigcirc$}}
\put(12.6,2){\makebox(0,0)[c]{$\bigcirc$}}
\put(10.4,2){\makebox(0,0)[c]{$\bigcirc$}}
\put(14.85,2){\makebox(0,0)[c]{$\bigcirc$}}
\put(6,2){\line(1,0){1.3}} \put(8.7,2){\line(1,0){1.3}} \put(10.82,2){\line(1,0){1.3}}
%
\put(13.7,2){\makebox(0,0)[c]{$\Longleftarrow$}}

\put(8,1.95){\makebox(0,0)[c]{$\cdots$}}
\put(5.6,1){\makebox(0,0)[c]{\tiny ${\alpha}^\dagger_1$}}
\put(12.7,1){\makebox(0,0)[c]{\tiny ${\alpha}^\dagger_{k-1}$}}
\put(10.4,1){\makebox(0,0)[c]{\tiny ${\alpha}^\dagger_{k-2}$}}
\put(15,1){\makebox(0,0)[c]{\tiny ${\alpha}^\dagger_k$}}

\put(20.5,2){\makebox(0,0)[c]{($\epsilon=2$).}}

\end{picture}
\end{center}\vskip 2mm

For $1\leq i\leq k$, let $\omega^{\mf x}_i$ be the $i$-th fundamental weight for ${\mf x}^\dagger_k$ with respect to $\{\,\alpha^{\dagger}_i\,|\,1\leq i\leq k\,\}$. Recall that $\omega^{\mf x}_i=\epsilon^\dagger_1+\cdots+\epsilon^\dagger_i$ for $1\leq i\leq k-1$ and $\omega^{\mf x}_k=(\epsilon/2)(\epsilon^\dagger_1+\cdots+\epsilon^\dagger_k)$. Then as an ${\mf x}^\dagger_k$-weight, we have
$\omega^{\mf x}_i=-\Lambda^{\mf x}_{k-i}$ for $1\leq i\leq k$, and $\Lambda^{\mf x}_i=0$ for $i\geq k$.
Hence, for a dominant integral ${\mf x}^\dagger_k$-weight $\omega=\sum_{1\leq i\leq k}a_i\omega^{\mf x}_i$ with $a_i\in\Z_{\geq 0}$, we have
\begin{equation}
\omega=-\Lambda^{\mf x}(\lambda(\omega), n(\omega)),
\end{equation}
where
\begin{equation}
\begin{split}
\lambda(\omega)&=
(a_1+\cdots+a_{k-1},a_1+\cdots+a_{k-2},\ldots,a_{1})
,\\
n(\omega)&=
\begin{cases}
2a_1+\cdots+2a_{k-1}+a_k, & \text{for ${\mf x}={\mf b}$}, \\
a_1+\cdots+a_k, & \text{for ${\mf x}={\mf c}$}.
\end{cases}
\end{split}
\end{equation}
Note that $\lambda(\omega)\in \cP$ with $\ell(\lambda(\omega))\leq k-1$ and $(\lambda(\omega),n(\omega))\in\cP({\mf x})$.
Then the map sending $\omega$ to $(\lambda(\omega),n(\omega))$ is a bijection from the set of dominant integral weights for ${\mf x}_k^\dagger$ to the set of $(\lambda,n)\in \cP({\mf x})$ with $\ell(\lambda)\leq k-1$.

Now, let $\B({\mf x}_k^{\dagger},\omega)$ be the crystal associated to the highest weight  $U_q({\mf x}_k^{\dagger})$-module with highest weight $\omega$. Then
$\B({\mf x}_k^{\dagger},\omega)$ is obtained from $\B({\mf x}_k,\Lambda^{\mf x}(\lambda(\omega), n(\omega)))$ by replacing an $i$-arrow with a $(k-i)$-arrow in opposite direction and negating the weight of each vertex modulo the weight lattice of ${\mf x}^\dagger_k$.

\begin{thm}\label{LR in opposite Dynkin}
We have
\begin{equation*}
\begin{split}
\B({\mf x}^\dagger_k,\omega') &\otimes \B({\mf x}^\dagger_k,\omega'')\simeq  \ \bigsqcup_{\omega}\B({\mf x}^\dagger_k,\omega)\times {{\rm \bf LR}}^{\omega}_{\omega'\, \omega''}({\mf x}_k^\dagger),
\end{split}
\end{equation*}
where $\omega$, $\omega'$ and $\omega''$ are dominant integral  ${\mf x}_k^\dagger$-weights and
{\rm \begin{equation*}
\begin{split}
& {{\rm \bf LR}}^{\omega}_{\omega'\, \omega''}({\mf x}_k^\dagger) =\\ &\left\{\, (U,V) \, \Bigg\vert
\begin{array}{l}
\text{(1) $(U,V) \in {{\rm \bf LR}}^{\eta}_{\lambda(\omega')\, \tau}\times {{\rm \bf LR}}^{\lambda(\omega)}_{\eta\, \lambda(\omega'')}$ for some $\tau\in \cP_{\mf x}$ and $\eta\in\cP$,}\\
\text{(2) $\nabla_{\mf x}(U)\leq  n(\omega')-(2\lambda(\omega')_1/\epsilon)$\,  and   $\Delta_{\mf x}\left(( V^\sharp \rightarrow  U^\sharp)_R\right)\leq  n(\omega'')$.}
\end{array} \right\}.
\end{split}
\end{equation*}}
\end{thm}
\pf This can be obtained by modifying the arguments in Theorem \ref{LR rule}. \qed\vskip 3mm

Let ${\mf l}_k^\dagger$ be the subalgebra of ${\mf x}_k^\dagger$ isomorphic to $\gl_k$ with positive simple roots $\{\,\alpha_i^\dagger \,|\,1\leq i\leq k-1\,\}$. Let $\left(\tfrac{1}{\epsilon}\Z\right)^k_+=\left\{\,(\mu_1,\ldots,\mu_k)\,\big\vert\,\mu_i\in\tfrac{1}{\epsilon}\Z,\, \mu_i-\mu_{i+1}\in \Z_{\geq 0}\,\right\}$. For $\mu\in \left(\tfrac{1}{\epsilon}\Z\right)^k_+$, let $\B({\mf l}_k^\dagger,\omega_\mu)$ be the ${\mf l}_k^\dagger$-crystal associated to the highest weight $U_q({\mf l}_k^\dagger)$-module with highest weight $\omega_\mu=\mu_1\epsilon^\dagger_1+\ldots +\mu_k\epsilon^\dagger_k$.

\begin{thm}\label{branching in opposite Dynkin}
For a dominant integral ${\mf x}_k^\dagger$-weight $\omega$,  we have as an ${\mf l}_k^\dagger$-crystal
\begin{equation*}
\begin{split}
\B({\mf x}^\dagger_k,\omega)  \simeq  \ \bigsqcup_{\mu}\B({\mf l}^\dagger_k,\omega_\mu)\times {{\rm \bf LR}}^{\omega_\mu}_{\omega}({\mf x}_k^\dagger),
\end{split}
\end{equation*}
where  the sum runs over all $\mu\in \left(\tfrac{1}{\epsilon}\Z\right)^k_+$ such that $\sigma:=\left(\mu_1+\tfrac{\epsilon n(\omega)}{ 2},\ldots,\mu_k+\tfrac{\epsilon n(\omega)}{ 2}\right)$ is a partition and
{\rm $${{\rm \bf LR}}^{\omega_\mu}_{\omega}({\mf x}_k^\dagger)=\{\,U\,|\,\text{$U\in {\rm \bf LR}^{\sigma}_{\tau\, \lambda(\omega)}$ for some $\tau\in \cP_{\mf x}$ and $\Delta_{\mf x}(U)\leq n(\omega)$}\,\}.$$}
\end{thm}
\pf Note that an ${\mf l}^\dagger_k$-highest (resp. lowest) weight element in $\B({\mf x}^\dagger_k,\omega)$ corresponds to an ${\mf l}_k$-lowest  (resp. highest) weight element in $\B({\mf x}_k,\Lambda^{\mf x}(\lambda(\omega), n(\omega)))$, where ${\mf l}_k$ is the subalgebra of ${\mf x}_k$ isomorphic to $\gl_k$ with positive simple roots $\{\,\alpha_i \,|\,1\leq i\leq k-1\,\}$.

Let $b$  be an ${\mf l}_k$-lowest weight element in $\B({\mf x}_k,\Lambda^{\mf x}(\lambda(\omega), n(\omega)))$ of weight $n(\omega)\Lambda^{\mf x}_{0}+\sigma_k\widehat{\epsilon}_1+\cdots+\sigma_1\widehat{\epsilon}_k$ for some partition $\sigma$ with  $\ell(\sigma)\leq k$.
Then as an ${\mf x}^\dagger_k$-weight,
\begin{equation*}
\begin{split}
n(\omega)\Lambda^{\mf x}_{0}+\sigma_k\widehat{\epsilon}_1+\cdots+\sigma_1\widehat{\epsilon}_k
&=-n(\omega)\omega^{\mf x}_k+\sigma_1\epsilon^\dagger_1+\cdots+\sigma_k\epsilon^\dagger_k \\
&=-\tfrac{\epsilon n(\omega)}{2}(\epsilon^\dagger_1+\cdots+\epsilon^\dagger_k)+\sigma_1\epsilon^\dagger_1+\cdots+\sigma_k\epsilon^\dagger_k\\
&=\left(\sigma_1- \tfrac{\epsilon n(\omega)}{2}\right)\epsilon^\dagger_1+\cdots+\left(\sigma_k- \tfrac{\epsilon n(\omega)}{2}\right)\epsilon^\dagger_k.
\end{split}
\end{equation*}
Hence $b$ corresponds to an ${\mf l}^\dagger_k$-highest weight element in $\B({\mf l}^\dagger_k,\omega_\mu)$ of weight $\omega_\mu$, where
\begin{equation*}
\mu=\left(\sigma_1- \frac{\epsilon n(\omega)}{2}, \ldots, \sigma_k- \frac{\epsilon n(\omega)}{2}\right).
\end{equation*}
By Theorem \ref{branching}, the number of connected components in $\B({\mf x}^\dagger_k,\omega) $ isomorphic to $\B({\mf l}^\dagger_k,\omega_\mu)$ is equal to the number of $U\in {\rm \bf LR}^{\sigma}_{\tau\, \lambda(\omega)}$ for some $\tau\in \cP_{\mf x}$ with $\Delta_{\mf x}(U)\leq n(\omega)$.
\qed

\begin{rem}{\rm
The formulas in Theorems \ref{LR in opposite Dynkin} and \ref{branching in opposite Dynkin} have stable limits, but they are different from the known results (cf. \cite[Remark 2.5 (a) and (c)]{HTW})
}
\end{rem}

\section{Character formula for orthosymplectic Lie superalgebras}\label{Character formula for orthosymplectic Lie superalgebras}

\subsection{Superization}
 Let $\A$ be a linearly
ordered $\mathbb{Z}_2$-graded set.
Let $R=\mathbb{Q}[\mathbb{C}]$ denote the group algebra of the additive group $\mathbb{C}$ with a $\mathbb{Q}$-basis $\{\,q^z\,|\,z\in\mathbb{C}\,\}$. We put $Sym_{\A;R}= R\otimes_\mathbb{Q}Sym_{\A}$. The map $\omega_\A$ in (\ref{Omega_A}) can be extended to an $R$-algebra homomorphism $Sym_{\Z_{>0};R} \rightarrow Sym_{\A;R}$. By abuse of notation, we still denote it by $\omega_\A$.

Let us introduce our main object in this section.
\begin{df}{\rm
For $(\lambda,n)\in\cP({\mf x})$, define
\begin{equation*}
\begin{split}
{\mathcal{T}}_\A^{\mf x}(\lambda,n)&=  \left\{\,(S,T)\,\Bigg\vert
\begin{array}{l}
\text{(1) $S\in SST_{\A}(\tau)$ for some $\tau\in\cP_{\mf x}$ and $T\in SST_{\A}(\lambda)$} \\
\text{(2) $\Delta_{\mf x}\left((T\rightarrow S)_R\right)\leq  n$}
\end{array} \right\}.
\end{split}
\end{equation*}
Let $\cP({\mf x})_\A=\{\,(\lambda,n)\in \cP({\mf x})\,|\,{\mathcal{T}}_\A^{\mf x}(\lambda,n)\neq \emptyset\,\}$ and for  $(\lambda,n)\in \cP({\mf x})_\A$ define
\begin{equation*}
S^{\mf x}_{(\lambda,n)}({\bf x}_\A)=q^n\sum_{(S,T)\in {\mathcal{T}}_\A^{\mf x}(\lambda,n)}x_\A^{S}x_A^T.
\end{equation*}
}
\end{df}\vskip 2mm

\begin{prop}\label{Branching and LR}\mbox{}
\begin{itemize}
\item[(1)] For $(\lambda,n)\in\cP({\mf x})_\A$, we have
\begin{equation*}
\begin{split}
S^{\mf x}_{(\lambda,n)}({\bf x}_\A)&=\omega_{\A}\left( S^{\mf x}_{(\lambda,n)}({\bf x}_{\Z_{>0}}) \right),\\
S^{\mf x}_{(\lambda,n)}({\bf x}_{\A})&=q^n\sum_{\sigma\in\cP}c^{\,\sigma}_{(\lambda,n)}({\mf x})  S_{\sigma}({\bf x}_{\A}),
\end{split}
\end{equation*}
where $c^{\,\sigma}_{(\lambda,n)}({\mf x}) $ is given in Definition \ref{branching coefficient}.

\item[(2)] For $(\mu,m), (\nu,n)\in\cP({\mf x})_\A$,
$$S^{\mf x}_{(\mu ,m)}({\bf x}_{\A})S^{\mf x}_{(\nu,n)}({\bf x}_{\A})=\sum_{(\lambda,m+n)\in\cP({\mf x})_\A}c^{(\lambda,m+n)}_{(\mu,m) (\nu,n)}({\mf x}) S^{\mf x}_{(\lambda,m+n)}({\bf x}_{\A}),$$
where  $c^{(\lambda,m+n)}_{(\mu,m) (\nu,n)}({\mf x})$ is given in Definition \ref{LR coefficient}.
\end{itemize}
\end{prop}
\pf  By Theorem \ref{main result}, $S^{\mf x}_{(\lambda,n)}({\bf x}_{\Z_{>0}})$ is the character of $\B({\mf x}_\infty,\La^{\mf x}(\lambda,n))$ where $q=e^{\Lambda^{\mf x}_0}$. Since $\omega_\A$ is an homomorphism of algebras,  the identities in (1) follow from (\ref{LR}) and Corollary \ref{branching-character}, repsectively.
The identity in (2) follows from Corollary \ref{LR rule-character}.
\qed\vskip 2mm

\begin{prop}\label{Linear independence of Sxlambda}
$\{\,S^{\mf x}_{(\lambda,n)}({\bf x}_\A)\,|\,(\lambda,n)\in \cP({\mf x})_\A\,\}$ is linearly independent over $\mathbb{Q}$.
\end{prop}
\pf  It suffices to show that  $\{\,S^{\mf x}_{(\lambda,n)}({\bf x}_\A)\,|\,(\lambda,n)\in \cP({\mf x})_\A\,\}$ is linearly independent, where $n$ is a fixed integer.

For $\mu,\nu\in\cP$, define $\mu > \nu$ if and only if there exists $k$ such that $\mu_i=\nu_i$ for $1\leq i\leq k$ and $\mu_k>\nu_k$. Then $>$ is a linear ordering on $\cP$. By Proposition \ref{Branching and LR} (1), we have
\begin{equation*}
S^{\mf x}_{(\lambda,n)}({\bf x}_{\A})=q^n\sum_{\sigma\geq \lambda} a_{\lambda \sigma } S_{\sigma}({\bf x}_{\A})
\end{equation*}
for some $a_{\lambda \sigma}\in\Z_{\geq 0}$ with $a_{\lambda \lambda}=1$. This together with the linear independence of $S_{\sigma}({\bf x}_\A)$ implies the linear independence of $\{\,S^{\mf x}_{(\lambda,n)}({\bf x}_\A)\,|\,(\lambda,n)\in \cP({\mf x})_\A\,\}$. \qed\vskip 2mm
\vskip 2mm

The main goal of this section is to show that $S^{\mf x}_{(\lambda,n)}({\bf x}_\A)$ is the character of an irreducible module over an orthosymplectic Lie superalgebra  using the  Cheng-Lam-Wang's super duality \cite{CLW}.

\subsection{Orthosymplectic Lie superalgebras}\label{Osp}
Let $V=V_0\oplus V_1$ be a complex superspace. For $v\in V_p$, we put $|v|=p$ ($p\in\Z_2$). Let $\gl(V)$ be the general linear Lie superalgebra of linear endomorphisms on $V$ which vanishes on a subspace of finite codimension.
We define $\frak{osp}(V)$ to be the subalgebra of $\gl(V)$ which preserves a non-degenerate supersymetric bilinear form $(\cdot|\cdot)$ on $V$, that is, $\frak{osp}(V)=\frak{osp}(V)_{{0}}\oplus\frak{osp}(V)_1$
with
\begin{equation}
\frak{osp}(V)_{p}=\{A\in\gl(V)_p\vert (Av|w)+(-1)^{p |v|}(v|Aw)=0 \text{ for homogeneous } v,w \in V \}.
\end{equation}
Similarly, we define $\frak{spo}(V)$ to be the subalgebra of $\gl(V)$ which preserves a non-degenerate skew-supersymetric bilinear form $(\cdot|\cdot)$ on $V$ (cf. \cite{K1}).\vskip 2mm

Now, let us consider the following $\Z_2$-graded linearly ordered sets ($m\in\Z_{\geq 0}$): \vskip 2mm
\begin{itemize}
\item[$\cdot$] $\td{\mathcal{I}}_m=\{\,\overline{k},-\overline{k}\,|\,1\leq k\leq m\,\}\cup \tfrac{1}{2}\Z$,\vskip 2mm

\item[$\cdot$] $\ov{\mathcal{I}}_m=\{\,\overline{k},-\overline{k} \,|\,1\leq k\leq m\,\}\cup \left(\tfrac{1}{2}+\Z \right)\cup\{0\}$,\vskip 2mm

\item[$\cdot$] $\mathcal{I}_m=\{\,\overline{k},-\overline{k} \,|\,1\leq k\leq m\,\}\cup \Z$,\vskip 2mm

\item[$\cdot$] $\td{\mathcal{I}}_m^{\times}=\td{\mathcal{I}}_m\setminus\{0\},\ \ \ov{\mathcal{I}}_m^{\times}=\ov{\mathcal{I}}_m\setminus\{0\},\ \ {\mathcal{I}}_m^{\times}={\mathcal{I}}_m\setminus\{0\}$,
\end{itemize}
where
\begin{equation}
\begin{split}
(\td{\mathcal{I}}_m)_0=\td{\mathcal{I}}_m\setminus \left(\tfrac{1}{2}+\Z\right), \ &\ \ (\td{\mathcal{I}}_m)_1=\tfrac{1}{2}+\Z,\\
\cdots<-\tfrac{3}{2} <-1 <-\tfrac{1}{2}  <-\ov{1} <\cdots<-\ov{m} <&0<\ov{m}<\cdots<\ov{1}<\tfrac{1}{2}<1 < \tfrac{3}{2}<\ldots \ \ ,
\end{split}
\end{equation}
and the $\Z_2$-gradings and linear orderings on the other sets are induced from those of $\td{\mathcal{I}}_m$. Here we assume that $\td{\mathcal{I}}_0=\tfrac{1}{2}\Z$, $\ov{\mathcal{I}}_0=\left(\frac{1}{2}+\Z \right)\cup\{0\}$ and $\mathcal{I}_0=\Z$.
In the rest of this section, we assume the following:
\begin{itemize}
\item[$\cdot$] $\mathcal{I}$ denotes one of $\td{\mathcal{I}}_m$, $\ov{\mathcal{I}}_m$, ${\mathcal{I}}_m$, $\td{\mathcal{I}}^\times_m$, $\ov{\mathcal{I}}^\times_m$ and ${\mathcal{I}}^\times_m$,

\item[$\cdot$] $\mathcal{I}^+=\{\,a\,|\,a\in \mathcal{I}, a>0\,\}$,

\item[$\cdot$] $V_{\mathcal{I}}$ is a superspace with basis $\{\,v_a\,|\,a\in \mathcal{I}\,\}$.
\end{itemize}
We identify $\gl(V_{\mathcal{I}})$ with the Lie superalgebra  of  matrices $(a_{i,j})_{i,j\in \mathcal{I}}$  with
$a_{i,j}=0$ for all but finitely many $a_{i,j}$'s. Denote by $E_{i,j}$
the elementary matrix with $1$ at the $i$-th row and $j$-th column
and zero elsewhere. Also, let
$\widehat{\gl}(V_{\td{\I}_m})=\gl(V_{\td{\I}_m})\oplus\C K$ the central extension
of $\gl(V_{\td{\I}_m})$ by a one-dimensional center $\C K$ given by the
$2$-cocycle $\beta(A,B)={\rm Str}([{J},A]B)$,
where ${\rm Str}$ is the supertrace   and ${J}=E_{0,0}+\sum_{r\le -\hf}E_{r,r}$.

First, suppose that $\mathcal{I}$ is one of $\td{\mathcal{I}}_m$, $\ov{\I}_m$ and ${\I}_m$. Define a supersymmetric bilinear form $(\cdot|\cdot)$ on $V_{\mathcal{I}}$ by
\begin{equation}
\begin{split}
&(v_{\pm a}|v_{\pm b})=0, \ \ (v_a|v_{-b})=(-1)^{|a||b|}(v_{-b}|v_a)=\delta_{a b}, \\
&(v_0|v_0)=1, \ \ (v_0|v_{\pm a})=0,
\end{split}
\end{equation}
for $a,b \in \mathcal{I}^+$. Let $\frak{osp}(V_\mathcal{I})$ be the subalgebra of $\gl(V_\mathcal{I})$ preserving $(\cdot|\cdot)$. The Cartan subalgebra   is spanned by $\{\,E_{a}=E_{a,a}-E_{-a,-a}\,|\,a\in \mathcal{I}^+ \,\}$.
We choose a Borel subalgebra spanned by the upper triangular matrices. With respect to  the corresponding dual basis $\{\,\delta_a\,|\,a\in \mathcal{I}^+\,\}$, the set of simple roots  are given by
\begin{equation}
\begin{split}
\td{\I}_m \,:& \, \{\, -\delta_{\ov{m}},\, \delta_{\ov{k}}-\delta_{\ov{k-1}}\  (2\leq k\leq m),\ \delta_{\ov{1}}-\delta_{\tfrac{1}{2}},\ \delta_r-\delta_{r+\tfrac{1}{2}} \ (r\in\tfrac{1}{2}\Z_{> 0}) \,\}  ,\\
\ov{\I}_m \,:& \, \{\, -\delta_{\ov{m}},\, \delta_{\ov{k}}-\delta_{\ov{k-1}}\  (2\leq k\leq m),\ \delta_{\ov{1}}-\delta_{\tfrac{1}{2}},\ \delta_r-\delta_{r+1} \ (r\in\tfrac{1}{2}+\Z_{\geq  0}) \,\}  ,\\
{\I}_m \,:& \, \{\, -\delta_{\ov{m}},\, \delta_{\ov{k}}-\delta_{\ov{k-1}}\  (2\leq k\leq m),\ \delta_{\ov{1}}-\delta_{1},\ \delta_r-\delta_{r+1} \ (r\in\Z_{> 0}) \,\},
\end{split}
\end{equation}
for $m\geq 1$, and
\begin{equation}
\begin{split}
\td{\I}_0 \,:& \, \{\, -\delta_{\hf}, \ \delta_r-\delta_{r+\hf} \ (r\in\tfrac{1}{2}\Z_{>0})\, \},\\
\ov{\I}_0 \,:& \, \{\, -\delta_{\hf},\,  \delta_r-\delta_{r+1} \ (r\in\tfrac{1}{2}+\Z_{\geq  0}) \,\},\\
{\I}_0 \,:& \, \{\, -\delta_{1},\,  \delta_r-\delta_{r+1} \ (r\in\Z_{> 0}) \,\}.
\end{split}
\end{equation}
Note that when $\mathcal{I}={{\I}_m}$ ($m\geq 0$),  the set of simple roots is of type $B_\infty$.
The associated  Dynkin diagrams are as follows.
\begin{center}
\hskip -3cm \setlength{\unitlength}{0.16in}
\begin{picture}(24,4)
\put(1.5,2){\makebox(0,0)[c]{$\td{\I}_{m}$ : }}
\put(30,2){\makebox(0,0)[c]{$(m\geq 1)$}}
\put(5.6,2){\makebox(0,0)[c]{$\bigcirc$}}
\put(8,2){\makebox(0,0)[c]{$\bigcirc$}}
\put(10.4,2){\makebox(0,0)[c]{$\bigcirc$}}
\put(14.85,2){\makebox(0,0)[c]{$\bigcirc$}}
\put(17.25,2){\makebox(0,0)[c]{$\bigotimes$}}
\put(19.4,2){\makebox(0,0)[c]{$\bigotimes$}}
\put(23.9,2){\makebox(0,0)[c]{$\bigotimes$}}
\put(8.35,2){\line(1,0){1.5}}
\put(10.82,2){\line(1,0){0.8}}
\put(13.2,2){\line(1,0){1.2}}
\put(15.28,2){\line(1,0){1.45}}
\put(17.7,2){\line(1,0){1.25}}
\put(19.81,2){\line(1,0){0.9}}
\put(22.1,2){\line(1,0){1.28}}
\put(24.3,2){\line(1,0){1.28}}
\put(6.8,2){\makebox(0,0)[c]{$\Longleftarrow$}}
\put(12.5,1.95){\makebox(0,0)[c]{$\cdots$}}
\put(21.5,1.95){\makebox(0,0)[c]{$\cdots$}}
\put(26.5,1.95){\makebox(0,0)[c]{$\cdots$}}
\put(5.4,0.8){\makebox(0,0)[c]{\tiny $\delta^\ast=-\delta_{\ov{m}}$}}
\put(7.8,3){\makebox(0,0)[c]{\tiny $\delta_{\ov{m}}-\delta_{\ov{m-1}}$}}
\put(10.4,0.8){\makebox(0,0)[c]{\tiny $\delta_{\ov{m-1}}-\delta_{\ov{m-2}}$}}
\put(14.5,3){\makebox(0,0)[c]{\tiny $\delta_{\ov{2}}-\delta_{\ov{1}}$}}
\put(17.15,0.8){\makebox(0,0)[c]{\tiny $\delta^{\times}=\delta_{\ov{1}}-\delta_{\tfrac{1}{2}}$}}
\put(19.5,3){\makebox(0,0)[c]{\tiny $\delta_{\tfrac{1}{2}}-\delta_{1}$}}
\put(23.8,0.8){\makebox(0,0)[c]{\tiny $\delta_r-\delta_{r+\tfrac{1}{2}}$}}
\end{picture}\vskip 5mm

\hskip -3cm \setlength{\unitlength}{0.16in}
\begin{picture}(24,4)
\put(1.5,2){\makebox(0,0)[c]{$\td{\I}_{0}$ : }}
\put(5.6,2){\makebox(0,0)[c]{$\bigcirc$}}
\put(8,2){\makebox(0,0)[c]{$\bigotimes$}}
\put(10.4,2){\makebox(0,0)[c]{$\bigotimes$}}
\put(14.85,2){\makebox(0,0)[c]{$\bigotimes$}}
\put(17.45,2){\makebox(0,0)[c]{$\cdots$}}
\put(8.35,2){\line(1,0){1.5}}
\put(10.82,2){\line(1,0){0.8}}
\put(13.2,2){\line(1,0){1.2}}
\put(15.28,2){\line(1,0){1.45}}
\put(6.8,2){\makebox(0,0)[c]{$\Longleftarrow$}}
\put(12.5,1.95){\makebox(0,0)[c]{$\cdots$}}
\put(5.4,0.8){\makebox(0,0)[c]{\tiny $\delta^\ast=-\delta_{\hf}$}}
\put(7.8,3){\makebox(0,0)[c]{\tiny $\delta_{\hf}-\delta_{1}$}}
\put(10.4,0.8){\makebox(0,0)[c]{\tiny $\delta_{1}-\delta_{\tfrac{3}{2}}$}}
\put(15,3){\makebox(0,0)[c]{\tiny $\delta_{r}-\delta_{r+\hf}$}}
\end{picture}\vskip 5mm

\hskip -3cm \setlength{\unitlength}{0.16in}
\begin{picture}(24,4)
\put(1.5,2){\makebox(0,0)[c]{$\ov{\I}_{m}$ : }}
\put(30,2){\makebox(0,0)[c]{$(m\geq 1)$}}
\put(5.6,2){\makebox(0,0)[c]{$\bigcirc$}}
\put(8,2){\makebox(0,0)[c]{$\bigcirc$}}
\put(10.4,2){\makebox(0,0)[c]{$\bigcirc$}}
\put(14.85,2){\makebox(0,0)[c]{$\bigcirc$}}
\put(17.25,2){\makebox(0,0)[c]{$\bigotimes$}}
\put(19.4,2){\makebox(0,0)[c]{$\bigcirc$}}
\put(23.9,2){\makebox(0,0)[c]{$\bigcirc$}}
\put(8.35,2){\line(1,0){1.5}}
\put(10.82,2){\line(1,0){0.8}}
\put(13.2,2){\line(1,0){1.2}}
\put(15.28,2){\line(1,0){1.45}}
\put(17.7,2){\line(1,0){1.25}}
\put(19.81,2){\line(1,0){0.9}}
\put(22.1,2){\line(1,0){1.28}}
\put(24.3,2){\line(1,0){1.28}}
\put(6.8,2){\makebox(0,0)[c]{$\Longleftarrow$}}
\put(12.5,1.95){\makebox(0,0)[c]{$\cdots$}}
\put(21.5,1.95){\makebox(0,0)[c]{$\cdots$}}
\put(26.5,1.95){\makebox(0,0)[c]{$\cdots$}}
\put(5.4,0.8){\makebox(0,0)[c]{\tiny $\delta^\ast=-\delta_{\ov{m}}$}}
\put(7.8,3){\makebox(0,0)[c]{\tiny $\delta_{\ov{m}}-\delta_{\ov{m-1}}$}}
\put(10.4,0.8){\makebox(0,0)[c]{\tiny $\delta_{\ov{m-1}}-\delta_{\ov{m-2}}$}}
\put(14.5,3){\makebox(0,0)[c]{\tiny $\delta_{\ov{2}}-\delta_{\ov{1}}$}}
\put(17.15,0.8){\makebox(0,0)[c]{\tiny $\delta^\times=\delta_{\ov{1}}-\delta_{\tfrac{1}{2}}$}}
\put(19.5,3){\makebox(0,0)[c]{\tiny $\delta_{\tfrac{1}{2}}-\delta_{\tfrac{3}{2}}$}}
\put(23.8,0.8){\makebox(0,0)[c]{\tiny $\delta_r-\delta_{r+1}$}}
\end{picture}\vskip 5mm

\hskip -3cm \setlength{\unitlength}{0.16in}
\begin{picture}(24,4)
\put(1.8,2){\makebox(0,0)[c]{$\ov{\I}_{0}$ : }}
\put(6,2){\makebox(0,0)[c]{\circle*{0.9}}}
\put(8,2){\makebox(0,0)[c]{$\bigcirc$}}
\put(10.4,2){\makebox(0,0)[c]{$\bigcirc$}}
\put(14.85,2){\makebox(0,0)[c]{$\bigcirc$}}
\put(17.45,2){\makebox(0,0)[c]{$\cdots$}}
\put(8.35,2){\line(1,0){1.5}}
\put(10.82,2){\line(1,0){0.8}}
\put(13.2,2){\line(1,0){1.2}}
\put(15.28,2){\line(1,0){1.45}}
\put(6.8,2){\makebox(0,0)[c]{$\Longleftarrow$}}
\put(12.5,1.95){\makebox(0,0)[c]{$\cdots$}}
\put(5.4,0.8){\makebox(0,0)[c]{\tiny $\delta^\ast=-\delta_{\hf}$}}
\put(7.8,3){\makebox(0,0)[c]{\tiny $\delta_{\hf}-\delta_{\tfrac{3}{2}}$}}
\put(10.4,0.8){\makebox(0,0)[c]{\tiny $\delta_{\tfrac{3}{2}}-\delta_{\tfrac{5}{2}}$}}
\put(15,3){\makebox(0,0)[c]{\tiny $\delta_{r}-\delta_{r+1}$}}
\end{picture}\vskip 5mm

\hskip -3cm \setlength{\unitlength}{0.16in}
\begin{picture}(24,4)
\put(1.5,2){\makebox(0,0)[c]{${\I}_{m}$ : }}
\put(30,2){\makebox(0,0)[c]{$(m\geq 1)$}}
\put(5.6,2){\makebox(0,0)[c]{$\bigcirc$}}
\put(8,2){\makebox(0,0)[c]{$\bigcirc$}}
\put(10.4,2){\makebox(0,0)[c]{$\bigcirc$}}
\put(14.85,2){\makebox(0,0)[c]{$\bigcirc$}}
\put(17.25,2){\makebox(0,0)[c]{$\bigcirc$}}
\put(19.4,2){\makebox(0,0)[c]{$\bigcirc$}}
\put(23.9,2){\makebox(0,0)[c]{$\bigcirc$}}
\put(8.35,2){\line(1,0){1.5}}
\put(10.82,2){\line(1,0){0.8}}
\put(13.2,2){\line(1,0){1.2}}
\put(15.28,2){\line(1,0){1.45}}
\put(17.7,2){\line(1,0){1.25}}
\put(19.81,2){\line(1,0){0.9}}
\put(22.1,2){\line(1,0){1.28}}
\put(24.3,2){\line(1,0){1.28}}
\put(6.8,2){\makebox(0,0)[c]{$\Longleftarrow$}}
\put(12.5,1.95){\makebox(0,0)[c]{$\cdots$}}
\put(21.5,1.95){\makebox(0,0)[c]{$\cdots$}}
\put(26.5,1.95){\makebox(0,0)[c]{$\cdots$}}
\put(5.4,0.8){\makebox(0,0)[c]{\tiny $\delta^\ast=-\delta_{\ov{m}}$}}
\put(7.8,3){\makebox(0,0)[c]{\tiny $\delta_{\ov{m}}-\delta_{\ov{m-1}}$}}
\put(10.4,0.8){\makebox(0,0)[c]{\tiny $\delta_{\ov{m-1}}-\delta_{\ov{m-2}}$}}
\put(14.5,3){\makebox(0,0)[c]{\tiny $\delta_{\ov{2}}-\delta_{\ov{1}}$}}
\put(17.15,0.8){\makebox(0,0)[c]{\tiny $\delta^\times=\delta_{\ov{1}}-\delta_{1}$}}
\put(19.5,3){\makebox(0,0)[c]{\tiny $\delta_{1}-\delta_{2}$}}
\put(23.8,0.8){\makebox(0,0)[c]{\tiny $\delta_r-\delta_{r+1}$}}
\end{picture}\vskip 5mm

\hskip -3cm \setlength{\unitlength}{0.16in}
\begin{picture}(24,4)
\put(1.8,2){\makebox(0,0)[c]{${\I}_{0}$ : }}
\put(5.6,2){\makebox(0,0)[c]{$\bigcirc$}}
\put(8,2){\makebox(0,0)[c]{$\bigcirc$}}
\put(10.4,2){\makebox(0,0)[c]{$\bigcirc$}}
\put(14.85,2){\makebox(0,0)[c]{$\bigcirc$}}
\put(17.45,2){\makebox(0,0)[c]{$\cdots$}}
\put(8.35,2){\line(1,0){1.5}}
\put(10.82,2){\line(1,0){0.8}}
\put(13.2,2){\line(1,0){1.2}}
\put(15.28,2){\line(1,0){1.45}}
\put(6.8,2){\makebox(0,0)[c]{$\Longleftarrow$}}
\put(12.5,1.95){\makebox(0,0)[c]{$\cdots$}}
\put(5.4,0.8){\makebox(0,0)[c]{\tiny $\delta^\ast=-\delta_{1}$}}
\put(7.8,3){\makebox(0,0)[c]{\tiny $\delta_{1}-\delta_{2}$}}
\put(10.4,0.8){\makebox(0,0)[c]{\tiny $\delta_{2}-\delta_{3}$}}
\put(15,3){\makebox(0,0)[c]{\tiny $\delta_{r}-\delta_{r+1}$}}
\end{picture}\vskip 5mm
\end{center}
\noindent Here $\delta^\ast$ and $\delta^\times$ denote  distinguished simple roots in each diagram, which are necessary for our later arguments.

Next, suppose that $\mathcal{I}$ is one of $\td{\I}^\times_m$, $\ov{\I}^\times_m$ and ${\I}^\times_m$.
Define a skew-supersymmetric bilinear form $(\cdot|\cdot)$ on $V_\mathcal{I}$ by
\begin{equation}
\begin{split}
&(v_{\pm a}|v_{\pm b})=0, \ \ (v_a|v_{-b})=-(-1)^{|a||b|}(v_{-b}|v_a)=\delta_{a b}, \\
\end{split}
\end{equation}
for $a,b \in \mathcal{I}^+$. Let $\frak{spo}(V_\mathcal{I})$ be the subalgebra of $\gl(V_\mathcal{I})$ preserving $(\cdot|\cdot)$. The Cartan subalgebra  is spanned by $\{\,E_{a}=E_{a,a}-E_{-a, -a}\,|\,a\in \mathcal{I}^+\,\}$. We choose a Borel subalgebra spanned by the upper triangular matrices. With respect to the corresponding dual basis $\{\,\delta_a\,|\,a\in \mathcal{I}^+\,\}$, the set of simple roots are given by
\begin{equation}
\begin{split}
\td{\I}^\times_m \,:& \, \{\, -2\delta_{\ov{m}},\, \delta_{\ov{k}}-\delta_{\ov{k-1}}\  (2\leq k\leq m),\ \delta_{\ov{1}}-\delta_{\tfrac{1}{2}},\ \delta_r-\delta_{r+\tfrac{1}{2}} \ (r\in\tfrac{1}{2}\Z_{> 0}) \,\}  ,\\
\ov{\I}^\times_m \,:& \, \{\, -2\delta_{\ov{m}},\, \delta_{\ov{k}}-\delta_{\ov{k-1}}\  (2\leq k\leq m),\ \delta_{\ov{1}}-\delta_{\tfrac{1}{2}},\ \delta_r-\delta_{r+1} \ (r\in\tfrac{1}{2}+\Z_{\geq  0}) \,\}  ,\\
{\I}^\times_m \,:& \, \{\, -2\delta_{\ov{m}},\, \delta_{\ov{k}}-\delta_{\ov{k-1}}\  (2\leq k\leq m),\ \delta_{\ov{1}}-\delta_{1},\ \delta_r-\delta_{r+1} \ (r\in\Z_{> 0}) \,\},
\end{split}
\end{equation}
for $m\geq 1$, and
\begin{equation}
\begin{split}
\td{\I}^\times_0 \,:& \, \{\,-\delta_{\hf}-\delta_1, \
\delta_r-\delta_{r+\hf} \ (r\in\tfrac{1}{2}\Z_{>0})\, \}  ,\\
\ov{\I}^\times_0 \,:& \, \{\, -\delta_{\tfrac{1}{2}}-\delta_{\tfrac{3}{2}},  \ \delta_r-\delta_{r+1} \ (r\in\tfrac{1}{2}+\Z_{\geq  0}) \,\},\\
{\I}^\times_0 \,:& \, \{\, -2\delta_{1},\ \delta_r-\delta_{r+1} \ (r\in\Z_{> 0}) \,\}.
\end{split}
\end{equation}
Note that when $\mathcal{I}={\I}^\times_m$ ($m\geq 0$),  the set of simple roots is of type $C_\infty$. The associated Dynkin diagrams  are as follows.
\begin{center}
\hskip -3cm \setlength{\unitlength}{0.16in}
\begin{picture}(24,4)
\put(1.5,2){\makebox(0,0)[c]{$\td{\I}^\times_{m}$ : }}
\put(30,2){\makebox(0,0)[c]{$(m\geq 1)$}}
\put(5.6,2){\makebox(0,0)[c]{$\bigcirc$}}
\put(8,2){\makebox(0,0)[c]{$\bigcirc$}}
\put(10.4,2){\makebox(0,0)[c]{$\bigcirc$}}
\put(14.85,2){\makebox(0,0)[c]{$\bigcirc$}}
\put(17.25,2){\makebox(0,0)[c]{$\bigotimes$}}
\put(19.4,2){\makebox(0,0)[c]{$\bigotimes$}}
\put(23.9,2){\makebox(0,0)[c]{$\bigotimes$}}
\put(8.35,2){\line(1,0){1.5}}
\put(10.82,2){\line(1,0){0.8}}
\put(13.2,2){\line(1,0){1.2}}
\put(15.28,2){\line(1,0){1.45}}
\put(17.7,2){\line(1,0){1.25}}
\put(19.81,2){\line(1,0){0.9}}
\put(22.1,2){\line(1,0){1.28}}
\put(24.3,2){\line(1,0){1.28}}
\put(6.8,2){\makebox(0,0)[c]{$\Longrightarrow$}}
\put(12.5,1.95){\makebox(0,0)[c]{$\cdots$}}
\put(21.5,1.95){\makebox(0,0)[c]{$\cdots$}}
\put(26.5,1.95){\makebox(0,0)[c]{$\cdots$}}
\put(5.4,0.8){\makebox(0,0)[c]{\tiny $\delta^\ast=-2\delta_{\ov{m}}$}}
\put(7.8,3){\makebox(0,0)[c]{\tiny $\delta_{\ov{m}}-\delta_{\ov{m-1}}$}}
\put(10.4,0.8){\makebox(0,0)[c]{\tiny $\delta_{\ov{m-1}}-\delta_{\ov{m-2}}$}}
\put(14.5,3){\makebox(0,0)[c]{\tiny $\delta_{\ov{2}}-\delta_{\ov{1}}$}}
\put(17.15,0.8){\makebox(0,0)[c]{\tiny $\delta^\times=\delta_{\ov{1}}-\delta_{\tfrac{1}{2}}$}}
\put(19.5,3){\makebox(0,0)[c]{\tiny $\delta_{\tfrac{1}{2}}-\delta_{1}$}}
\put(23.8,0.8){\makebox(0,0)[c]{\tiny $\delta_r-\delta_{r+\tfrac{1}{2}}$}}
\end{picture}\vskip 5mm

\hskip -3cm \setlength{\unitlength}{0.16in} \medskip
\begin{picture}(24,5.5)
\put(1.5,2){\makebox(0,0)[c]{$\td{\I}^\times_{0} $ : }}
\put(6,0){\makebox(0,0)[c]{$\bigotimes$}}
\put(6,4){\makebox(0,0)[c]{$\bigotimes$}}
\put(8,2){\makebox(0,0)[c]{$\bigotimes$}}
\put(10.4,2){\makebox(0,0)[c]{$\bigotimes$}}
\put(14.85,2){\makebox(0,0)[c]{$\bigotimes$}}
\put(17.65,2){\makebox(0,0)[c]{$\cdots$}}

\put(6.35,0.3){\line(1,1){1.35}}
\put(6,0.4){\line(0,1){3.25}} \put(6.35,3.7){\line(1,-1){1.35}}
\put(8.4,2){\line(1,0){1.55}} \put(10.82,2){\line(1,0){0.8}}
\put(13.2,2){\line(1,0){1.2}} \put(15.28,2){\line(1,0){1.45}}

%
\put(12.5,1.95){\makebox(0,0)[c]{$\cdots$}}
%
\put(6,5){\makebox(0,0)[c]{\tiny $\delta^\ast=-\delta_{\hf}-\delta_1$}}
\put(6,-1.2){\makebox(0,0)[c]{\tiny $\delta_{\hf}-\delta_{1},$}}
\put(8.2,3.5){\makebox(0,0)[c]{\tiny $\delta_{1}-\delta_{\tfrac{3}{2}},$}}
\put(10.4,1){\makebox(0,0)[c]{\tiny $\delta_{\tfrac{3}{2}}-\delta_2$}}
\put(14.7,.8){\makebox(0,0)[c]{\tiny $\delta_r-\delta_{r+\hf}$}}
\end{picture}\vskip 6mm

\hskip -3cm \setlength{\unitlength}{0.16in}
\begin{picture}(24,4)
\put(1.5,2){\makebox(0,0)[c]{$\ov{\I}^\times_{m}$ : }}
\put(30,2){\makebox(0,0)[c]{$(m\geq 1)$}}
\put(5.6,2){\makebox(0,0)[c]{$\bigcirc$}}
\put(8,2){\makebox(0,0)[c]{$\bigcirc$}}
\put(10.4,2){\makebox(0,0)[c]{$\bigcirc$}}
\put(14.85,2){\makebox(0,0)[c]{$\bigcirc$}}
\put(17.25,2){\makebox(0,0)[c]{$\bigotimes$}}
\put(19.4,2){\makebox(0,0)[c]{$\bigcirc$}}
\put(23.9,2){\makebox(0,0)[c]{$\bigcirc$}}
\put(8.35,2){\line(1,0){1.5}}
\put(10.82,2){\line(1,0){0.8}}
\put(13.2,2){\line(1,0){1.2}}
\put(15.28,2){\line(1,0){1.45}}
\put(17.7,2){\line(1,0){1.25}}
\put(19.81,2){\line(1,0){0.9}}
\put(22.1,2){\line(1,0){1.28}}
\put(24.3,2){\line(1,0){1.28}}
\put(6.8,2){\makebox(0,0)[c]{$\Longrightarrow$}}
\put(12.5,1.95){\makebox(0,0)[c]{$\cdots$}}
\put(21.5,1.95){\makebox(0,0)[c]{$\cdots$}}
\put(26.5,1.95){\makebox(0,0)[c]{$\cdots$}}
\put(5.4,0.8){\makebox(0,0)[c]{\tiny $\delta^\ast=-2\delta_{\ov{m}}$}}
\put(7.8,3){\makebox(0,0)[c]{\tiny $\delta_{\ov{m}}-\delta_{\ov{m-1}}$}}
\put(10.4,0.8){\makebox(0,0)[c]{\tiny $\delta_{\ov{m-1}}-\delta_{\ov{m-2}}$}}
\put(14.5,3){\makebox(0,0)[c]{\tiny $\delta_{\ov{2}}-\delta_{\ov{1}}$}}
\put(17.15,0.8){\makebox(0,0)[c]{\tiny $\delta^\times=\delta_{\ov{1}}-\delta_{\tfrac{1}{2}}$}}
\put(19.5,3){\makebox(0,0)[c]{\tiny $\delta_{\tfrac{1}{2}}-\delta_{\tfrac{3}{2}}$}}
\put(23.8,0.8){\makebox(0,0)[c]{\tiny $\delta_r-\delta_{r+1}$}}
\end{picture}\vskip 6mm

\hskip -3cm \setlength{\unitlength}{0.16in} \medskip
\begin{picture}(24,5.5)
\put(1.5,2){\makebox(0,0)[c]{$\ov{\I}^\times_{0} $ : }}
\put(6,0){\makebox(0,0)[c]{$\bigcirc$}}
\put(6,4){\makebox(0,0)[c]{$\bigcirc$}}
\put(8,2){\makebox(0,0)[c]{$\bigcirc$}}
\put(10.4,2){\makebox(0,0)[c]{$\bigcirc$}}
\put(14.85,2){\makebox(0,0)[c]{$\bigcirc$}}
\put(17.65,2){\makebox(0,0)[c]{$\cdots$}}

\put(6.35,0.3){\line(1,1){1.35}} \put(6.35,3.7){\line(1,-1){1.35}}
\put(8.4,2){\line(1,0){1.55}} \put(10.82,2){\line(1,0){0.8}}
\put(13.2,2){\line(1,0){1.2}} \put(15.28,2){\line(1,0){1.45}}

%
\put(12.5,1.95){\makebox(0,0)[c]{$\cdots$}}
%
\put(6,5){\makebox(0,0)[c]{\tiny $\delta^\ast=-\delta_{\hf}-\delta_{\tfrac{3}{2}}$}}
\put(6,-1.2){\makebox(0,0)[c]{\tiny $\delta_{\hf}-\delta_{\tfrac{3}{2}},$}}
\put(8.2,3.5){\makebox(0,0)[c]{\tiny $\delta_{\tfrac{3}{2}}-\delta_{\tfrac{5}{2}},$}}
\put(10.4,1){\makebox(0,0)[c]{\tiny $\delta_{\tfrac{5}{2}}-\delta_{\tfrac{7}{2}}$}}
\put(14.7,.8){\makebox(0,0)[c]{\tiny $\delta_r-\delta_{r+1}$}}
\end{picture}\vskip 5mm

\hskip -3cm \setlength{\unitlength}{0.16in}
\begin{picture}(24,4)
\put(1.5,2){\makebox(0,0)[c]{${\I}^\times_{m}$ : }}
\put(30,2){\makebox(0,0)[c]{$(m\geq 1)$}}
\put(5.6,2){\makebox(0,0)[c]{$\bigcirc$}}
\put(8,2){\makebox(0,0)[c]{$\bigcirc$}}
\put(10.4,2){\makebox(0,0)[c]{$\bigcirc$}}
\put(14.85,2){\makebox(0,0)[c]{$\bigcirc$}}
\put(17.25,2){\makebox(0,0)[c]{$\bigcirc$}}
\put(19.4,2){\makebox(0,0)[c]{$\bigcirc$}}
\put(23.9,2){\makebox(0,0)[c]{$\bigcirc$}}
\put(8.35,2){\line(1,0){1.5}}
\put(10.82,2){\line(1,0){0.8}}
\put(13.2,2){\line(1,0){1.2}}
\put(15.28,2){\line(1,0){1.45}}
\put(17.7,2){\line(1,0){1.25}}
\put(19.81,2){\line(1,0){0.9}}
\put(22.1,2){\line(1,0){1.28}}
\put(24.3,2){\line(1,0){1.28}}
\put(6.8,2){\makebox(0,0)[c]{$\Longrightarrow$}}
\put(12.5,1.95){\makebox(0,0)[c]{$\cdots$}}
\put(21.5,1.95){\makebox(0,0)[c]{$\cdots$}}
\put(26.5,1.95){\makebox(0,0)[c]{$\cdots$}}
\put(5.4,0.8){\makebox(0,0)[c]{\tiny $\delta^\ast=-2\delta_{\ov{m}}$}}
\put(7.8,3){\makebox(0,0)[c]{\tiny $\delta_{\ov{m}}-\delta_{\ov{m-1}}$}}
\put(10.4,0.8){\makebox(0,0)[c]{\tiny $\delta_{\ov{m-1}}-\delta_{\ov{m-2}}$}}
\put(14.5,3){\makebox(0,0)[c]{\tiny $\delta_{\ov{2}}-\delta_{\ov{1}}$}}
\put(17.15,0.8){\makebox(0,0)[c]{\tiny $\delta^\times=\delta_{\ov{1}}-\delta_{1}$}}
\put(19.5,3){\makebox(0,0)[c]{\tiny $\delta_{1}-\delta_{2}$}}
\put(23.8,0.8){\makebox(0,0)[c]{\tiny $\delta_r-\delta_{r+1}$}}
\end{picture}\vskip 6mm

\hskip -3cm \setlength{\unitlength}{0.16in}
\begin{picture}(24,4)
\put(1.8,2){\makebox(0,0)[c]{${\I}_{0}^\times$ : }}
\put(5.6,2){\makebox(0,0)[c]{$\bigcirc$}}
\put(8,2){\makebox(0,0)[c]{$\bigcirc$}}
\put(10.4,2){\makebox(0,0)[c]{$\bigcirc$}}
\put(14.85,2){\makebox(0,0)[c]{$\bigcirc$}}
\put(17.45,2){\makebox(0,0)[c]{$\cdots$}}
\put(8.35,2){\line(1,0){1.5}}
\put(10.82,2){\line(1,0){0.8}}
\put(13.2,2){\line(1,0){1.2}}
\put(15.28,2){\line(1,0){1.45}}
\put(6.8,2){\makebox(0,0)[c]{$\Longrightarrow$}}
\put(12.5,1.95){\makebox(0,0)[c]{$\cdots$}}
\put(5.4,0.8){\makebox(0,0)[c]{\tiny $\delta^\ast=-2\delta_{1}$}}
\put(7.8,3){\makebox(0,0)[c]{\tiny $\delta_{1}-\delta_{2}$}}
\put(10.4,0.8){\makebox(0,0)[c]{\tiny $\delta_{2}-\delta_{3}$}}
\put(15,3){\makebox(0,0)[c]{\tiny $\delta_{r}-\delta_{r+1}$}}
\end{picture}\vskip 5mm
\end{center}
\noindent Here $\delta^\ast$  and $\delta^\times$ also denote distinguished simple roots.\vskip 2mm

Now we define $\td{\mf b}_\infty$, $\ov{\mf b}_\infty$ and ${\mf b}_\infty$ to be the central extensions of $\frak{osp}(V_{\mathcal{I}})$  induced from $\widehat{\gl}(V_{\td{\I}_m})$, where $\mathcal{I}=\td{\I}_m$, $\ov{\I}_m$ and ${\I}_m$  respectively. Similarly, we  define $\td{\mf c}_\infty$, $\ov{\mf c}_\infty$ and ${\mf c}_\infty$  to be the central extensions of $\frak{spo}(V_{\mathcal{I}})$ induced from $\widehat{\gl}(V_{\td{\I}_m})$, where $\mathcal{I}=\td{\I}^\times_m$, $\td{\I}^\times_m$ and $\td{\I}^\times_m$ respectively. For convenience,  the $\Z_2$-graded set $\mathcal{I}$ corresponding to a given Lie superalgebra is summarized in the  table below.\vskip 5mm

\begin{center}
\begin{tabular}{|c||c|c|c|c|c|c|}
\hline \raisebox{-.0ex}{Lie superalgebra}   & \raisebox{-.4ex}{$\td{\mf b}_\infty$} & \raisebox{-.3ex}{$\ov{\mf b}_\infty$} & \raisebox{-.3ex}{${\mf b}_\infty$} & \raisebox{-.2ex}{$\td{\mf c}_\infty$} & \raisebox{-.2ex}{$ \ov{\mf c}_\infty$} & \raisebox{-.2ex}{${\mf c}_\infty$} \\
\hline \raisebox{-.3ex}{$\mathcal{I}$} &  \raisebox{-.4ex}{$\td{\I}_m$} & \raisebox{-.3ex}{$\ov{\I}_m$} & \raisebox{-.2ex}{${\I}_m$} &\raisebox{-.3ex}{$\td{\I}^\times_m$} & \raisebox{-.3ex}{$\ov{\I}^\times_m$} & \raisebox{-.2ex}{${\I}^\times_m$} \\ \hline
\end{tabular}

\vspace{.2cm} Table 1
\end{center}\vskip 3mm
We also assume the following notations:
\begin{itemize}
\item[$\cdot$] $\td{\frak{h}}$, $\ov{\frak{h}}$, ${\frak{h}}$ :  the Cartan subalgebra spanned by $K$ and $E_a$ ($a\in \mathcal{I}^+$),

\item[$\cdot$] $\td{\frak{h}}^\ast$, $\ov{\frak{h}}^\ast$, ${\frak{h}}^\ast$ : the dual Cartan subalgebra spanned by $\Lambda^{\mf x}_0$ and $\delta_a\,(a\in \mathcal{I}^+)$,

\item[$\cdot$] $\td{\Pi}_{\mf x}$, $\ov{\Pi}_{\mf x}$, ${\Pi}_{\mf x}$ : the set of simple roots,


\end{itemize}
of $\td{\mf x}_\infty$, $\ov{\mf x}_\infty$, ${\mf x}_\infty$, respectively.  Here $\Lambda^{\mf x}_0$ is defined by
$\Lambda^{\mf x}_0(K)=1$ and $\Lambda^{\mf x}_0(E_{a})=0$  for $a\in \mathcal{I}^+$. Note that $\ov{\mf x}_\infty, {\mf x}_\infty \subset \td{\mf x}_\infty$ and $\ov{\mf h}^\ast, {\mf h}^\ast\subset \td{\mf h}^\ast$.

\begin{rem}\label{remark on Lambda_0}{\rm
For a positive simple root $\delta$, denote by $\delta^\vee$  the corresponding simple coroot.
We can check that $\Lambda^{\mf x}_0$ is a fundamental weight of ${\mf x}_\infty$ with respect to $\delta^\ast$, that is,  $\Lambda^{\mf x}_0(\delta^\vee)=1$  if $\delta=\delta^\ast$, and $0$ otherwise. Hence $\Lambda^{\mf x}_0$ coincides with the one given in Section \ref{affine Lie algebra}, and  $\delta_{\ov{i}}$ (resp. $\delta_i$) corresponds to  $\widehat{\epsilon}_{m-i+1}$  for $1\leq i\leq m$ (resp. $\widehat{\epsilon}_{i+m}$  for $i\geq 1$).
Also, note that in case of $\ov{\mf x}_\infty$ with $m=0$,  $\Lambda^{\mf x}_0$ is the negative of the fundamental weight of $\ov{\mf x}_\infty$ with respect to $\delta^\ast$.
}
\end{rem}

\subsection{Super duality}\label{review on super duality}
Let us briefly recall super duality in \cite{CLW}. We should remark here that our exposition of super duality deals only  with types ${\mf b}$ and ${\mf c}$, while there are two more cases of type ${\mf d}$ and ${\mf b}^{\bullet}$ in \cite{CLW}.

Fix an arbitrary subset $Y_0$ of $\{\,\delta^\ast,\delta_{\ov{m}}-\delta_{\ov{m-1}},\ldots, \delta_{\ov{2}}-\delta_{\ov{1}}\,\}$. We assume that $Y_0$ is empty when $m=0$. Let $\td{Y}$, $\ov{Y}$ and $Y$ be the union of $Y_0$ and the simple roots which are placed to the right of $\delta^\times$ in the Dynkin diagram of $\td{\mf x}_\infty$, $\ov{\mf x}_\infty$, ${\mf x}_\infty$, respectively, and let $\td{{\mf l}}$, $\ov{{\mf l}}$ and ${\mf l}$ be the standard Levi subalgebras of $\td{\mf x}_\infty$, $\ov{\mf x}_\infty$, ${\mf x}_\infty$ corresponding to $\td{Y}$, $\ov{Y}$ and $Y$, respectively.

Let $\cP_{Y_0}$ be the set of sequences $\lambda=(\lambda_{\ov{m}},\ldots,\lambda_{\ov{1}}\,;\lambda_1,\lambda_2,\ldots)$ such that
\begin{itemize}
\item[(1)] $(\lambda_1,\lambda_2,\ldots) \in \cP$,

\item[(2)] $\lambda_{\ov{m}},\ldots,\lambda_{\ov{1}}\in \mathbb{C}$ and $\left\langle \sum_{a=\ov{m}}^{\ov{1}}\lambda_a\delta_a, \delta^\vee \right\rangle\in\mathbb{Z}_{\geq 0}$ for $\delta\in Y_0$,
\end{itemize}
and put $\lambda^+=(\lambda_1,\lambda_2,\ldots)$ for $\lambda\in \cP_{Y_0}$. We assume that $\lambda=\lambda^+$ when $m=0$. Also for  $\lambda\in\cP$, let $\theta(\lambda)=(\theta(\lambda)_r)_{r\in\hf\Z_{>0}}$ be a sequence determined by
\begin{equation}
\theta(\lambda)_{i-\hf}=\max\{\lambda'_i-i+1,0 \}, \ \ \ \theta(\lambda)_i=\max\{\lambda_i-i,0\}
\end{equation}
for $i\in\Z_{>0}$.
Let
\begin{equation}
\begin{split}
P_{{\mf l}}^+&=\left\{\,\Lambda=\sum_{a=\ov{m}}^{\ov{1}}\lambda_a\delta_a +\sum_{b\in \Z_{>0}}\lambda^+_b\delta_b+ c\Lambda^{\mf x}_0\  \Bigg\vert\ \lambda\in\cP_{Y_0},\, c\in\mathbb{C}\,\right\}\subset \mf{h}^\ast,\\
P_{\td{\mf l}}^+&=\{\,\Lambda^\theta\,|\,\Lambda\in P_{{\mf l}}^+\,\}\subset \td{\mf h}^\ast, \\ P_{\ov{\mf l}}^+&=\{\,\Lambda^\natural\,|\,\Lambda\in P_{{\mf l}}^+\,\}\subset \ov{\mf h}^\ast,
\end{split}
\end{equation}
where
\begin{equation}
\begin{split}
\Lambda^\theta&=\sum_{a=\ov{m}}^{\ov{1}}\lambda_a\delta_a +\sum_{b\in \hf\Z_{>0}}\theta(\lambda^+)_b\delta_b+ c\Lambda^{\mf x}_0,\\
\Lambda^\natural&= \sum_{a=\ov{m}}^{\ov{1}}\lambda_a\delta_a +\sum_{b\in \hf+\Z_{\geq 0}}(\lambda^+)'_b\delta_b+ c\Lambda^{\mf x}_0.
\end{split}
\end{equation}

For $\Lambda\in P_{{\mf l}}^+$, let  $L({\mf l},\Lambda)$ be the irreducible ${\mf l}$-module with highest weight $\Lambda$, and  $L({\mf x}_\infty,\Lambda)$  the irreducible quotient of  $K(\mf{x}_\infty,\Lambda):={\rm Ind}^{{\mf x}_\infty}_{{\mf p}}L({\mf l},\Lambda)$, where ${\mf b}$ is the Boreal subalgebra spanned by $K$ and upper triangular matrices, ${\mf p}={\mf b}+{\mf l}$ and $L({\mf l},\Lambda)$ is extended to a ${\mf p}$-module in a standard way.  We define  $L(\td{\mf l},\Lambda^\theta)$, $K(\td{\mf{x}}_\infty,\Lambda^\theta)$, $L(\td{\mf x}_\infty,\Lambda^\theta)$ and $L(\ov{\mf l},\Lambda^\natural)$, $K(\ov{\mf{x}}_\infty,\Lambda^\natural)$, $L(\ov{\mf x}_\infty,\Lambda^\natural)$ in the same way.

Let $\mathcal{O}$ be the category of ${\mf x}_\infty$-modules $M$ satisfying
\begin{itemize}
\item[(1)] $M=\bigoplus_{\gamma\in {\mf h}^\ast}M_\gamma$, where $M_\gamma=\{\,m\,|\,h\cdot m=\gamma(h)m\ (h\in {\mf h})\,\}$ and $\dim M_\gamma<\infty$,

\item[(2)] $M$ decomposes into direct sum of $L({\mf l},\Lambda)$ for $\Lambda\in P_{\mf l}^+$,

\item[(3)] there exist $\Lambda_1,\ldots,\Lambda_\ell\in P_{\mf l}^+$ such that  if $M_\gamma\neq \{0\}$, then $\gamma\in \Lambda_i-\sum_{\alpha\in \Pi_{\mf x}}c_\alpha\alpha$ for some $i$ and $c_\alpha\in\Z_{\geq 0}$.
\end{itemize}
Note that it is a parabolic analogue of  Berstein-Gelfand-Gelfand category. In a similar way, we define the categories $\td{\mathcal{O}}$ and $\ov{\mathcal{O}}$ of modules over $\td{\mf x}_\infty$ and $\ov{\mf x}_\infty$, respectively.

Now, for $M=\bigoplus_{\gamma\in {\td{\mf h}}^\ast}M_\gamma$ in $\td{\mathcal{O}}$, we define
\begin{equation}
\begin{split}
{T}(M)=\bigoplus_{\gamma\in {{\mf h}}^\ast}M_\gamma,\ \ \ \
\ov{T}(M)=\bigoplus_{\gamma\in {\ov{\mf h}}^\ast}M_\gamma.\\
\end{split}
\end{equation}
Note that $T(M), \ov{T}(M)\subset M$. If $f : M \rightarrow N$ is a homomorphism of $\td{\mf x}_\infty$-modules for $M,N\in \td{\mathcal{O}}$, then $T$ and $\ov{T}$ naturally induce homomorphisms $T(f) : T(M) \rightarrow T(N)$ of ${\mf x}_\infty$-modules and $\ov{T}(f) : \ov{T}(M) \rightarrow \ov{T}(N)$ of $\ov{\mf x}_\infty$-modules, respectively. Hence    $T$ and $\ov{T}$ define functors
\begin{equation}
\begin{CD}
\mathcal{O} @  < T <<    \td{\mathcal{O}} @ >\ov{T} >> \ov{\mathcal{O}}.
\end{CD}
\end{equation}\vskip 3mm

\begin{thm}[Theorems 4.6 and 5.4 in \cite{CLW}]\label{super duality} \mbox{}
\begin{itemize}
\item[(1)] For $\Lambda\in P_{\mf l}^+$, we have
\begin{equation*}
\begin{split}
&   \begin{cases} \ \ \ T(L(\td{\mf l},\Lambda^\theta))=L({\mf l},\Lambda),\\ T(K(\td{\mf{x}}_\infty,\Lambda^\theta))=K(\mf{x}_\infty,\Lambda),\\ \ \! T(L(\td{\mf x}_\infty,\Lambda^\theta))=L({\mf x}_\infty,\Lambda) \end{cases}\ \
 \begin{cases} \ \ \ \ov{T}(L(\td{\mf l},\Lambda^\theta))=L(\ov{\mf l},\Lambda^\natural),\\ \ov{T}(K(\td{\mf{x}}_\infty,\Lambda^\theta))=K(\ov{\mf{x}}_\infty,\Lambda^\natural), \\ \ \! \ov{T}(L(\td{\mf x}_\infty,\Lambda^\theta))=L(\ov{\mf x}_\infty,\Lambda^\natural).\end{cases}
\end{split}
\end{equation*}

\item[(2)] $T$ and $\ov{T}$ are equivalences of categories.  Hence  the categories $\mathcal{O}$ and $\ov{\mathcal{O}}$ are equivalent.
\end{itemize}
\end{thm}\vskip 2mm

\begin{rem}{\rm
The equivalence between $\mathcal{O}$ and $\ov{\mathcal{O}}$ is called {\it super duality}, while $\td{\mathcal{O}}$ is used as an intermediate category to establish it. But in later arguments, we will consider a simple module in $\td{\mathcal{O}}$ together with one in $\ov{\mathcal{O}}$ since both  have a meaning for a representation theoretical interpretation of $S^{\mf x}_{(\lambda,n)}({\bf x}_\A)$.}
\end{rem}

\begin{rem}{\rm
For each $M$ in $\td{\mathcal{O}}$, $\ov{\mathcal{O}}$ and $\mathcal{O}$,  we define the character of $M$ to be the trace of the operator $\prod_{a\in \mathcal{I}^+}q^Kx_a^{E_{a}}$. By Theorem \ref{super duality}, if  for $\Lambda\in P_{\mf l}^+$
\begin{equation*}
{\rm ch}L({\mf x}_\infty,\Lambda)=\sum_{\Lambda'}a_{\Lambda,\Lambda'}{\rm ch}K(\mf{x}_\infty,\Lambda'),
\end{equation*}
for some $a_{\Lambda,\Lambda'}$, then we have
\begin{equation*}
\begin{split}
{\rm ch}L(\td{\mf x}_\infty,\Lambda^\theta)&=\sum_{\Lambda'}a_{\Lambda,\Lambda'}{\rm ch}K(\td{\mf{x}}_\infty,(\Lambda')^\theta),\\
{\rm ch}L(\ov{\mf x}_\infty,\Lambda^\sharp)&=\sum_{\Lambda'}a_{\Lambda,\Lambda'}{\rm ch}K(\ov{\mf{x}}_\infty,(\Lambda')^\sharp).
\end{split}
\end{equation*}
Hence,  the characters of $L(\td{\mf x}_\infty,\Lambda^\theta)$ and $L(\ov{\mf x}_\infty,\Lambda^\natural)$ follow immediately from  a Kazhdan-Lusztig type formula of $L({\mf x}_\infty,\Lambda)$ (see \cite[Theorem 4.8 and Remark 4.9]{CLW}).}
\end{rem}

\subsection{Combinatorial character formula}
Let us write $\lambda=(\lambda_{\ov{m}},\ldots,\lambda_{\ov{1}}\,;\lambda_1,\lambda_2,\ldots)$ for $(\lambda,n)\in\cP({\mf x})$, and put
\begin{equation}\label{Lambda lambda n-2}
\La^{\mf x}({\lambda,n})=
\sum_{a=\ov{m}}^{\ov{1}}\lambda_a\delta_a +\sum_{b\in \Z_{>0}}\lambda_b\delta_b+ n\Lambda^{\mf x}_0,
\end{equation}
which is a dominant integral weight for ${\mf x}_\infty$. Note that it coincides with our previous definition in (\ref{Lambda lambda n}) (see Remark \ref{remark on Lambda_0}). Let $\mathcal{O}_{\rm int}$ be the full subcategory of $\mathcal{O}$, whose objects are integrable ${\mf x}_\infty$-modules. It is well-known that $\mathcal{O}_{\rm int}$ is a semisimple tensor category and every object in $\mathcal{O}_{\rm int}$ is isomorphic to a direct sum of $L({\mf x}_\infty,\Lambda^{\mf x}(\lambda,n))$ with finite multiplicity for each $(\lambda,n)\in \cP({\mf x})$.

Let ${\mathcal{N}}=\{\,\ov{m},\ldots,\ov{1}\,\}\cup \{\,1,2,3,\ldots\,\,\,\}$. We have
\begin{equation}
{\rm ch}L({{\mf x}}_\infty,{\Lambda}^{\mf x}(\lambda,n))=S^{\mf x}_{(\lambda,n)}({\bf x}_{\mathcal{N}}),
\end{equation}
 for $(\lambda,n)\in\cP({\mf x})$  by Theorem \ref{main result}, and
\begin{equation}\label{LR for x}
\begin{split}
L({{\mf x}}_\infty,{\Lambda}^{\mf x}(\mu,m))  \otimes L(&{{\mf x}}_\infty,{\Lambda}^{\mf x}(\nu,n))\\ &\simeq\bigoplus_{(\lambda,m+n)}L({{\mf x}}_\infty,{\Lambda}^{\mf x}(\lambda,m+n))^{\oplus c^{(\lambda,m+n)}_{(\mu,m) (\nu,n)}({\mf x}) },
\end{split}
\end{equation}
for $(\mu,m), (\nu,n)\in\cP({\mf x})$ by Theorem \ref{LR rule}.
\qed\vskip 3mm

Let
\begin{equation}
\begin{split}
\td{\mathcal{N}}&=\{\,\ov{m},\ldots,\ov{1}\,\}\cup \{\,\tfrac{1}{2},1,\tfrac{3}{2},\ldots\,\,\,\}, \\
\ov{\mathcal{N}}&=\{\,\ov{m},\ldots,\ov{1}\,\}\cup \{\,\tfrac{1}{2},\tfrac{3}{2},\tfrac{5}{2},\ldots\,\,\,\}.\\
\end{split}
\end{equation}
Then we obtain the following combinatorial character formula, which is the main result in this section.
\begin{thm} \label{application to super duality}
For $(\lambda,n)\in\cP({\mf x})$, we have
\begin{equation*}
\begin{split}
&{\rm ch}L(\td{{\mf x}}_\infty,{\Lambda}^{\mf x}(\lambda,n)^\theta)=  S^{\mf x}_{(\lambda,n)}({\bf x}_{\td{\mathcal{N}}}),\\
&{\rm ch}L(\ov{{\mf x}}_\infty,{\Lambda}^{\mf x}(\lambda,n)^\natural)= S^{\mf x}_{(\lambda,n)}({\bf x}_{\ov{\mathcal{N}}}).
\end{split}
\end{equation*}
\end{thm}
\pf We keep the notations in Section \ref{review on super duality} and suppose that $Y_0=\{\,\delta_{\ov{m}}-\delta_{\ov{m-1}},\ldots, \delta_{\ov{2}}-\delta_{\ov{1}}\,\}$.

Let $(\lambda,n)\in\cP({\mf x})$ be given. Note that $L({{\mf x}}_\infty,{\Lambda}^{\mf x}(\lambda,n))$ (resp. $L({\td {\mf x}}_\infty,{\Lambda}^{\mf x}(\lambda,n)^\theta)$) is completely reducible as an ${\mf l}$-module (resp. $\td{\mf l}$-module).  By Theorem \ref{super duality} (1), if $${\rm ch}L({{\mf x}}_\infty,{\Lambda}^{\mf x}(\lambda,n))=\sum_{\Lambda \in P^+_{\mf l}} a_{\Lambda}{\rm ch}L({\mf l},\Lambda),$$ for some $a_{\Lambda}\in \mathbb{Z}_{\geq 0}$, then $${\rm ch}L({\td{\mf x}}_\infty,{\Lambda}^{\mf x}(\lambda,n)^\theta)=\sum_{\Lambda \in P^+_{\mf l}} a_{\Lambda}{\rm ch}L(\td{\mf l},\Lambda^\theta).$$
Note that
\begin{equation*}
\begin{split}
{\rm ch}L({{\mf x}}_\infty,{\Lambda}^{\mf x}(\lambda,n))&=S^{\mf x}_{(\lambda,n)}({\bf x}_{{\mathcal{N}}})=q^n\sum_{\sigma\in\cP}c^{\,\sigma}_{(\lambda,n)}({\mf x})  S_{\sigma}({\bf x}_{{\mathcal{N}}}),\\
&=q^n\sum_{\sigma\in\cP}c^{\,\sigma}_{(\lambda,n)}({\mf x})\sum_{\tau\subset \sigma} S_{\tau}({\bf x}_{\{\ov{m},\ldots,\ov{1}\}}) S_{\sigma/\tau}({\bf x}_{\mathbb{Z}_{>0}})\\
&=q^n\sum_{\sigma\in\cP}c^{\,\sigma}_{(\lambda,n)}({\mf x})\sum_{\tau,\gamma} c^\sigma_{\tau \gamma}S_{\tau}({\bf x}_{\{\ov{m},\ldots,\ov{1}\}}) S_{\gamma}({\bf x}_{\mathbb{Z}_{>0}}),
\end{split}
\end{equation*}
Since ${\rm ch}L({\mf l},\Lambda)=q^nS_{\tau}({\bf x}_{\{\ov{m},\ldots,\ov{1}\}}) S_{\gamma}({\bf x}_{\mathbb{Z}_{>0}})$ with $\Lambda=\sum_{a=1}^m\tau_{m-a+1}\delta_{\ov{a}}+\sum_{b\in\mathbb{Z}_{>0}}\gamma_b\delta_b+n\Lambda^{\mf x}_0$ and $${\rm ch}L(\td{\mf l},\Lambda^\theta)=q^nS_{\tau}({\bf x}_{\{\ov{m},\ldots,\ov{1}\}}) S_{\gamma}({\bf x}_{\hf\mathbb{Z}_{>0}}),$$
we have
\begin{equation*}
\begin{split}
{\rm ch}L({\td {\mf x}}_\infty,{\Lambda}^{\mf x}(\lambda,n)^\theta)&=q^n\sum_{\sigma\in\cP}c^{\,\sigma}_{(\lambda,n)}({\mf x})\sum_{\tau,\gamma} c^\sigma_{\tau \gamma}S_{\tau}({\bf x}_{\{\ov{m},\ldots,\ov{1}\}}) S_{\gamma}({\bf x}_{\hf\mathbb{Z}_{>0}}) \\
&=q^n\sum_{\sigma\in\cP}c^{\,\sigma}_{(\lambda,n)}({\mf x})\sum_{\tau\subset \sigma} S_{\tau}({\bf x}_{\{\ov{m},\ldots,\ov{1}\}}) S_{\sigma/\tau}({\bf x}_{\hf\mathbb{Z}_{>0}})\\
&=q^n\sum_{\sigma\in\cP}c^{\,\sigma}_{(\lambda,n)}({\mf x})  S_{\sigma}({\bf x}_{\td{\mathcal{N}}})=S^{\mf x}_{(\lambda,n)}({\bf x}_{\td{\mathcal{N}}}).
\end{split}
\end{equation*}
By the same arguments, we also have  ${\rm ch}L(\ov{{\mf x}}_\infty,{\Lambda}^{\mf x}(\lambda,n)^\natural)= S^{\mf x}_{(\lambda,n)}({\bf x}_{\ov{\mathcal{N}}})$.

\qed

\begin{cor}\label{semisimplicity of O}  For $(\mu,m), (\nu,n)\in\cP({\mf x})$, we have

\begin{equation*}
\begin{split}
L({\td{\mf x}}_\infty,{\Lambda}^{\mf x}(\mu,m)^\theta)\otimes L(&{\td{\mf x}}_\infty,{\Lambda}^{\mf x}(\nu,n)^\theta)\\ &\simeq\bigoplus_{(\lambda,m+n)}L({\td{\mf x}}_\infty,{\Lambda}^{\mf x}(\lambda,m+n)^\theta)^{\oplus c^{(\lambda,m+n)}_{(\mu,m) (\nu,n)}({\mf x}) },\\
L({\ov{\mf x}}_\infty,{\Lambda}^{\mf x}(\mu,m)^\natural)\otimes L(&{\ov{\mf x}}_\infty,{\Lambda}^{\mf x} (\nu,n)^\natural)\\ &\simeq\bigoplus_{(\lambda,m+n)}L({\ov{\mf x}}_\infty,{\Lambda}^{\mf x}(\lambda,m+n)^\natural)^{\oplus c^{(\lambda,m+n)}_{(\mu,m) (\nu,n)}({\mf x}) }.
\end{split}
\end{equation*}
\end{cor}
\pf It suffices to prove the first identity since the proof of the second one is identical. Write $M=L({\td{\mf x}}_\infty,{\Lambda}^{\mf x}(\mu,m)^\theta)$ and $N=L({\td{\mf x}}_\infty,{\Lambda}^{\mf x}(\nu,n)^\theta)$. Assume that
\begin{equation*}
\begin{split}
{\rm ch}M&=q^m\sum_{\sigma}a_\sigma S_{\sigma}({\bf x}_{\td{\mathcal{N}}}), \ \ {\rm ch}N =q^n\sum_{\tau}b_\sigma S_{\tau}({\bf x}_{\td{\mathcal{N}}}).
\end{split}
\end{equation*}
Then we have
\begin{equation*}
\begin{split}
&{\rm ch}\left(M\otimes N\right)=q^{m+n}\sum_{\sigma, \tau }a_\sigma b_\tau\sum_{\gamma}c^\gamma_{\sigma\tau}S_{\gamma}({\bf x}_{\td{\mathcal{N}}})\\
\end{split}
\end{equation*}
By the same arguments as in Theorem \ref{application to super duality}, we have
\begin{equation}\label{tensor}
\begin{split}
&{\rm ch}\,{T}\left( M\otimes N\right)=q^{m+n}\sum_{\sigma, \tau }a_\sigma b_\tau\sum_{\gamma}c^\gamma_{\sigma\tau}S_{\gamma}({\bf x}_{{\mathcal{N}}})\\
\end{split}
\end{equation}
Since
\begin{equation*}
\begin{split}
{\rm ch}\,{T}\left( M\right)={\rm ch}\left(L({{\mf x}}_\infty,{\Lambda}^{\mf x}(\mu,m))\right)=q^m\sum_{\sigma}a_\sigma S_{\sigma}({\bf x}_{{\mathcal{N}}}),\\
{\rm ch}\,{T}\left( N \right)={\rm ch}\left(L({{\mf x}}_\infty,{\Lambda}^{\mf x}(\nu,n))\right)=q^n\sum_{\tau}b_\tau S_{\tau}({\bf x}_{{\mathcal{N}}}),\\
\end{split}
\end{equation*}
we have
\begin{equation}\label{equality of ch}
{\rm ch}\,{T}\left( M\otimes N\right)={\rm ch}\,{T}\left( M\right) {\rm ch}\,{T}\left( N\right)={\rm ch}\left({T}( M) \otimes {T}(N)\right)
\end{equation}
On the other hand, by definition of ${T}$, we have
${T}( M) \otimes {T}(N) \subset {T}\left( M\otimes N\right)$.
Therefore, (\ref{equality of ch}) implies that ${T}( M) \otimes{T}(N) ={T}\left( M\otimes N\right)$.
Finally, (\ref{LR for x}) combined with the equivalence of ${T}$ implies that
\begin{equation*}
M\otimes N \simeq \bigoplus_{(\lambda,m+n)}L({\td{\mf x}}_\infty,{\Lambda}^{\mf x}(\lambda,m+n)^\theta)^{\oplus c^{(\lambda,m+n)}_{(\mu,m) (\nu,n)}({\mf x}) }.
\end{equation*}
\qed

\begin{rem}{\rm
Let $\A$ be an arbitrary $\Z_2$-graded linearly ordered set with infinitely many odd elements. Then $S^{\mf x}_{(\lambda,n)}({\bf x}_{\A})$ is equal to one of the characters $S^{\mf x}_{(\lambda,n)}({\bf x}_{\td{\mathcal{N}}})$
and $S^{\mf x}_{(\lambda,n)}({\bf x}_{\ov{\mathcal{N}}})$ since $S^{\mf x}_{(\lambda,n)}({\bf x}_{\A})$ does not depend on the choice of a linear ordering on $\A$.
When $\A$ is a finite set, $S^{\mf x}_{(\lambda,n)}({\bf x}_{\A})$ is also an irreducible character of a finite dimensional orthosymplectic Lie superalgebra. This can be deduced by applying a truncation functor on $\ov{\mathcal{O}}$ introduced in \cite{CLW}. But we would like to remark that the irreducible representation associated to $S^{\mf x}_{(\lambda,n)}({\bf x}_{\A})$ is infinite dimensional called an oscillator representation, which was also studied in \cite{CKW} using Howe duality.}
\end{rem}

\section{Proof of Theorem \ref{main result}  }\label{Proof for BC}
\subsection{RSK correspondence}\label{MAB}
Let $\A$ and $\cB$ be two linearly ordered sets. Let $\M_{\A,\cB}$ be the set of $\A\times \cB$ matrices $A=(a_{i,j})$ with entries in $\Z_{\geq 0}$ such that $\sum_{i\in\A,j\in\cB}a_{i,j}<\infty$. Let $\Omega_{\A,\cB}$ be the set of biwords
$(\bi,\bj)\in \W_{\A}\times \W_{\cB}$ with $\bi=i_1\cdots i_r$ and $\bj=j_1\cdots j_r$  for some $r\geq 0$ such that $(i_1,j_1)\leq \cdots \leq (i_r,j_r)$. Here, we assume that $(i,j)< (k,l)$ if and only if $(j< l)$ or $(j=l\ \text{and} \ i>k)$ for $(i,j)$ and $(k,l)\in \A\times \cB$.
Then the map from $\Omega_{\A,\cB}$ to $\M_{\A,\cB}$ sending
$(\bi,\bj)$ to $A(\bi,\bj)=(a_{i,j})$ with
$a_{i,j}=|\{\,k\,|\,(i_k,j_k)=(i,j) \,\}|$ is a bijection. Note that the pair of
empty words $(\emptyset,\emptyset)$ corresponds to zero matrix $\mathbb{O}$.

For $A\in \M_{\A,\cB}$, suppose that $A=A(\bi,\bj)$ and $A^t=A(\bk,\bl)$ for $(\bi,\bj)\in\Omega_{\A,\cB}$ and $(\bk,\bl)\in\Omega_{\cB,\A}$.
Then we have a bijection called the Robinson-Schensted-Knuth (simply RSK) correspondence
\begin{equation}\label{Knuth}
\begin{split}
\M_{\A,\cB}& \stackrel{\kappa}{\longrightarrow} \bigsqcup_{\lambda\in\cP} SST_{\A}(\lambda)\times  SST_\cB(\lambda),\\
\end{split}
\end{equation}
where $\kappa(A)=((\bi\rightarrow \emptyset),(\bk\rightarrow \emptyset))$ \cite{Kn}.

\subsection{RSK map as a  $\gl_\infty$-crystal isomorphism}\label{gl crystal on MAB}
Let us recall a $\gl_\infty$-crystal structure associated with the RSK correspondence  \cite{K09} (with a slight modification of notations).

Let $\M=\M_{\Z_{>0},\Z_{<0}^\vee}$ and $\Omega=\Omega_{\Z_{>0},\Z_{<0}^\vee}$ where $\Z_{<0}^\vee=\{\,-k^\vee\,|\,k\in\Z_{>0}\,\}$ is the dual $\gl_{<0}$-crystal of $\Z_{<0}$ with a linear ordering given by $-1^\vee<-2^\vee<\cdots$. Note that $\M$ is a normal $\gl_{>0}$-crystal, where $\widetilde{x}_iA=A(\widetilde{x}_i\bi,\bj)$ for $A\in\M$ with $A=A(\bi,\bj)$
($x=e,f$ and $i\in \Z_{>0}$). Here, we assume that $\widetilde{x}_i
A={\bf 0}$ if $\widetilde{x}_i\bi={\bf 0}$. In a similar way, we may view $\M$ as a normal $\gl_{<0}$-crystal with respect to $\te_i$, $\tf_i$ ($i\in\Z_{<0}$) by considering the transpose of $A\in\M$. (Note that $\M$ was defined as $\M_{\Z_{<0}^\vee,\Z_{>0}}$ in \cite{K09}.) Then $\M$ is a $(\gl_{<0},\gl_{>0})$-bicrystal, that is, $\te_i$, $\tf_i$ ($i\in\Z_{>0}$) commute with $\te_j$, $\tf_j$ ($j\in\Z_{<0}$).

Let $\mathcal{T}=\bigsqcup_{\lambda\in\cP}SST_{\Z_{>0}}(\lambda)\times SST_{\Z_{<0}^\vee}(\lambda)$. Then $\mathcal{T}$ is clearly a $(\gl_{<0},\gl_{>0})$-bicrystal, and $\kappa$ is a  $(\gl_{<0},\gl_{>0})$-bicrystal isomorphism \cite{DK,La}.

Now, for $A=(a_{i, -j^\vee})\in\M$, we define
\begin{equation}
\begin{split}
&\te_0 A=
\begin{cases}
A-E_{1,-1^\vee}, & \text{if $a_{1,-1^\vee}\neq 0$}, \\
{\bf 0}, & \text{otherwise},
\end{cases}\ \ \ \text{and} \ \ \
\tf_0 A= A+E_{1,-1^\vee},
\end{split}
\end{equation}
where $E_{i,-j^\vee}\in \M$ denotes the elementary matrix with $1$ at
the position  $(i,-j^\vee)$. Put
\begin{equation}
\begin{split}
{\rm wt}(A)&=\sum_{i,j\geq 1} a_{i,-j^\vee}(\epsilon_{i}-\epsilon_{-j}),\\
\varepsilon_0(A)&=a_{1,-1^\vee},\\
\varphi_0(A)&=\langle {\rm
wt}(A),h_0 \rangle + \varepsilon_0(A).
\end{split}
\end{equation}
Then
$\M$ is a connected $\gl_{\infty}$-crystal with a unique highest weight element
$\mathbb{O}$ \cite[Proposition 3.1]{K09}.

Next, let us describe $\te_0$ and $\tf_0$ on $\mathcal{T}$. Let $(S,T)\in \mathcal{T}$ be given. For $k\geq 1$, let $s_k$ and $t_k$ be the entries in the top of the $k$-th column of $S$ and $T$, respectively. We assign
\begin{equation}\label{signs}
\sigma_k=
\begin{cases}
+ \ , & \text{if the $k$-th column is empty}, \\
+ \ ,& \text{if $s_k>-1^\vee$ and $t_k>1$}, \\
- \ ,& \text{if $s_k=-1^\vee$ and $t_k=1$,}\\
\ \cdot \ \, ,& \text{otherwise}.
\end{cases}
\end{equation}
Let $\td{\sigma}$ be the sequence obtained from $\sigma=(\ldots,\sigma_2,\sigma_1)$ by the same method as in Section \ref{crystal Txlambda}.

We define $\te_0 (S,T)$ to be the bitableaux
obtained from $(S,T)$  by removing $\boxed{1}$ and $\boxed{\tiny{-1^\vee}}$  in the columns of $S$ and $T$ corresponding to the left-most $-$ in $\td{\sigma}$.
If there is no such $-$ sign, then we define $\te_0 (S,T)={\bf
0}$. We define $\tf_0 (S,T)$ to be the
bitableaux obtained from $(S,T)$ by adding $\boxed{1}$ and $\boxed{-1^\vee}$  on top of the
columns of $S$ and $T$ corresponding to the right-most $+$ in $\td{\sigma}$. If there is no such $+$ sign, then we
define $\tf_0 (S,T)={\bf 0}$.

Put
\begin{equation}
\begin{split}
{\rm wt}(S,T)&=
\sum_{i\geq 1}\left(m_i(S)\epsilon_i-m_{-i^\vee}(T)\epsilon_{-i}\right),\\
{\varepsilon}_0(S,T)&=\max\{\,k\,|\,\td{e}_0^k(S,T)\neq {\bf 0}\,\},\\
{\varphi}_0(S,T)&=\langle {\rm wt}(S,T), h_0\rangle+ {\varepsilon}_0(S,T).
\end{split}
\end{equation}
Then $\mathcal{T}$ is also a connected $\gl_{\infty}$-crystal with a unique highest weight element
$(\emptyset,\emptyset)$ \cite[Proposition 3.5]{K09}.

\begin{prop}[Theorem 3.6 in \cite{K09}]\label{KnuthIso}
$\kappa : \M \longrightarrow \mathcal{T}$ is a $\gl_\infty$-crystal isomorphism.
\end{prop}

Given  $\mu,\nu\in \cP$, we put
\begin{equation}
\begin{split}
\M({\mu,\nu})&=\M \times
 SST_{\Z_{>0}}(\mu)\times SST_{\Z_{<0}}(\nu)^\vee\\
\end{split}
\end{equation}
where $ SST_{\Z_{<0}}(\nu)^\vee$ is the dual $\gl_{<0}$-crystal of $ SST_{\Z_{<0}}(\nu)$.
For ${\bf m}=(A,S,T^\vee)\in\M(\mu,\nu)$ and $i\in\Z$, we define
\begin{equation}
\td{x}_i{\bf m}=
\begin{cases}
(A',U,T^\vee), & \text{if $i>0$ and $\td{x}_i(A
\otimes S )=A'\otimes U$,} \\
(A'',S,V^\vee), & \text{if $i<0$ and
$\td{x}_i(A\otimes T^\vee)=A''\otimes V^\vee$},\\
(\td{x}_0A,S,T^\vee), & \text{if $i=0$},
\end{cases}
\end{equation}
where $x=e,f$, and $\widetilde{x}_i(A,S,T^\vee)={\bf 0}$ if any
of its components is ${\bf 0}$. Put
\begin{equation}
\begin{split}
{\rm wt}({\bf m})&={\rm wt}(A)+{\rm wt}(S)+{\rm wt}(T^\vee),\\
\varepsilon_i({\bf m})&=\max\{\,k\,|\,\te^k_i{\bf m}\neq {\bf 0}\,\},\\
\varphi_i({\bf m})&=\langle {\rm wt}({\bf m}),h_i\rangle-\varepsilon_i({\bf m}),
\end{split}
\end{equation}
for $i\in\Z$. Then $\M(\mu,\nu)$ is a connected
$\gl_{\infty}$-crystal with a unique highest weight element
$\mathbb{O}_{\mu,\nu}=(\mathbb{O},H_{\mu},H_{\nu}^\vee)\in
\M(\mu,\nu)$
\cite[Proposition 4.5]{K09}.

For $n\geq \ell(\mu)+\ell(\nu)$, let
\begin{equation}
\begin{split}
\Lambda(\mu,\nu,n)&=n\Lambda_0+\sum_{i\geq 1}\mu_i\epsilon_{i}-\sum_{j\geq 1}\nu_j\epsilon_{-j}
\\
&=(n-\mu_1-\nu_1)\Lambda_0+\sum_{i=1}^{\mu_1}\Lambda_{\mu'_i}+\sum_{j=1}^{\nu_1}\Lambda_{-\nu'_j}.
\end{split}
\end{equation}
Then $\M(\mu,\nu)$ can be viewed as a crystal of the generalized Verma module in the following sense.
\begin{prop}[Proposition 4.6 in \cite{K09}]\label{Psilambdamu}
For $n\geq \ell(\mu)+\ell(\nu)$, there exists a unique embedding
of a $\gl_{\infty}$-crystal
\begin{equation*}
\Psi_{\mu,\nu,n} :
\B(\gl_{\infty},{\Lambda(\mu,\nu,n)}) \longrightarrow
\M(\mu,\nu)\otimes T_{n\Lambda_0}
\end{equation*}
sending $u_{\Lambda(\mu,\nu,n)}$ to $\mathbb{O}_{\mu,\nu}\otimes
t_{n\Lambda_0}$ and commuting with $\te_i$ $(i\in \Z)$.
\end{prop}

\begin{rem}{\rm
The embedding in Proposition \ref{Psilambdamu} is given by applying the Sagan and Stanley's RSK algorithm for  skew tableaux to $\B(\gl_{\infty},{\Lambda(\mu,\nu,n)})$, which is realized in terms of pairs of skew tableaux in \cite{K09}.
}
\end{rem}

\subsection{An ${\mf x}_\infty$-analogue of RSK correspondence}\label{subalgebras}
Let us consider an ${\mf x}_\infty$-analogue of Proposition \ref{KnuthIso}.
Let
\begin{equation}\label{hatMlambda}
\begin{split}
{\M}^{\mf x}& =\{\,A\in\M\ |\ a_{i,-j^\vee}=a_{j,-i^\vee}, \
  \text{$\epsilon\big\vert a_{i,-i^\vee}$
($i,j\in\Z_{>0}$)} \,\} \subset \M, \\
{\M}^{\mf x}(\lambda)& =\{\,(A,T,T^\ast)\ |\,A\in{\M}^{\mf x}, \ T\in  SST_{\Z_{>0}}(\lambda) \,\}\subset \M({\lambda,\lambda}),
\end{split}
\end{equation}
for $\lambda\in\cP$.
Here $T^\ast = \tf_{-i_1}\cdots\tf_{-i_r}H_\lambda^\vee\in  SST_{\Z_{<0}}(\lambda)^\vee$ when $T = \tf_{i_1}\cdots\tf_{i_r}H_\lambda$. It is not difficult to see that $T^\ast\neq {\bf 0}$ if and only if $T\neq {\bf 0}$. For $i\in\mathbb{Z}_{\geq 0}$, let
\begin{equation}\label{folded operators}
\begin{split}
&\td{E}_0 =(\te_0)^\epsilon, \ \ \td{F}_0=(\tf_0)^\epsilon, \\
&\td{E}_i =\te_i\te_{-i},\ \ \td{F}_i =\tf_i\tf_{-i} \ \
(i> 0).
\end{split}
\end{equation}
Then we have the following, whose proof is almost the same in \cite[Proposition 5.14]{K09}.
\begin{lem}\label{Connectivity for MhatLambda}
$ {\M}^{\mf x}(\lambda)\cup\{{\bf 0}\}$ is invariant under $\td{E}_i$ and $\td{F}_i$ for $i\in\Z_{\geq 0}$. In particular,
$$ {\M}^{\mf x}(\lambda)=\{\,\td{F}_{i_1}\cdots\td{F}_{i_r}{\mathbb{O}}_{\lambda,\lambda}\,|\,r\geq
0,\  i_1,\ldots,i_r\in\mathbb{Z}_{\geq 0}\,\}\setminus\{{\bf 0}\} .$$
Hence, $ {\M}^{\mf x}(\lambda)$ is an ${\mf x}_\infty$-crystal with respect to ${\rm wt}^{\mf x}$,  ${\varepsilon}^{\mf x}_i$, ${\varphi}^{\mf x}_i$, $\td{E}_i$, $\td{F}_i$  $(i\in\Z_{\geq 0})$, where for ${\bf m}=(A,T,T^\ast)\in{\M}^{\mf x}(\lambda)$,
\begin{equation}
\begin{split}
&{\rm wt}^{\mf x}({\bf m})={\rm wt}({\bf m})\in {P}_{\mf x},\\
&{\varepsilon}^{\mf x}_i({\bf m})=\varepsilon_i({\bf m}),\ \
{\varphi}^{\mf x}_i({\bf m})=\varphi_i({\bf m}) \ \ (i>0),\\
&{\varepsilon}^{\mf x}_0({\bf m})=\tfrac{1}{\epsilon}\varepsilon_0({\bf m}), \ \
{\varphi}^{\mf x}_0({\bf m})=\tfrac{1}{\epsilon}\varphi_0({\bf m}).
\end{split}
\end{equation}
\end{lem}

Since the map sending $(A,T,T^*)$ to $(A,T)$ is a bijection
\begin{equation}
 {\M}^{\mf x}(\lambda)\longrightarrow{\M}^{\mf x}\times  SST_{\Z_{>0}}(\lambda),
\end{equation}
we identify $ {\M}^{\mf x}(\lambda)$ with ${\M}^{\mf x}\times  SST_{\Z_{>0}}(\lambda)$ as an ${\mf x}_\infty$-crystal, where $\td{E}_i$ and $\td{F}_i$ ($i>0$) act on $ SST_{\Z_{>0}}(\lambda)$ of the righthand side in the same way as $\te_i$ and $\tf_i$ ($i>0$), and $\td{E}_0$, $\td{F}_0$ act on ${\M}^{\mf x}$.

Let ${\bf m}=(A,T)\in {\M}^{\mf x}({\lambda})$ be given. Note that  $\kappa(A)=(S,S^\ast)$ for some $S\in SST_{\Z_{>0}}(\tau)$,
 with  $\tau\in \cP_{\mf x}$  (see \cite{Kn}).
Hence the map sending ${\bf m}=(A,T)$ to $(S,T)$ gives a weight preserving bijection
\begin{equation}
\kappa^{\mf x}_{\lambda}\, : \, {\M}^{\mf x}({\lambda}) \longrightarrow {\mathcal{T}}^{\mf x}(\lambda).
\end{equation}
Since $\kappa$ is a $(\gl_{<0},\gl_{>0})$-bicrystal isomorphism, it is not difficult to see that $\kappa^{\mf x}_{(0)}$ is a $\gl_{>0}$-crystal isomorphism with respect to $\td{E}_i$ and $\td{F}_i$ ($i\in\Z_{>0}$), and hence so is $\kappa^{\mf x}_\lambda$.

\begin{prop}\label{kappa x}
$\kappa^{\mf x}_\lambda :  {\M}^{\mf x}({\lambda}) \longrightarrow {\mathcal{T}}^{\mf x}(\lambda)$ is an ${\mf x}_\infty$-crystal isomorphism.
\end{prop}
\pf It suffices to show that $\kappa^{\mf x}_\lambda$ commutes with $\td{E}_0, \td{F}_0$. 
Recall that $\td{E}_0=\td{e}^\epsilon_0$ and $\td{F}_0=\td{f}_0^\epsilon$, which act on $\M^{\mf x}\subset \M$. Let $A\in \M^{\mf x}$ be given with $\kappa(A)=(S,S^\ast)$ for some $S\in SST_{\Z_{>0}}(\tau)$. It is straightforward to check that $\td{X}_0(S,S^\ast)=\td{x}^\epsilon_0(S,S^\ast)=(S',(S')^\ast)$ ($\td{x}_0(S,S^\ast)$ defined in Section \ref{gl crystal on MAB}) if and only if $\td{X}_0S=S'$ ($\td{X}_0S$ defined in Section \ref{crystal Txlambda}) for some $S'\in SST_{\Z_{>0}}(\tau)$  ($x=e,f$, $X=E,F$).

Then, it follows from Proposition \ref{KnuthIso} that $\kappa^{\mf x}_\lambda$ commutes with $\td{E}_0, \td{F}_0$, and hence an ${\mf x}_\infty$-crystal isomorphism.
\qed\vskip 2mm

\subsection{A grading on $\M^{\mf x}(\lambda)$}
In this subsection, we introduce a $\Z_{\geq 0}$-valued function on $\M^{\mf x}(\lambda)$, which is constant on each connected $\mathfrak{l}_\infty$-subcrystal of $\M^{\mf x}(\lambda)$.

\begin{df}\label{def of L}{\rm
Let ${\bf m}=(A,T)\in  {\M}^{\mf x}(\lambda)$ be given. 
\begin{itemize}
\item[(1)] A word  $\w=\w_1\cdot\w_2$ is called a {\it weakly decreasing sequence for ${\bf m}$} if

\hskip -7mm (a) $\w_1=w(A')$ for some $A'\in{\M}^{\mf x}$ with $A'\leq A$,

\hskip -7mm (b) $\w_2$ is a subword of $w(T)$,

\hskip -7mm (c) $\w_1\cdot\w_2$ is weakly decreasing.

\hskip -7mm Here,  $A'\leq A$ means $a'_{i,j}\leq a_{i,j}$ for all $i,j$, and $w(A')=\bi$ when $A'=A(\bi,\bj)$.

\item[(2) ]A weakly decreasing sequence  $\w_1\cdot\w_2$ for ${\bf m}$ is called  {\it maximal} if $\ell(\w_1)+2\ell(\w_2)$ is maximal, where $\ell(\w)$ denotes the length of a word $\w$.
In this case, we define
\begin{equation*}\label{c(m)}
L_{\mf x}({\bf m})=\frac{\ell(\w_1)+2\ell(\w_2)}{\epsilon}.
\end{equation*}
\end{itemize}
}
\end{df}

\begin{rem}\label{remark on L}{\rm (1) Let $A'\in{\M}^{\mf x}$ be a matrix as in Definition \ref{def of L} (1) (a) with $\w_1=w(A')$.  Let $${\rm supp}(A')=\{\,(i,-j^\vee)\,|\,a'_{i,-j^\vee}\neq 0\,\}\subset \Z_{>0}\times \Z_{<0}^\vee$$ be the support of $A'=(a'_{i,j})$. Then
(i) $(i,-j^\vee)\in {\rm supp}(A')$ if and only if $(j,-i^\vee)\in {\rm supp}(A')$, and (ii)  if $(i,-j^\vee)\in {\rm supp}(A')$, then $(k,-l^\vee)\not \in {\rm supp}(A')$ for $k<i$ and $l<j$ since $\w_1$ is weakly decreasing.
Note that $\w_1$ is a subword of $w(A)$, and $\ell(\w_1)$ is even when $\epsilon=2$ since $A'$  is symmetric and  its diagonal entries are even. Hence $L_{\mf x}({\bf m})$ is an integer.

(2) Let $\w_2=w_1\ldots w_r$ be a subword of $w(T)$ as in Definition \ref{def of L} (1) (b). Suppose that $w_k$ is the entry of $T$ at the $(i_k,j_k)$ position. By the same arguments as in Remark \ref{remark on Delta}, we may assume that $j_1< \ldots < j_r$ and $i_1\geq \ldots \geq i_r$.}
\end{rem}

\begin{lem}\label{normality-0} Let ${\bf m}=(A,T)\in {\M}^{\mf x}({\lambda})$ be given with $A=(a_{i,-j^\vee})$.
\begin{itemize}
\item[(1)] If $\td{E}^k_0{\bf m}\neq {\bf 0}$ for some $k\geq 1$, then $L_{\mf x}(\td{E}^k_0{\bf m})\leq L_{\mf x}({\bf m})$.

\item[(2)]If  $\td{F}^k_0{\bf m}\neq {\bf 0}$ for some $k\geq 1$, then
\begin{equation*}
L_{\mf x}(\td{F}_0^k{\bf m})=\max \left\{\,L_{\mf x}({\bf m})\ ,\ \frac{1}{\epsilon}a_{1,-1^\vee} +\frac{2}{\epsilon}\sum_{i\geq 2}a_{1,-i^\vee}+\frac{2}{\epsilon}m_1(T)+ k \,\right\}.
\end{equation*}
\end{itemize}
\end{lem}
\pf (1) It suffices to show for the case $k=1$.
Suppose that $\td{E}_0{\bf m}\neq {\bf 0}$, that is, $a_{1,-1^\vee}\geq \epsilon$. Let  $\td{\w}_1\cdot\td{\w}_2$ be a maximal  weakly decreasing sequence for $\td{E}_0{\bf m}$ with $\td{\w}_1=w(\td{A})$ for some $\td{A}\in{\M}^{\mf x}$ with $\td{A}\leq \td{E}_0A$.

If $(1,-1^\vee)\not\in {\rm supp}(\td{A})$, then  $\td{\w}=\td{\w}_1\cdot\td{\w}_2$ is a weakly decreasing sequence for ${\bf m}$, which implies that $L_{\mf x}(\td{E}_0{\bf m})\leq L_{\mf x}({\bf m})$. If $(1,-1^\vee)\in {\rm supp}(\td{A})$, then choose $\td{\w }'_1$ of maximal length such that $\td{\w}'_1=w(\td{A}')$ for some $\td{A}'\in{\M}^{\mf x}$ with $\td{A}'\leq A$ and ${\rm supp}(\td{A}')={\rm supp}(\td{A})$. Then $\td{\w}'_1\cdot\td{\w}_2$ is a weakly decreasing sequence for ${\bf m}$ and we have $L_{\mf x}(\td{E}_0{\bf m})+1\leq L_{\mf x}({\bf m})$ since $\ell(\td{\w}'_1)+\epsilon=\ell(\td{\w}_1)$.

(2) Let $\w_1(k)\cdot\w_2$ be a weakly decreasing sequence for $\td{F}_0^k{\bf m}$ ($k\geq 1$), where $\w_1(k)=w(A'(k))$ for some $A'(k)\in{\M}^{\mf x}$ with $${\rm supp}(A'(k))\subset \bigcup_{j\geq 1}\{\,(1,-j^\vee) ,(j,-1^\vee)\,\}$$ and $\w_2=1^{m_1(T)}$.
If $\w_1(k)$ is of maximal length, then $$\ell(k):=\frac{\ell(\w_1(k))+2\ell(\w_2)}{\epsilon}=\frac{1}{\epsilon}a_{1,-1^\vee} +\tfrac{2}{\epsilon}\sum_{i\geq 2}a_{1,-i^\vee}+\tfrac{2}{\epsilon}m_1(T)+k,$$ and $\ell(k)\leq L_{\mf x}(\td{F}_0^k{\bf m})$. Also by (1),  $L_{\mf x}({\bf m})\leq L_{\mf x}(\td{F}_0^k{\bf m})$. Hence we have $$\max\{\,L_{\mf x}({\bf m}),\ell(k)\,\}\leq L_{\mf x}(\td{F}_0^k{\bf m}).$$

Conversely, let $\td{\w}_1\cdot\td{\w}_2$ be a maximal weakly decreasing sequence for $\td{F}_0^k{\bf m}$, where $\td{\w}_1=w(A')$ for some $A'\in{\M}^{\mf x}$. If $(1,-1^\vee)\in {\rm supp}(A')$, then $\td{\w}_1=\w_1(k)$ and $\td{\w}_2=1^{m_1(T)}$, which implies that $L_{\mf x}(\td{F}_0^k{\bf m})=\ell(k)$. If $(1,-1^\vee)\not\in {\rm supp}(A')$, then $\td{\w}_1\cdot\td{\w}_2$ is also a weakly decreasing sequence for ${\bf m}$, and $L_{\mf x}(\td{F}_0^k{\bf m})\leq L_{\mf x}({\bf m})$. Therefore, we have $$\max\{\,L_{\mf x}({\bf m}),\ell(k)\,\}\geq L_{\mf x}(\td{F}_0^k{\bf m}).$$
\qed

\begin{lem}\label{normality-1.1}
If $\td{F}_i{\bf m}\neq {\bf 0}$ for ${\bf m}\in {\M}^{\mf x}({\lambda})$ and $i\in\Z_{>0}$, then $L_{\mf x}(\td{F}_i{\bf m})\leq L_{\mf x}({\bf m})$.
\end{lem}
\pf Let ${\bf m}=(A,T)\in {\M}^{\mf x}({\lambda})$ with $A=(a_{i,-j^\vee})$.

\textsc{Case 1.} Suppose that $\td{F}_i{\bf m}=(\td{F}_iA,T)\neq {\bf 0}$ and
$\td{f}_iA= A-E_{i,-j^\vee}+E_{i+1,-j^\vee}$ for some $j$, equivalently $\td{f}_{-i}A= A-E_{j,-i^\vee}+E_{j,-(i+1)^\vee}$.
Note that
\begin{equation*}
\begin{split}
&\td{F}_iA=\td{f}_i\td{f}_{-i}A= \\
&\begin{cases}
A-E_{i,-i^\vee}+E_{i+1,-(i+1)^\vee}, & \text{if $i=j$ and $a_{i,-i^\vee}=1$}, \\
A-E_{i,-j^\vee}+E_{i+1,-j^\vee}-E_{j,-i^\vee}+E_{j,-(i+1)^\vee}, & \text{otherwise}.
\end{cases}
\end{split}
\end{equation*}
Let $\td{\w}_1\cdot\td{\w}_2$ be a maximal weakly decreasing sequence for $\td{F}_i{\bf m}$, where  $\td{\w}_1=w(\td{A})$ for some $\td{A}\in {\M}^{\mf x}$ with $\td{A}\leq \td{F}_i A$ and $\td{\w}_2=\td{w}_{2,1}\cdots \td{w}_{2,p}$. Let $P=(i,-j^\vee)$ and $Q=(i+1,-j^\vee)$.

(i) Suppose that either $P, Q\not\in {\rm supp}(\td{A})$ or $P,Q\in {\rm supp}(\td{A})$. In this case, $\td{\w}_1\cdot\td{\w}_2$ is also a weakly decreasing sequence for ${\bf m}$. Hence we have $L_{\mf x}(\td{F}_i{\bf m})\leq L_{\mf x}({\bf m})$.

(ii) Suppose that $P\in {\rm supp}(\td{A})$ and $Q\not\in {\rm supp}(\td{A})$.  Let $\td{\w}'_1$ be of maximal length such that $\td{\w}'_1=w (\td{A}')$ for some $\td{A}'\in \M^{\mf x}$ with $\td{A}'\leq A$ and ${\rm supp}(\td{A}')={\rm supp}(\td{A})$.
Then  $\td{\w}'_1\cdot\td{\w}_2$ is a weakly decreasing sequence for ${\bf m}$ and
\begin{equation*}
\ell(\td{\w}'_1)=
\begin{cases}
\ell(\td{\w}_1)+1, & \text{if $i=j$ and $a_{i,-i^\vee}=1$}, \\
\ell(\td{\w}_1)+2, & \text{otherwise},
\end{cases}
\end{equation*}
which  imply that $L_{\mf x}(\td{F}_i{\bf m})\leq L_{\mf x}({\bf m})-1<L_{\mf x}({\bf m})$.

(iii) Suppose that $P\not\in {\rm supp}(\td{A})$ and $Q\in {\rm supp}(\td{A})$. Let ${\rm supp}(\td{A})=\{\,Q_1=(i_1,-j_1^\vee),\ldots,Q_r=(i_r,-j_r^\vee) \,\}$, where $Q_1<\ldots<Q_r$ with respect to the lexicographic orderng on $\Z_{>0}\times\Z_{<0}^\vee$ (see Section \ref{MAB}), and $Q_s=Q$ for some $1\leq s\leq  r$.

If  $i_{s+1}<i_s=i+1$ with $s<r$, then choose $\td{\w }'_1$ of maximal length such that $\td{\w}'_1=w(\td{A}')$ for some $\td{A}'\in {\M}^{\mf x}$  with $\td{A}'\leq A$ and
$${\rm supp}(\td{A}')\subset{\rm supp}(\td{A})\cup \{\,(i,-j^\vee), (j,-i^\vee)\,\}.$$ Then $\td{\w}'_1\cdot\td{\w}_2$ is a weakly decreasing sequence for ${\bf m}$ and $\ell(\td{\w}_1)\leq\ell(\td{\w}'_1)$, which imply that
$L_{\mf x}(\td{F}_i{\bf m})\leq L_{\mf x}({\bf m})$.

If $i_{s+1}=i+1$ with $s<r$, then choose a maxiaml $t\geq s+1$ such that $i_t=i+1$. Consider the subword of $w(A)$ consisting of $i$ and $i+1$. Then by tensor product rule of crystals with respect to $\tf_i$ and $\tf_{-i}$ (see \cite[Proposition 2.1.1 (i)]{KN})  we necessarily have
\begin{equation}\label{inequality}
\begin{split}
&\sum_{j_{s}+1\leq k\leq j_t}a_{i+1,-k^\vee}< \sum_{j_{s}\leq l\leq j_t-1}a_{i,-l^\vee},\\
&\sum_{j_{s}+1\leq k\leq j_t}a_{k,-(i+1)^\vee}< \sum_{j_{s}\leq l\leq j_t-1}a_{l,-i^\vee}.
\end{split}
\end{equation}
If $t<r$ or $\td{w}_{2,1}\leq i$, then choose $\td{\w }'_1$ of maximal length such that $\td{\w}'_1=w(\td{A}')$  for some $\td{A}'\in {\M}^{\mf x}$  with $\td{A}'\leq A$ and
\begin{equation}\label{support}
\begin{split}
&{\rm supp}(\td{A}')\subset   \left({\rm supp}(\td{A})\setminus X\right)\cup Y,
\end{split}
\end{equation}
where
\begin{equation*}
\begin{split}
& X=  \{ \,(i+1,-k^\vee), (k,-(i+1)^\vee)\,|\,j_s+1=j+1\leq k\leq j_t\,\},\\
& Y=  \{\,(i,-l^\vee), (l,-i^\vee)\,|\,j_s=j\leq l\leq j_t-1\,\}.
\end{split}
\end{equation*}
Then
 $\td{\w}'_1\cdot\td{\w}_2$ is a weakly decreasing sequence for ${\bf m}$ by (\ref{support}),  and  $\ell(\td{\w}_1)\leq\ell(\td{\w}'_1)$ by (\ref{inequality}), which imply that
$L_{\mf x}(\td{F}_i{\bf m})\leq L_{\mf x}({\bf m})$. If $t=r$ and $\td{w}_{2,1}=i+1$, then choose a maximal $u$ such that $\td{w}_{2,u}=i+1$, where $\td{w}_{2,u}=w_v$ for some $1\leq v\leq p$. Also by
considering the subword of $w(A)\cdot w(T)$ consisting of $i$ and $i+1$ and then applying \cite[Proposition 2.1.1 (i)]{KN}, we necessarily have

\begin{equation}\label{inequality-2'}
\begin{split}
&\sum_{j_{s}+1\leq k}a_{i+1,-k^\vee}+m_{i+1}< \sum_{j_{s}\leq l}a_{i,-l^\vee}+m_{i},\\
\end{split}
\end{equation}
where $m_i$ is the number of occurrences of $i$ in $w_1\ldots w_{v-1}$ and $m_{i+1}$ is the number of occurrences of $i+1$ in $w_1\ldots w_{v}$. Note that $m_{i+1}=u$ by the maximality of $\td{\w}_1\cdot\td{\w}_2$. Choose $\td{\w }'_1$ of maximal length such that $\td{\w}'_1=w(\td{A}')$ for some $\td{A}'\in {\M}^{\mf x}$ with $\td{A}'\leq A$ and
\begin{equation}\label{support-2}
\begin{split}
&{\rm supp}(\td{A}')\subset   \left({\rm supp}(\td{A})\setminus X\right)\cup Y,
\end{split}
\end{equation}
where
\begin{equation*}
\begin{split}
& X=  \{ \,(i+1,-k^\vee), (k,-(i+1)^\vee)\,|\,j_s+1=j+1\leq k\,\},\\
& Y=  \{\,(i,-l^\vee), (l,-i^\vee)\,|\,j_s=j\leq l\,\}.
\end{split}
\end{equation*}
Let $\td{\w}'_2$ be the subword obtained from $\td{\w}_2$ by replacing $\td{w}_{2,1}\ldots \td{w}_{2,u}=(i+1)^{m_{i+1}}$ with $i^{m_i}$. Then $\td{\w}'_1\cdot\td{\w}'_2$ is a weakly decreasing sequence for ${\bf m}$ by (\ref{support-2}), and
\begin{equation*}\label{length difference}
\begin{split}
&\left(\ell(\td{\w}'_1)+2\ell(\td{\w}'_2)\right)-\left(\ell(\td{\w}_1)+2\ell(\td{\w}_2)\right) \\
&=\sum_{j_{s}\leq l}(a_{i,-l^\vee}+a_{l,-i^\vee}) +2m_i-2-\sum_{j_{s}+1\leq k}(a_{i+1,-k^\vee}+a_{k,-(i+1)^\vee})-2m_{i+1}\\
&=2\sum_{j_{s}\leq l}a_{i,-l^\vee}+2m_i-2\sum_{j_{s}+1\leq k}a_{i+1,-k^\vee}-2m_{i+1}-2\geq 0.
\end{split}
\end{equation*}
by (\ref{inequality-2'}). Hence
$L_{\mf x}(\td{F}_i{\bf m})\leq L_{\mf x}({\bf m})$.
\vskip 2mm

\textsc{Case 2.} Suppose that $\td{F}_i{\bf m}=(A,\td{F}_iT)\neq {\bf 0}$. Let $\td{\w}_1\cdot\td{\w}_2$ be a maximal weakly decreasing sequence for $\td{F}_i{\bf m}$. Let $w(T)=w_{1}\ldots w_{p}$ and $\td{\w}_2=\td{w}_{2,1}\ldots \td{w}_{2,q}$. Assume that $w_r=i$ for some $1\leq r\leq p$ in $T$ is replaced by $w'_r=i+1$ by applying $\td{F}_i$ to $T$.
If $\td{\w}_2$ does not contain $i+1$ corresponding to $w'_r$, then $\td{\w}_1\cdot\td{\w}_2$ is also a weakly decreasing sequence for ${\bf m}$ and hence $L_{\mf x}(\td{F}_i{\bf m})\leq L_{\mf x}({\bf m})$. Now we assume that $\td{\w}_2$  contains $i+1$  corresponding to $w'_r$, say $\td{w}_{2,s}=i+1$.

(i) Suppose that $s=q$ or $\td{w}_{2,s+1}<i+1$. Let $\td{\w}'_2$ be obtained from $\td{\w}_2$ by replacing $\td{w}_{2,s}$ with $i$. Then $\td{\w}_1\cdot\td{\w}'_2$ is a weakly decreasing sequence for ${\bf m}$ and hence $L_{\mf x}(\td{F}_i{\bf m})\leq L_{\mf x}({\bf m})$.

(ii) Suppose that $\td{w}_{2,s+1}=i+1$. Choose a maximal $t>s$ such that $\td{w}_{2,t}=i+1$. If $\td{w}_{2,t}=w_{u}$ for some $u>r$, then by  \cite[Proposition 2.1.1 (i)]{KN} we have $u>r+1$ and there exists a subword $i^{t-s}$ of $w_{r+1}\ldots w_{u-1}$. Let $\td{\w}'_2$ be obtained from $\td{\w}_2$ by replacing $\td{w}_{2,s}\ldots \td{w}_{2,t}=(i+1)^{t-s+1}$ with $i^{t-s+1}$. Then $\td{\w}_1\cdot\td{\w}'_2$ is a weakly decreasing sequence for ${\bf m}$, and hence $L_{\mf x}(\td{F}_i{\bf m})\leq L_{\mf x}({\bf m})$.\qed\vskip 3mm

\begin{lem}\label{normality-1.2}
If $\td{E}_i{\bf m}\neq {\bf 0}$ for ${\bf m}\in {\M}^{\mf x}({\lambda})$ and $i\in\Z_{>0}$, then $L_{\mf x}(\td{E}_i{\bf m})\leq L_{\mf x}({\bf m})$.
\end{lem}
\pf The proof is similar to that of Lemma \ref{normality-1.1}. Let ${\bf m}=(A,T)\in {\M}^{\mf x}({\lambda})$ with $A=(a_{i,-j^\vee})$.

\textsc{Case 1.} Suppose that $\td{E}_i{\bf m}=(\td{E}_iA,T)\neq {\bf 0}$ and $\te_iA= A+E_{i,-j^\vee}-E_{i+1,-j^\vee}$ for some $j$, equivalently $\te_{-i}A=E_{j,-i^\vee}-E_{j,-(i+1)^\vee}$. Let $P=(i,-j^\vee)$ and $Q=(i+1,-j^\vee)$.

Let $\td{\w}_1\cdot\td{\w}_2$ be a maximal weakly decreasng sequence for $\td{E}_i{\bf m}$, where $\td{\w}_1=w(\td{A})$ for some $\td{A}\in {\M}^{\mf x}$ with $\td{A}\leq \td{E}_i A$.

(i) Suppose that either $P, Q\not\in {\rm supp}(\td{A})$ or $P,Q\in {\rm supp}(\td{A})$. In this case, $\td{\w}_1\cdot\td{\w}_2$ is also a weakly decreasing sequence for ${\bf m}$, and hence $L_{\mf x}(\td{E}_i{\bf m})\leq L_{\mf x}({\bf m})$.

(ii) Suppose that $P\not\in {\rm supp}(\td{A})$ and $Q\in {\rm supp}(\td{A})$.  Let $\td{\w}'_1$ be of maximal length such that $\td{\w}'_1=w (\td{A}')$ with $\td{A}'\leq A$ and ${\rm supp}(\td{A}')={\rm supp}(\td{A})$.
Then  $\td{\w}'_1\cdot\td{\w}_2$ is a weakly decreasing sequence for ${\bf m}$ and
\begin{equation*}
\ell(\td{\w}'_1)=
\begin{cases}
\ell(\td{\w}_1)+1, & \text{if $i=j$ and $a_{i+1,-(i+1)^\vee}=1$}, \\
\ell(\td{\w}_1)+2, & \text{otherwise},
\end{cases}
\end{equation*}
which  imply that $L_{\mf x}(\td{E}_i{\bf m})\leq L_{\mf x}({\bf m})-1<L_{\mf x}({\bf m})$.

(iii) Suppose that $P\in {\rm supp}(\td{A})$ and $Q\not\in {\rm supp}(\td{A})$. Let ${\rm supp}(\td{A})=\{\,P_1=(i_1,-j_1^\vee),\ldots,P_r=(i_r,-j_r^\vee) \,\}$, where $P_1<\ldots<P_r$ with respect to the lexicographic orderng on $\Z_{>0}\times\Z_{<0}^\vee$, and $P_s=P$ for some $1\leq s\leq  r$.

If $i_{s-1}>i_s=i$ or $s=1$, choose $\td{\w }'_1$ of maximal length such that $\td{\w}'_1=w(\td{A}')$ for some $\td{A}'\in\M^{\mf x}$ with $\td{A}'\leq A$ and $${\rm supp}(\td{A}')\subset{\rm supp}(\td{A})\cup \{\,(i+1,-j^\vee), (j,-(i+1)^\vee)\,\}.$$ Then $\td{\w}'_1\cdot\td{\w}_2$ is a weakly decreasing sequence for ${\bf m}$ and $\ell(\td{\w}_1)\leq\ell(\td{\w}'_1)$, which imply that
$L_{\mf x}(\td{E}_i{\bf m})\leq L_{\mf x}({\bf m})$.

If $s\geq 2$ and $i_{s-1}=i$, then choose a minimal $t\leq s-1$ such that $i_t=i$.
Considering the subword of $w(A)\cdot w(T)$ consisting of $i$ and $i+1$ and then applying \cite[Proposition 2.1.1 (i)]{KN}, we have
\begin{equation}\label{inequality-3}
\begin{split}
&\sum_{j_{t}\leq k\leq j_s-1}a_{i,-k^\vee}< \sum_{j_{t}+1\leq l\leq j_s}a_{i+1,-l^\vee}, \\
&\sum_{j_{t}\leq k\leq j_s-1}a_{k,-i^\vee}< \sum_{j_{t}+1\leq l\leq j_s}a_{l,-(i+1)^\vee}.
\end{split}
\end{equation}
Choose $\td{\w }'_1$ of maximal length such that $\td{\w}'_1=w(\td{A}')$ for some $\td{A}'\in\M^{\mf x}$ with $\td{A}'\leq A$ and
\begin{equation}\label{support-3}
\begin{split}
&{\rm supp}(\td{A}')\subset   \left({\rm supp}(\td{A})\setminus X\right)\cup Y,
\end{split}
\end{equation}
where
\begin{equation*}
\begin{split}
& X=  \{ \,(i,-k^\vee), (k,-i^\vee)\,|\,j_t\leq k\leq j_s-1=j-1\,\},\\
& Y=  \{\,(i+1,-l^\vee), (l,-(i+1)^\vee)\,|\,j_t+1\leq l\leq j_s=j\,\}.
\end{split}
\end{equation*}
Then $\td{\w}'_1\cdot\td{\w}_2$ is a weakly decresing sequence for ${\bf m}$ by (\ref{support-3}),  and  $\ell(\td{\w}_1)\leq \ell(\td{\w}'_1)$ by (\ref{inequality-3}), which also imply that
$L_{\mf x}(\td{E}_i{\bf m})\leq L_{\mf x}({\bf m})$. \vskip 2mm

\textsc{Case 2.} Suppose that $\td{E}_i{\bf m}=(A,\td{E}_iT)\neq {\bf 0}$. Let $\td{\w}_1\cdot\td{\w}_2$ be a maximal weakly decreasng sequence for $\td{E}_i{\bf m}$  with $\td{\w}_1=\td{w}_{1,1}\ldots \td{w}_{1,r}$ and $\td{\w}_2=\td{w}_{2,1}\ldots \td{w}_{2,s}$.  Assume that $w(T)=w_{1}\ldots w_{u}$ and $w_v=i+1$ ($1\leq v\leq u$) in $T$ is replaced by $w'_v=i$ by applying  $\td{E}_i$ to $T$.
If $\td{\w}_2$ does not contain $i$ corresponding to $w'_v$, then $\td{\w}_1\cdot\td{\w}_2$ is also a weakly decreasing sequence for ${\bf m}$ and hence $L_{\mf x}(\td{E}_i{\bf m})\leq L_{\mf x}({\bf m})$. Now we assume that $\td{\w}_2$  contains $i$  corresponding to $w'_v$, say $\td{w}_{2,p}=i$.

(i) Suppose that $\td{w}_{1,r}\geq i+1$. If  $p=1$ or $\td{w}_{2,p-1}>i$, then $\td{\w}_1\cdot\td{\w}'_2$ is a weakly decreasing sequence for ${\bf m}$, where $\td{\w}'_2$ is obtained from $\td{\w}_2$ by replacing $\td{w}_{2,p}$ with $i+1$, and hence $L_{\mf x}(\td{E}_i{\bf m})\leq L_{\mf x}({\bf m})$. Otherwise, we have $p\geq 2$ and $\td{w}_{2,p-1}=i$. In this case, choose a minimal $1\leq q<p$ such that $\td{w}_{2,q}=i$. If $\td{w}_{2,q}=w_t$ for some $1\leq t<v$, then by \cite[Proposition 2.1.1 (ii)]{KN} there exists a subword $(i+1)^{p-q}$ of $w_{t+1}\ldots w_{v-1}$. Let $\td{\w}'_2$ be obtained from $\td{\w}_2$ by replacing $\td{w}_{2,q}\ldots \td{w}_{2,p}=i^{p-q+1}$ with $(i+1)^{p-q+1}$. Then $\td{\w}_1\cdot\td{\w}'_2$ is a weakly decreasing sequence for ${\bf m}$. Hence   $L_{\mf x}(\td{E}_i{\bf m})\leq L_{\mf x}({\bf m})$.

(ii) Suppose that $\td{w}_{1,r}=i$. We have $\td{w}_{2,1}=\cdots=\td{w}_{2,p}=i$. Note that $p$ is the number of occurrences of $i$ in $w_{1}\ldots w_{v-1}$ by the maximality of $\td{\w}_1\cdot\td{\w}_2$.  Next, let ${\rm supp}({A})=\{\,P_1=(i_1,-j_1^\vee),\ldots,P_a=(i_a,-j_a^\vee) \,\}$ with $P_1<\ldots<P_a$ and let $b$ be the smallest such that $i_b=i$. Note that $i_k=i$ for $b\leq k\leq a$. Again by \cite[Proposition 2.1.1 (ii)]{KN}, we have
\begin{equation}\label{inequality-4}
\sum_{j_b\leq j} a_{i,-j^\vee}+m_{i}< \sum_{j_b\leq j} a_{i+1,-(j+1)^\vee}+m_{i+1},
\end{equation}
where $m_i$ is the number of occurrences of $i$ in $w_1\ldots w_{v-1}$ and  $m_{i+1}$ is the number of occurrences of $i+1$ in $w_1\ldots w_v$.
Choose $\td{\w }'_1$ of maximal length such that $\td{\w}'_1=w({A}')$ for some $A'\in \M^{\mf x}$ with ${A}'\leq A$ and
\begin{equation}\label{support-4}
\begin{split}
&{\rm supp}({A}')\subset   \left({\rm supp}({A})\setminus X\right)\cup Y,
\end{split}
\end{equation}
where
\begin{equation*}
\begin{split}
& X=  \{ \,(i,-k^\vee), (k,-i^\vee)\,|\,j_b\leq k\,\},\\
& Y=  \{\,(i+1,-l^\vee), (l,-(i+1)^\vee)\,|\,j_b+1\leq l\,\}.
\end{split}
\end{equation*}
Let $\td{\w}'_2$ be obtained from $\td{\w}_2$ by replacing $\td{w}_{2,1}\ldots\td{w}_{2,p}=i^{m_i+1}$ with $(i+1)^{m_{i+1}}$. Then $\td{\w}'_1\cdot\td{\w}'_2$ is a weakly decreasing sequence for ${\bf m}$  by (\ref{support-4}),  and
\begin{equation*}\label{length difference-2}
\begin{split}
&\left(\ell(\td{\w}'_1)+2\ell(\td{\w}'_2)\right)-\left(\ell(\td{\w}_1)+2\ell(\td{\w}_2)\right) \\
&=\sum_{j_{b}+1\leq l}(a_{i+1,-l^\vee}+a_{l,-(i+1)^\vee})+2m_{i+1}-\sum_{j_{b}\leq k}(a_{i,-k^\vee}+a_{k,-i^\vee}) -2(m_i+1)\\
&=2\sum_{j_{b}+1\leq l}a_{i+1,-l^\vee}+2m_{i+1}-2\sum_{j_{b}\leq k}a_{i,-k^\vee}  -2m_i-2\geq 0,
\end{split}
\end{equation*}
by (\ref{inequality-4}). Hence, $L_{\mf x}(\td{E}_i{\bf m})\leq L_{\mf x}({\bf m})$.
\qed\vskip 2mm

\begin{prop}\label{normality}
$L_{\mf x}(\cdot)$ is constant on each connected $\mathfrak{l}_{\infty}$-subcrystal of $\M^{\mf x}(\lambda)$.
\end{prop}
\pf Suppose that ${\bf m}'=\td{F}_i{\bf m}\neq {\bf 0}$ for ${\bf m}, {\bf m}'\in \M^{\mf x}(\lambda)$ and $i\in\Z_{>0}$. By Lemmas \ref{normality-1.1} and \ref{normality-1.2}, it follows that
$$L_{\mf x}({\bf m})= L_{\mf x}(\td{E}_i{\bf m}')\leq L_{\mf x}({\bf m}')=L_{\mf x}(\td{F}_i{\bf m})\leq L_{\mf x}({\bf m}).$$  Hence we have $L_{\mf x}({\bf m})= L_{\mf x}({\bf m}')$.
\qed

\subsection{Restriction of $ {\M}^{\mf x}(\lambda)$ to a normal highest weight crystal}
\begin{df}{\rm
For  $(\lambda,n)\in\cP({\mf x})$, define
\begin{equation*}\label{Mlambda,n}
\begin{split}
\M^{\mf x}(\lambda,n)& = \{\,{\bf m}\otimes t_{n{\Lambda}^{\mf x}_{0}}\in {\M}^{\mf x}({\lambda})\otimes T_{n{\Lambda}^{\mf x}_{0}} \,|\,L_{\mf x}({\bf m})\leq\,n\,\}.
\end{split}
\end{equation*}}
\end{df}

Regarding $T_{n{\Lambda}^{\mf x}_{0}}=\{\,t_{n{\Lambda}^{\mf x}_{0}}\,\}$ as an ${\mf x}_\infty$-crystal with ${\rm wt}^{\mf x}(t_{n{\Lambda}^{\mf x}_{0}})=n{\Lambda}^{\mf x}_{0}$, ${\varepsilon}^{\mf x}_i(t_{n{\Lambda}^{\mf x}_{0}})=\varphi_i^{\mf x}(t_{n{\Lambda}^{\mf x}_{0}})=-\infty$ ($i\in\Z_{\geq 0}$), we may view $\M^{\mf x}(\lambda,n)$ as an ${\mf x}_\infty$-subcrystal of ${\M}^{\mf x}({\lambda})\otimes T_{n{\Lambda}^{\mf x}_{0}}$.

\begin{lem}\label{normality-2} $\M^{\mf x}(\lambda,n)$ is  a connected normal ${\mf x}_\infty$-crystal  with a unique highest weight element ${\mathbb{O}}^{\mf x}_{(\lambda,n)}=(\mathbb{O}, H_\lambda)\otimes t_{n{\Lambda}^{\mf x}_{0}} \in \M^{\mf x}(\lambda,n)$ of weight ${\Lambda}^{\mf x}(\lambda,n)$.
\end{lem}
\pf Let us first show that  $\M^{\mf x}(\lambda,n)$ is a normal ${\mf x}_\infty$-crystal.

Let ${\bf b}={\bf m}\otimes t_{n{\Lambda}^{\mf x}_{0}} \in \M^{\mf x}(\lambda,n)$
be given, where ${\bf m}=(A,T)$ with $A=(a_{i,-j^\vee})$.  Let $S_i({\bf b})=\{\,\td{X}^k_i{\bf b}\,|\,k\geq 0, X=E,F\,\}\setminus\{{\bf 0}\}$ for $i\in\Z_{\geq 0}$, where we regard ${\bf b}$ as an element in the crystal ${\M}^{\mf x}({\lambda})\otimes T_{n{\Lambda}^{\mf x}_{0}}$. Note that $\td{X}_i^k{\bf b}\neq {\bf 0}$ if and only if $\td{X}_i^k{\bf m}\neq {\bf 0}$ ($X=E, F$) since $\td{X}_i^k{\bf b}=\left( \td{X}_i^k{\bf m} \right)\otimes  t_{n{\Lambda}^{\mf x}_{0}}$.

Suppose that $i\neq 0$.  Since  ${\M}^{\mf x}({\lambda})$ (and hence ${\M}^{\mf x}({\lambda})\otimes T_{n{\Lambda}^{\mf x}_{0}}$) is a normal $\mathfrak{l}_{\infty}$-crystal by (\ref{l infty decomposition}) and ${\varepsilon}^{{\mf x}}_i({\bf b})={\varepsilon}^{{\mf x}}_i({\bf m})$,
${\varphi}^{{\mf x}}_i({\bf b})={\varphi}^{{\mf x}}_i({\bf m})$, we have
$$S_i({\bf b})=\{\,\td{E}_i^k{\bf b}, \td{F}_i^l{\bf b}\,|\,0\leq k\leq {\varepsilon}^{\mf x}_i({\bf b}), 0\leq l\leq {\varphi}^{{\mf x}}_i({\bf b})\,\}.$$
By Proposition \ref{normality}, we have $S_i({\bf b})\subset\M^{\mf x}(\lambda,n)$. Hence, $\M^{\mf x}(\lambda,n)$ is a normal $\mathfrak{l}_{\infty}$-crystal.

Suppose that $i=0$. By Lemma \ref{normality-0}, $\td{E}_0^k{\bf b}\in \M^{\mf x}(\lambda,n)\cup\{{\bf 0}\}$ for $k\geq 0$ and $\varepsilon_0^{\mf x}({\bf b})=\max\{\,k\,|\,\td{E}_0^k{\bf b}\neq {\bf 0}\,\}=\varepsilon^{\mf x}_0({\bf m})$. Next, consider $\td{F}_0^l{\bf b}$ for $l\geq 0$. Since
\begin{equation*}
{\rm wt}^{{\mf x}}({\bf b})=  n{\Lambda}^{\mf x}_0 +\sum_{i\geq 1}\left(m_i({\bi})+m_i(T)\right)\widehat{\epsilon}_i,
\end{equation*}
where $\bi=w(A)$ and  $m_i(\bi)$ (resp. $m_i(T)$) is the number of occurrences of $i$ in $\bi$ (resp. $T$), we have
\begin{equation*}
\langle {\rm wt}^{{\mf x}}({\bf b}), \widehat{h}_0 \rangle =n - \frac{2}{\epsilon}\left(m_1(\bi)+m_1(T)\right)=n -\frac{2}{\epsilon}\left(\sum_{i\geq 1}a_{1,-i^\vee}+m_1(T)\right)
\end{equation*}
and
\begin{equation*}
\begin{split}
{\varphi}^{{\mf x}}_0({\bf b})
&=\langle  {\rm wt}^{{\mf x}}({\bf b}), \widehat{h}_0 \rangle + {\varepsilon}^{{\mf x}}_0({\bf b})
=\langle  {\rm wt}^{{\mf x}}({\bf b}), \widehat{h}_0 \rangle + \frac{1}{\epsilon}a_{1,-1^\vee} \\
&=n -\frac{1}{\epsilon}a_{1,-1^\vee} -\frac{2}{\epsilon}\sum_{i\geq 2}a_{1,-i^\vee} -\frac{2}{\epsilon} m_1(T).
\end{split}
\end{equation*}
By Lemma \ref{normality-0}, we have
\begin{equation*}
L_{\mf x}(\td{F}_0^l{\bf m})=\max \left\{\,L_{\mf x}({\bf m})\ ,\ \frac{1}{\epsilon}a_{1,-1^\vee} +\frac{2}{\epsilon}\sum_{i\geq 2}a_{1-i^\vee}+\frac{2}{\epsilon}m_1(T)+ l \,\right\}.
\end{equation*}
Then
\begin{equation*}
\begin{split}
\td{F}_0^l{\bf b}\in \M^{\mf x}(\lambda,n) & \  \Longleftrightarrow L_{\mf x}(\td{F}^l_0{\bf m})\leq n \\
& \  \Longleftrightarrow
\frac{1}{\epsilon}a_{1,-1^\vee} +\frac{2}{\epsilon}\sum_{i\geq 2}a_{1-i^\vee}+\frac{2}{\epsilon}m_1(T)+ l \leq n \\
& \  \Longleftrightarrow l\leq {\varphi}^{{\mf x}}_0({\bf b}).
\end{split}
\end{equation*}
Therefore, $$S_0({\bf b})\cap \M^{\mf x}(\lambda,n)=\{\,\td{E}_0^k{\bf b}, \td{F}_0^l{\bf b}\,|\,0\leq k\leq {\varepsilon}^{\mf x}_0({\bf b}), 0\leq l\leq {\varphi}^{{\mf x}}_0({\bf b})\,\}.$$
Equivalently,
\begin{equation*}
\begin{split}
{\varepsilon}^{\mf x}_0({\bf b})&=\max\{\,k\,|\,\td{E}_0^k{\bf b}\neq {\bf 0} \text{ in $\M^{\mf x}(\lambda,n)$}\,\},\\
{\varphi}^{\mf x}_0({\bf b})&=\max\{\,k\,|\,\td{F}_0^k{\bf b}\neq {\bf 0} \text{ in $\M^{\mf x}(\lambda,n)$}\,\}.
\end{split}
\end{equation*}
Hence $\M^{\mf x}(\lambda,n)$ is a normal ${\mf x}_\infty$-crystal.

Next, we claim that $\M^{\mf x}(\lambda,n)$ is  connected. Let ${\bf b}={\bf m}\otimes t_{n{\Lambda}^{\mf x}_{0}} \in \M^{\mf x}(\lambda,n)$
be given, where ${\bf m}=(A,T)$ with $A=(a_{i,-j^\vee})$.  Since $\td{E}_i{\bf b}\in \M^{\mf x}(\lambda,n)\cup\{{\bf 0}\}$, we may assume that $\td{E}_i{\bf b}={\bf 0}$ for all $i\in\Z_{\geq 0}$.
Then $\td{E}_iA={\bf 0}$ for $i\in\Z_{\geq 0}$ by tensor product rule, and in particular $A$ is a diagonal matrix with $a_{i,-i^\vee}\geq a_{i+1,-(i+1)^\vee}$ for $i\geq 1$. If $a_{1,-1^\vee}\neq 0$, then we have $\td{E}_0{\bf m}\neq {\bf 0}$, which is a contrdiction. Hence we have $A=\mathbb{O}$, and $T=H_\lambda$, that is, ${\bf b}={\mathbb{O}}^{\mf x}_{(\lambda,n)}$. This implies that $\M^{\mf x}(\lambda,n)=\{\,\td{F}_{i_1}\cdots\td{F}_{i_r}{\mathbb{O}}^{\mf x}_{(\lambda,n)}\,|\,r\geq 0, i_1,\ \ldots,i_r\in\Z_{\geq 0}\,\}$, and it is connected.
\qed

\begin{thm}\label{main-1} For  $(\lambda,n)\in\cP({\mf x})$, we have
$$\M^{\mf x}(\lambda,n)\simeq\B({\mf x}_\infty,{\Lambda}^{\mf x}(\lambda,n)).$$
\end{thm}
\pf Consider the following diagram
$$
\begin{CD}
\M^{\mf x}(\lambda,n)  @  . \ \ \ \ \ \ \ \ \ \ \ \ \ \   \B(\frak{gl}_{\infty},\Lambda(\lambda,\lambda,\epsilon n)) \\
 \cap@.  \ \ \ \ \ \ \ \ \ \ \  \ \  @ VV\Psi_{\lambda,\lambda,n}V \\
 {\M}^{\mf x}({\lambda})\otimes T_{n{\Lambda}^{\mf x}_0} @.\ \ \  \subset \ \ \ \ \ \  \
 \M(\lambda,\lambda)\otimes T_{\epsilon n{\Lambda}_0}
 \end{CD}
$$
where $\Psi_{\lambda,\lambda,n}$ is the embedding in Proposition \ref{Psilambdamu}.
First, we claim that
\begin{equation*}\label{iota-n}
\M^{\mf x}(\lambda,n) \subset {\rm Im}\Psi_{\lambda,\lambda,n}
\end{equation*}
which implies that there exists an injective map
\begin{equation*}
\psi_{\lambda,n} :  \M^{\mf x}(\lambda,n)\longrightarrow \B(\frak{gl}_{\infty},\Lambda(\lambda,\lambda,\epsilon n))
\end{equation*}
such that $\Psi_{\lambda,\lambda,n}\circ \psi_{\lambda,n}({\bf b})={\bf b}$ for ${\bf b}\in \M^{\mf x}(\lambda,n)$.

Let ${\bf b}=\td{F}_{i_r}\ldots\td{F}_{i_1}{\mathbb{O}}^{\mf x}_{(\lambda,n)}\in \M^{\mf x}(\lambda,n)$ be given for $r\geq 0$ and $i_1,\ldots, i_r\in\Z_{\geq 0}$.
We use induction on $r$.
If $r=0$, then
${\mathbb{O}}^{\mf x}_{(\lambda,n)}=\mathbb{O}_{\lambda,\lambda}\otimes t_{n{\Lambda}^{\mf x}_0}=\Psi_{\lambda,\lambda,n}(u_{\Lambda(\lambda,\lambda,n)})$.
Suppose that $r\geq 1$. Put ${\bf b}'=\td{F}_{i_{r-1}}\ldots\td{F}_{i_1}{\mathbb{O}}^{\mf x}_{(\lambda,n)}$.
By induction hypothesis, ${\bf b}'=\Psi_{\lambda,\lambda,n}(b)$ for some $b\in  \B(\frak{gl}_{\infty},\Lambda(\lambda,\lambda,\epsilon n))$.
Since $\Psi_{\lambda,\lambda,n}$ preserves the weights, and commutes with $\td{E}_i$ for $i\in\Z_{\geq 0}$,
we have \begin{equation*}
\begin{split}
&{\varphi}^{\mf x}_{i_r}\left({\bf b}'\right)={\varphi}^{\mf x}_{i_r}\left(\Psi_{\lambda,\lambda,n}(b)\right)  = 
\begin{cases}
{\varphi}_{i_r}(b) ={\varphi}_{-i_r} (b), & \text{if $i_r>0$}, \\
\tfrac{1}{\epsilon}{\varphi}_{0}(b), & \text{if $i_r=0$}.
\end{cases}
\end{split}
\end{equation*}
We have ${\varphi}^{\mf x}_{i_r}\left({\bf b}'\right)\geq 1$ by Lemma \ref{normality-2}, and  $\td{F}_{i_r}b\neq {\bf 0}$ since ${\varphi}_{i_r}(b) ={\varphi}_{-i_r} (b) \geq 1$ and $\B(\frak{gl}_{\infty},\Lambda(\lambda,\lambda,\epsilon n))$ is normal. Hence
\begin{equation*}
{\bf b}=\td{F}_{i_r}\Psi_{\lambda,\lambda,n}(b)=\Psi_{\lambda,\lambda,n}(\td{F}_{i_r}b)\in  {\rm Im}\Psi_{\lambda,\lambda,n}.
\end{equation*}
This completes the induction, and proves the claim.

Recall that $\td{E}_i$ and $\td{F}_i$ ($i\in\Z_{\geq 0}$) on $\B(\frak{gl}_{\infty},\Lambda(\lambda,\lambda,\epsilon n))$ are defined as in  (\ref{folded operators}).
Then we can check that $\psi_{\lambda,n}$ commutes with $\td{E}_i$ and
\begin{equation*}
\begin{split}
{\rm wt}\left(\psi_{\lambda,n}({\bf b})\right)&={\rm wt}^{{\mf x}}({\bf b}), \\
{\varepsilon}_i(\psi_{\lambda,n}({\bf b}))&={\varepsilon}_{-i}(\psi_{\lambda,n}({\bf b}))=
\varepsilon^{{\mf x}}_i({\bf b}) \ \ \ (i>0),\\
{\varphi}_i(\psi_{\lambda,n}({\bf b}))&={\varphi}_{-i}(\psi_{\lambda,n}({\bf b}))=\varphi^{{\mf x}}_i({\bf b})\ \ \ (i>0),\\
{\varepsilon}_0(\psi_{\lambda,n}({\bf b}))&=\epsilon \varepsilon^{{\mf x}}_0({\bf b}), \ \ {\varphi}_0(\psi_{\lambda,n}({\bf b}))=\epsilon \varphi^{{\mf x}}_0({\bf b}),
\end{split}
\end{equation*}
for ${\bf b}\in  \M^{\mf x}(\lambda,n)$ and $i\in\Z_{\geq 0}$.
This implies that $\psi_{\lambda,n}$  commutes with $\td{F}_i$ ($i\in\Z_{\geq 0}$) since $\M^{\mf x}(\lambda,n)$ and $\B(\frak{gl}_{\infty},\Lambda(\lambda,\lambda,\epsilon n))$ are normal.

Therefore,  as an ${\mf x}_\infty$-crystal we have $$\M^{\mf x}(\lambda,n)\simeq {C}(u_{\Lambda(\lambda,\lambda,\epsilon n)}),$$where ${C}(u_{\Lambda(\lambda,\lambda,\epsilon n)})$ is the connected component in $\B(\frak{gl}_{\infty},\Lambda(\lambda,\lambda,\epsilon n))$ generated by $u_{\Lambda(\lambda,\lambda,\epsilon n)}$ with respect to $\td{E}_i$ and $\td{F}_i$ for $i\in\Z_{\geq 0}$. On the other hand, by \cite[Theorem 5.1]{Kas96}, $${C}(u_{\Lambda(\lambda,\lambda,\epsilon n)})\simeq \B({\mf x}_\infty,\Lambda^{\mf x}(\lambda,n))$$ as an ${\mf x}_\infty$-crystal.
This completes the proof. \qed

\subsection{Proof of Theorem \ref{main result}}
Recall that we have an ${\mf x}_\infty$-crystal isomorphism
$$\kappa^{\mf x}_\lambda :  {\M}^{\mf x}({\lambda}) \longrightarrow {\mathcal{T}}^{\mf x}(\lambda)$$
by Proposition \ref{kappa x}.
The following lemma shows that $L_{\mf x}$ and $\Delta_{\mf x}$ coincide under $\kappa^{\mf x}_\lambda$.

\begin{lem}\label{L=Delta}
For ${\bf m}=(A,T)\in\M^{\mf x}(\lambda)$ with $\kappa^{\mf x}_\lambda({\bf m})=(S,T)$, we have
$$L_{\mf x}({\bf m})=\Delta_{\mf x}((T\rightarrow S)_R).$$
\end{lem}
\pf Let ${\bf m}'\in \M^{\mf x}(\lambda)$ be the  highest weight element  of the connected $\mathfrak{l}_{\infty}$-subcrystal of $\M^{\mf x}(\lambda)$ including ${\bf m}$. By Proposition \ref{normality}, $L_{\mf x}({\bf m})=L_{\mf x}({\bf m}')$. On the other hand, the recording tableau $(T\rightarrow S)_R$ is invariant under the actions of $\td{E}_i$ and $\td{F}_i$ on $(S,T)$ for $i\in \Z_{>0}$, that is $(T\rightarrow S)_R=(T'\rightarrow S')_R$, where $\kappa^{\mf x}_\lambda({\bf m}')=(S',T')$ (see for example \cite[Proposition 4.17]{KK}). So, it suffices to show for the case when ${\bf m}$ is an $\mathfrak{l}_{\infty}$-highest weight element.

Since ${\bf m}=(A,T)$ is an $\mathfrak{l}_{\infty}$-highest weight element in ${\M}^{\mf x}({\lambda})$, $S=H_\tau$ for some $\tau\in \cP_{\mf x}$ and $A$ is a diagonal matrix with $a_{i,-i^\vee}=\tau_i$. Also $\left(T \rightarrow H_{\tau} \right)=H_{\sigma}$ for some $\sigma\in\cP$. This implies that $T\in {\rm LR}^{\sigma}_{\tau \lambda}$ or $U=\left(T \rightarrow H_{ \tau} \right)_R\in \texttt{\bf LR}^{\sigma}_{\tau \lambda}$.

For a subword $\w$ of $w(U)$, let us denote by $\w^\sharp$  the subword of $w(T)$ corresponding to $\w$ under the insertion $T\rightarrow S$.
Suppose that $\w=w_1\ldots w_r$ is weakly decreasing and $w_1$ is in the $k$-th row of $\tau$ such that $(\tau_k+2\ell(\w))/\epsilon=\Delta_{\mf x}(U)$.
Then by Remark \ref{remark on Delta}, we may assume that $\w^\sharp=w^\sharp_1\ldots w^\sharp_r$ is a weakly decreasing subword of $w(T)$
with $w^\sharp_1=k$, and $k^{\tau_k}\cdot \w^\sharp$ forms a weakly decreasing sequence for ${\bf m}$, which implies that $\Delta_{\mf x}(U)\leq L_{\mf x}({\bf m})$.

Conversely, let $\w_1\cdot\w_2$ be a maximal weakly decreasing sequence for ${\bf m}$. Since $A$ is a diagonal matrix,  $\w_1=k^{\tau_k}$ for some $k$. Let $\w$ be a subword of $w(U)$ such that $\w^\sharp=\w_2$. Then by Remark \ref{remark on L} (2), we may also assume that $\w=w_1\ldots w_r$ is weakly decreasing and $w_1$ is in the $k$-th row of $\tau$. Since $L_{\mf x}({\bf m})=(\tau_k+2\ell(\w_2))/\epsilon=(\tau_k+2\ell(\w))/\epsilon$, we have $L_{\mf x}({\bf m})\leq \Delta_{\mf x}(U)$.
\qed

By Lemma \ref{L=Delta}, the map ${\bf m}\otimes t_{n{\Lambda}^{\mf x}_{0}} \mapsto \kappa^{\mf x}_\lambda(A)$ gives  an
${\mf x}_\infty$-crystal isomorphism
$$\M^{\mf x}(\lambda,n) \longrightarrow \mathcal{T}^{\mf x}(\lambda,n).$$
Therefore by Theorem \ref{main-1}, we have
\begin{equation*}
\mathcal{T}^{\mf x}(\lambda,n) \simeq \B({\mf x}_\infty,\Lambda^{\mf x}(\lambda,n)).
\end{equation*}
This proves Theorem \ref{main result}.\qed

\end{document}